\documentclass[12pt]{article}

\usepackage{amsthm}
\usepackage[dvips]{graphicx}

\usepackage{epsfig}

\usepackage{latexsym}
\usepackage{amsmath}
\usepackage{amssymb}

\theoremstyle{definition}

\newtheorem{dfn}{Definition}[section]

\newtheorem{ex}[dfn]{Example}

\newtheorem{rem}[dfn]{Remark}

\theoremstyle{plain}

\newtheorem{add}[dfn]{Addendum}

\newtheorem{ass}[dfn]{Assumption}
\newtheorem{cor}[dfn]{Corollary}

\newtheorem{lem}[dfn]{Lemma}
\newtheorem{prop}[dfn]{Proposition}

\newtheorem{sublem}[dfn]{Sublemma}
\newtheorem{thm}[dfn]{Theorem}

\newcommand{\co}{\colon\thinspace}    

\def\lra{\longrightarrow}

\def\ra{\rightarrow}

\newenvironment{ftheorem}{\noindent{\bf Fibration theorem}\em}{}

\newenvironment{sphtheorem}{\noindent{\bf Spherical uniformization}\em}{}

\newenvironment{maintheorem}{\noindent{\bf Main Theorem (Uniformization of
small 3-orbifolds).}\em}{}

\def\JJ{{\mathcal J}}
\def\II{{\mathcal I}}
\def\TT{{\mathcal T}}

\def\NN{{\mathcal N}}
\def\CC{{\mathcal C}}
\def\OO{{\mathcal O}}

\def\XX{{\mathcal X}}

\def\L{{\Lambda}}

\def\R{{\mathbb R}}
\def\H{{\mathbb H}}
\def\Z{{\mathbb Z}}

\def\N{{\mathbb N}}
\def\M{{\mathbb M}}
\def\S{{\mathbb S}}

\def\Diam{\operatorname{diam}}
\def\Hom{\operatorname{Hom}}

\def\Int{\operatorname{int}}
\def\trace{\operatorname{trace}}
\def\V{\operatorname{v}}
\def\Vol{\operatorname{vol}}

\def\Star#1{\operatorname{star}({#1})}

\def\al{\alpha}
\def\de{\delta}
\def\De{\Delta}
\def\eps{\varepsilon}
\def\la{\lambda}
\def\La{\Lambda}
\def\om{\omega}

\def\si{\sigma}
\def\Si{\Sigma}

\def\ga{\gamma}

\def\Ga{\Gamma}

\def\D{\partial}
\def\half{\frac{1}{2}}
\def\ol{\overline}

\def\3{\ss}
\def\8{\infty}

\def\BI{\begin{itemize}}
\def\EI{\end{itemize}}

\begin{document}

\title{Geometrization of $3$-dimensional orbifolds}
\author{Michel Boileau, Bernhard Leeb, Joan Porti}
\date{November 28, 2003}

\maketitle

\begin{abstract}
This paper is devoted to the proof of the orbifold theorem: If $\OO$ is a
compact connected orientable irreducible and topologically atoroidal
3-orbifold with non-empty ramification locus, then $\OO$ is geometric
(i.e.\ has a metric of constant curvature or is Seifert fibred). As a
corollary, any smooth orientation preserving non-free finite group action
on $S^3$ is conjugate to an orthogonal action.
\end{abstract}

\tableofcontents

\section{Introduction}

A 3-dimensional orbifold is a metrizable space equipped with an atlas of
compatible local models given by quotients of $\mathbb R^3$ by finite
subgroups of $O(3)$. For example, the quotient of a 3-manifold by a
properly discontinuous smooth group action naturally inherits a structure
of a 3-orbifold. When the group action is finite, such an orbifold is
called very good. We will consider in this paper only orientable
orbifolds. The ramification locus, i.e.\ the set of points with
non-trivial local isotropy group, is then a trivalent graph.

In 1982, Thurston \cite{ThuBull,Thu6} announced the
geometrization theorem for 3-orbifolds with non-empty ramification
locus and lectured about it.
Several partial results have been obtained in the meantime,
see \cite{BoP}.
The purpose of this article is to give a
complete proof of the orbifold theorem,
compare \cite{cras} for an outline.
A different proof has been announced in \cite{CHK}.

The main result of this article is the following uniformization
theorem which implies the orbifold theorem for compact orientable
3-orbifolds. A 3-orbifold $\OO$ is said to be \emph{geometric} if
either its interior has one of Thurston's eight geometries or
$\OO$ is the quotient of a ball by a finite orthogonal action.

\medskip
\begin{maintheorem}
\label{thm:unismall} Let $\OO$ be a compact connected orientable small
$3$-orbifold with non-empty ramification locus. Then $\OO$ is geometric.
\end{maintheorem}
\medskip

An orientable compact $3$-orbifold $\OO$ is \emph{small} if it is
irreducible,  its boundary $\partial \OO$ is a (perhaps empty)
collection of turnovers (i.e. 2-spheres with three branching points),
and  it does not contain any other closed incompressible orientable
$2$-suborbifold.

An application of the main theorem concerns non-free finite group
actions on the 3-sphere $S^3$, see Section \ref{sec:redgroup}. It
recovers all the previously known partial results (cf. \cite{DaM},
\cite{Fei}, \cite{MB}, \cite{Mor}), as well as the results about
finite group actions on the 3-ball   (cf. \cite{MYTwo},
\cite{KS}).

\begin{cor}
\label{cor:finact}
An orientation preserving smooth non-free
finite group action on $S^3$ is smoothly conjugate to an orthogonal action.
\end{cor}

Every compact orientable irreducible and atoroidal 3-orbifold can
be canonically split along a maximal (perhaps empty) collection of
disjoint and pairwise non-parallel hyperbolic turnovers. The
resulting pieces are either \emph{Haken} or \emph{small}
3-orbifolds (cf. Section \ref{sec:reduction}).
Using an extension of Thurston's hyperbolization theorem to the
case of Haken orbifolds (cf. \cite[Ch. 8]{BoP}), we show that the
main theorem
implies the orientable case of the orbifold theorem:

\begin{cor}[Orbifold Theorem]\label{cor:TOT} Let $\OO$ be a
compact connected orientable irreducible $3$-orbifold with
non-empty ramification locus. If $\OO$ is topologically atoroidal,
then $\OO$ is geometric.
\end{cor}

Any compact connected orientable  $3$-orbifold, that does not contain any
bad 2-suborbifold (i.e. a 2-sphere with one branching point or with two
branching points having different branching index), can be split
 along a finite collection of disjoint embedded spherical and
toric 2-suborbifolds  (\cite[Ch. 3]{BMP})  into irreducible and
atoroidal 3-orbifolds, which are geometric by Corollary \ref{cor:TOT}.
The fact that 3-orbifolds with a geometric decomposition are finitely
covered by a manifold \cite{McCMi} implies:

\begin{cor}
Every compact connected orientable $3$-orbifold which does not contain any
bad 2-suborbifold is the quotient of a compact orientable $3$-manifold by
a finite group action.
\end{cor}

The paper is organized as follows. In Section~\ref{sec:reduction}
we recall some basic terminology about orbifolds. Then we deduce
the orbifold theorem from our main theorem.

The proof of the main theorem is based on some geometric
properties of cone manifolds, which are presented in
Sections~\ref{sec:conemanifolds}-\ref{sec:conelowdiambound}. This
geometric approach is one of the main differences with \cite{BoP}.

In Section \ref{sec:conemanifolds}, we define cone manifolds and
develop some basic geometric concepts. Motivating examples are
geometric orbifolds which arise as quotients of model spaces by
properly discontinuous group actions. These have cone angles
$\leq\pi$, and only cone manifolds with cone angles $\leq\pi$ will
be relevant for the approach to geometrizing orbifolds pursued in
this paper. The main result of Section \ref{sec:conemanifolds} is
a compactness result for spaces of cone manifolds with cone angles
$\leq\pi$ which are thick in a certain sense.

In Section \ref{sec:euc} we classify noncompact Euclidean cone
3-manifolds with cone angles $\leq \pi$. This classification is
needed for the proof of the fibration theorem in Section
\ref{sec:Seifert}. It also motivates our results in Section
\ref{sec:conelowdiambound} where we study the local geometry of
cone 3-manifolds with cone angles $\leq\pi$; there, a lower
diameter bound plays the role of the noncompactness condition in
the flat case. Our main result, cf.\ Section \ref{sec:thin}, is a
geometric description of the thin part in the case when cone angles are
bounded away from $\pi$ and $0$ (Theorem \ref{descthin}). As
consequences, we obtain thickness (Theorem \ref{thick}) and, when
the volume is finite, the existence of a goemetric compact core
(Theorem \ref{finiteness}). The other results relevant for the
proof of the main theorem are the geometric fibration theorem for
thin cone manifolds with totally geodesic boundary (Corollary
\ref{cor:Ibundle}) and the thick vertex lemma
(Lemma~\ref{lem:thickvertex}) which is a simple result useful in
the case of platonic vertices.

We give the proof of the main theorem in Section~\ref{sec:main}.
Firstly we reduce to the case when the smooth part of the orbifold
is hyperbolic. We view the (complete) hyperbolic structure on the
smooth part  as a hyperbolic cone structure on the
 orbifold with cone angles zero.
The goal is to increase the cone angles of this hyperbolic cone
structure as much as possible. In Section~\ref{sec:deform} we
prove first that there exist such deformations which change the
cone angles (\emph{openness theorem}).

Next we consider a sequence of hyperbolic cone structures on the orbifold
whose cone angles   converge to the supremum of the cone angles
in the deformation space. We have the following dichotomy: either the
sequence collapses (i.e. the supremum of the injectivity radius for
each cone structure goes to zero) or not (i.e. each cone structure
contains a point with injectivity radius uniformly bounded away from
zero.

In the \emph{non-collapsing} case we show in Section~\ref{sec:deg}  that the
orbifold angles can be reached in the deformation space of
hyperbolic cone structures, and therefore the orbifold is
hyperbolic. This step uses a \emph{stability
theorem} which shows that a non-collapsing sequences of hyperbolic cone
structures on the orbifold has a subsequence converging to a
hyperbolic cone structure on the orbifold. We prove this theorem in Section~\ref{sec:noncollapse}.

Then we analyze the case where the sequence of cone structures collapses.
If the diameters of the collapsing cone structures
are bounded away from zero, then we conclude that the orbifold is Seifert
fibred, using the   \emph{fibration theorem} which is proved in
Section \ref{sec:Seifert}.
Otherwise the diameter of the sequence of cone structures  converges to zero.
Then  we show that the orbifold is geometric, unless the
following
situation
occurs:
the orbifold is closed and admits a Euclidean
cone structure with cone angles strictly less than its orbifold angles.

We deal with this last case in Sections~\ref{sec:sphericaluniformization}
and \ref{sec:sphericalcone} prove that then the orbifold is spherical
(\emph{spherical uniformization} theorem). For orbifolds with cyclic or
dihedral stabilizer, the proof relies on Hamilton's theorem \cite{HamOne}
about the Ricci flow on 3-manifolds. In the
general case the proof is by induction on the number of platonic vertices
and involves deformations of spherical cone structures.

\medskip
{\bf Acknowledgements.} We wish to thank D.\ Cooper for useful
conversations and J.\ Alze for reading an earlier version of our manuscript.
We thank the RiP-program at the Mathematisches
Forschungsinstitut Oberwolfach, as well as DAAD, MCYT (Grants
HA2000-0053 and BFM2000-0007) and DURSI (ACI2000-17) for financial
support.


%
%
%
%
%
%
%

\section{3-dimensional orbifolds}
\label{sec:reduction}

\subsection{Basic definitions}

For a general background about orbifolds we refer to \cite{BMP},
\cite{BSOne,BSTwo}, \cite{CHK}, \cite{DaM}, \cite[Ch.\ 7]{Kap},
\cite{Scott}, and \cite[Ch.\ 13]{ThuNotes}. We begin with
recalling some terminology from there.

A compact 2-orbifold $F^2$ is said to be \emph{spherical},
\emph{discal}, \emph{toric} or \emph{annular} if it is the
quotient by a finite smooth group action of respectively the
2-sphere $S^2$, the 2-disc $D^2$, the 2-torus $T^2$ or the annulus
$S^1\times [0,1]$.

A compact 2-orbifold is \emph{bad} if it is not good (i.e. it is
not covered by a surface). Such a 2-orbifold is the union of two
non-isomorphic discal 2-orbifolds along their boundaries.

A compact 3-orbifold $\OO$ is \emph{irreducible} if it does not
contain any bad 2-suborbifold and if every orientable spherical
2-suborbifold bounds in $\OO$ a discal 3-suborbifold, where a
\emph{discal} 3-orbifold is a finite quotient of the 3-ball by an
orthogonal action.

A connected 2-suborbifold $F^2$ in an orientable 3-orbifold $\OO$ is
\emph{compressible} if either $F^2$ bounds a discal 3-suborbifold in
$\OO$ or there is a discal 2-suborbifold $\Delta^2$ which intersects
transversally $F^2$ in $\partial\Delta^2=\Delta^2\cap F^2$ and such
that $\partial\Delta^2$ does not bound a discal 2-suborbifold in
$F^2$.

A 2-suborbifold $F^2$ in an orientable 3-orbifold $\OO$ is
\emph{incompressible} if no connected component of $F^2$ is
compressible in $\OO$.

A properly embedded 2-suborbifold $F^2$  is
\emph{$\partial$-parallel} if it co-bounds a product with a
suborbifold of the boundary (i.e.\ an embedded product $\overline
F\times [0,1]\subset\mathcal O$ with $\overline F\times 0=F^2$ and
$\overline F\times 1\subset
\partial\mathcal O$), so that $\partial \overline F\times
[0,1]\subset
\partial\mathcal O $.

A properly embedded 2-suborbifold $(F,\partial F) \hookrightarrow
(\OO, \partial \OO)$ is \emph{$\partial$-compressible} if:
\begin{itemize}
 \item [--] either $(F,\partial F)$ is a discal 2-suborbifold
$(D^2,\partial D^2)$ which is $\partial$-parallel,
 \item [--] or there is a discal 2-suborbifold $\Delta \subset \OO$
 such that $\partial \Delta \cap F$ is a simple arc $\alpha$ which does not
cobound a discal suborbifold of $F$ with an arc in $\partial F$, and
$\Delta
 \cap \partial \OO$ is a simple arc $\beta$ with $\partial \Delta =
 \alpha \cup \beta$ and $\alpha \cap \beta =
 \partial \alpha = \partial \beta$.
 \end{itemize}

A properly embedded $2$-suborbifold $F^2$ is \emph{essential} in a
compact orientable irreducible $3$-orbifold, if it is incompressible,
$\partial$-incompressible and not $\partial$-parallel.

A compact 3-orbifold is \emph{topologically atoroidal} if it does not
contain any embedded essential orientable toric $2$-suborbifold.

A \emph{turnover} is a $2$-orbifold with underlying space a
$2$-sphere and ramification locus three points. In an irreducible
orientable 3-orbifold, an embedded turnover either bounds a discal
$3$-suborbifold or is incompressible and of non-positive Euler
characteristic.

An orientable compact $3$-orbifold $\OO$ is \emph{Haken} if it is
irreducible, if every embedded turnover is either compressible or
$\partial$-parallel, and if it contains an embedded orientable
incompressible 2-suborbifold which is not a turnover.

\begin{rem} The word Haken may lead to confusion, since it is
not true that a compact orientable irreducible $3$-orbifold
containing an orientable incompressible properly embedded
$2$-suborbifold is Haken in our meaning (cf. \cite[Ch. 4]{BMP},
\cite{DunOne}, \cite[Ch. 8]{BoP}).
\end{rem}

An orientable compact $3$-orbifold $\OO$ is \emph{small} if it is
irreducible, its boundary $\partial \OO$ is a (perhaps empty)
collection of turnovers, and $\OO$ does not contain any essential
orientable $2$-suborbifold.  It follows from Dunbar's theorem
\cite{DunOne} that the hypothesis about the boundary is
automatically satisfied once $\mathcal O$ does not contain any
essential 2-suborbifold.

\begin{rem}
By irreducibility, if a small orbifold $\OO$ has  non-empty boundary,
then either $\OO$ is a discal $3$-orbifold, or $\partial \OO$ is a
collection of Euclidean and hyperbolic turnovers.
\end{rem}

A 3-orbifold $\OO$ is \emph{geometric} if either is the quotient
of a ball by an orthogonal action, or its interior has one of the
eight Thurston's geometries. We quickly review those geometries.

A compact orientable 3-orbifold $\OO$ is \emph{hyperbolic} if its
interior is orbifold-diffeo\-morphic to the quotient of the
hyperbolic space $\mathbb H^3$ by a non-elementary discrete group of
isometries. In particular $I$-bundles over hyperbolic 2-orbifolds are
hyperbolic, since their interiors are quotients of $\mathbb H^3$ by
non-elementary Fuchsian groups.

A compact orientable 3-orbifold is \emph{Euclidean} if its interior
has a complete Euclidean structure. Thus, if a compact orientable and
$\partial$-incompressible 3-orbifold $\OO$ is Euclidean, then either
$\OO$ is a $I$-bundle over a 2-dimensional Euclidean closed orbifold
or $\OO$ is closed.

A compact orientable 3-orbifold is \emph{spherical} when it is the
quotient of the standard sphere $\mathbb S^3$ or the round ball $
B^3$ by the orthogonal action of a finite group.

A \emph{Seifert fibration} on a 3-orbifold $\OO$ is a partition of
$\OO$ into closed 1-suborbifolds (circles or intervals with silvered
boundary) called fibers, such that each fiber has a saturated
neighborhood diffeomorphic to $S^1\times D^2/G$, where $G$ is a
finite group which acts smoothly, preserves both factors, and acts
orthogonally on each factor and effectively on $D^2$; moreover the
fibers of the saturated neighborhood correspond to the quotients of
the circles $S^1\times\{*\}$. On the boundary $\partial\OO$, the
local model of the Seifert fibration is $S^1\times D^2_+/G$, where
$D^2_+$ is a half disc.

A 3-orbifold that admits a Seifert fibration is called Seifert fibred. A
Seifert fibred 3-orbifold which does not contain a bad 2-suborbifold is
geometric (cf. \cite[Chap. 1, 2]{BMP}, \cite{Scott}, \cite{Thu7}).

Besides the constant curvature geometries $\mathbb E^3$ and $\mathbb S^3$,
there are four other possible 3-dimensional homogeneous geometries for a
Seifert fibred 3-orbifold: $\mathbb H^2\times\mathbb R$, $\mathbb
S^2\times \mathbb R$, $\widetilde{SL_2(\mathbb R)}$ and $Nil$.

The  geometric but non Seifert fibred 3-orbifolds require either a
constant curvature geometry or $Sol$. Compact 3-orbifolds with $Sol$
geometry are fibred over a closed 1-dimensional orbifold with toric fiber
and thus they are not topologically atoroidal (cf. \cite{DunTwo}).
\bigskip

\subsection{Spherical and toric decompositions}
\label{sec:decomp}

Thurston's geometrization conjecture asserts that any compact,
orientable, 3-orbifold, which does not contain any  bad
$2$-suborbifold, can be decomposed along a finite collection of
disjoint, non-parallel, essential, embedded spherical and toric
2-suborbifolds into geometric suborbifolds.

The topological background for Thurston's geometrization
conjecture is given by the spherical and toric decompositions.

Given a compact  orientable $3$-orbifold  without bad
$2$-suborbifolds, the first stage of the splitting is called
\emph{spherical} or \emph{prime}  decomposition, and it expresses
the $3$-orbifold as the connected sum of $3$-orbifolds which are
either homeomorphic to a finite quotient of $S^1\times S^2$ or
irreducible. We refer to \cite[Chap. 3]{BMP}, \cite{ty:decomp} for
details.

The second stage (\emph{toric} splitting) is a more subtle
decomposition of each irreducible factor along a finite (maybe
empty) collection of disjoint and non-parallel essential, toric
$2$-suborbifold. This collection  of essential toric
$2$-suborbifolds is unique up to isotopy. It cuts the irreducible
orbifold into topologically atoroidal or Seifert fibered
$3$-suborbifolds; see \cite{BSOne}, \cite[Chap. 3]{BMP}.

By these spherical and toric decompositions, Thurston's
geometrization conjecture reduces to the case of a compact,
orientable $3$-orbifold which is irreducible and topologically
atoroidal.

Our proof requires a further decomposition along turnovers due to
Dunbar (\cite[Chap 3]{BMP}, \cite[Thm. 12]{DunOne}).
A compact irreducible and topologically atoroidal 3-orbifold
has a maximal family of non-parallel essential turnovers, which
may be empty. This family is unique up to isotopy and
cuts the orbifold into pieces without essential turnovers.

\bigskip

\subsection{Finite group actions on spheres with fixed points}
\label{sec:redgroup}

\begin{proof}[Proof of Corollary~\ref{cor:finact} from the main theorem.]
Consider a non-free action of a finite group $\Gamma$ on $S^3$
by orientation preserving diffeomorphisms.
Let $\OO=\Gamma\backslash S^3$ be the quotient orbifold.

If $\OO$ is irreducible then the equivariant Dehn lemma implies
that any 2-suborbifold with infinite fundamental group has a
compression disc. Hence $\OO$ is small and we apply the main
theorem.

Suppose that $\OO$ is reducible. Since $\OO$ does not contain a
bad $2$-suborbifold, there is a prime decomposition along a family
of spherical 2-suborbifolds, see Section~\ref{sec:decomp}. These
lift to a family of 2-spheres in $S^3$. Consider an innermost
2-sphere; it bounds a ball $B\subset S^3$. The quotient $Q$ of $B$
by its stabilizer $\Gamma'$ in $\Gamma$ has one boundary component
which is a spherical 2-orbifold. We close it by attaching a discal
3-orbifold. The resulting closed 3-orbifold $\OO'$ is a prime
factor of $\OO$.  The orbifold $\OO'$ is irreducible, and hence
spherical. The action of $\Gamma'$ on $\widetilde{\OO'}\cong S^3$
is standard and preserves the sphere $\D B$. Thus the action is a
suspension and $Q$ is discal. This contradicts the minimality of
the prime decomposition.
\end{proof}

\subsection{Proof of the orbifold theorem from the main theorem}

This step of the proof is based on the following extension of Thurston's
hyperbolization theorem to Haken orbifolds (cf. \cite[Ch. 8]{BoP}):

\begin{thm}[Hyperbolization theorem of Haken orbifolds]\label{thm:Haken}
Let $\OO$ be a compact orientable connected Haken $3$-orbifold. If $\OO$
is topologically atoroidal and not Seifert fibred, nor Euclidean, then
$\OO$ is hyperbolic.
\end{thm}

\begin{rem} The proof of this theorem follows exactly the scheme of the
proof for Haken manifolds \cite{ThuBull,Thu3,Thu5}, \cite{McMOne},
\cite{Kap}, \cite{OtaOne,OtaTwo} (cf. \cite[Ch. 8]{BoP} for a precise
overview).
\end{rem}

\begin{proof}[Proof of Corollary~\ref{cor:TOT} (the orbifold theorem).]
Let $\OO$ be a compact orientable connected irreducible topologically
atoroidal $3$-orbifold. By \cite[Chap 3]{BMP}, \cite[Thm. 12]{DunOne}
there exists in $\OO$ a (possibly empty) maximal collection $\TT$ of
disjoint embedded pairwise non-parallel essential turnovers. Since $\OO$
is irreducible and topologically atoroidal, any turnover in $\TT$ is
hyperbolic (i.e. has negative Euler characteristic).

When $\TT$ is empty, the orbifold theorem reduces either to the
main theorem
or to Theorem~\ref{thm:Haken} according to whether $\OO$ is small or
Haken.

When $\TT$ is not empty, we first cut open the orbifold $\OO$ along
the turnovers of the family $\TT$. By maximality of the family $\TT$,
the closure of each component of $\OO-\TT$ is a compact orientable
irreducible topologically atoroidal $3$-orbifold that does not
contain any essential embedded turnover.
 Let $\OO'$ be one of these
connected components. By the previous case $\OO'$ is either hyperbolic,
Euclidean or Seifert fibred. Since, by construction, $\partial\OO'$
contains at least one hyperbolic turnover $T$, $\OO'$ must be hyperbolic.
Moreover any such hyperbolic turnover $T$ in $\partial\OO'$ is a Fuchsian
$2$-suborbifold, because there is a unique conjugacy class of faithful
representations of the fundamental group of a turnover in $PSL_2(\mathbb
C)$.

We assume first that all the connected components of
$\OO-{\mathcal T}$ have $3$-dimensional convex cores. In this case
the totally geodesic hyperbolic turnovers are the boundary
components of the convex cores. Hence the hyperbolic structures on
the components of $\OO-{\TT}$ can be glued together along the
hyperbolic turnovers of the family ${\mathcal T}$ to give a
hyperbolic structure on the $3$-orbifold $\OO$.

If the convex core of $\OO'$ is 2-dimensional, then $\OO'$ is either a
product $T\times [0,1]$, where $T$ is a hyperbolic turnover, or a quotient
of $T\times [0,1]$ by an involution. When $\OO'=T\times [0,1]$, then the
$3$-orbifold $\OO$ is Seifert fibred, because the mapping class group of a
turnover is finite. When $\OO'$ is the quotient of $T\times [0,1]$, then
it has only one boundary component and it is glued either to another
quotient of $T\times [0,1]$ or to a component with $3$-dimensional convex
core. When we glue two quotients of $T\times [0,1]$ by an involution, we
obtain a Seifert fibred orbifold. Finally, gluing $\OO'$ to a hyperbolic
orbifold with totally geodesic boundary is equivalent to quotient this
boundary by an isometric involution.
\end{proof}

\section {3-dimensional cone manifolds}
\label{sec:conemanifolds}

\subsection{Basic definitions}

We start by recalling the construction of metric cones.

Let $k$ and $r>0$ be real numbers;
if $k>0$ we assume in addition that $r\leq\frac{\pi}{\sqrt{k}}$.
Suppose that $Y$ is a metric space with $diam(Y)\leq\pi$.
On the set $Y\times [0,r]$ we define a pseudo-metric as follows.
Given $(y_1,t_1),(y_2,t_2)\in Y\times [0,r]$, let $p_0p_1p_2$ be a
triangle in the 2-dimensional model space $\M^2_k$ of constant curvature $k$
with $d(p_0,p_1)=t_1$, $d(p_0,p_2)=t_2$ and
$\angle_{p_0}=d_Y(y_1,y_2)$.
We put
\[ d_{Y\times [0,r]}\bigl((y_1,t_1),(y_2,t_2)\bigr):=d_{\M^2_k}(p_1,p_2) \]
The metric space $C_{k,r}(Y)$ obtained from collapsing the subset
$Y\times\{0\}$ to a point is called the
{\em metric cone of curvature k} or {\em $k$-cone} of radius $r$ over $Y$.
In the special case when $k>0$ and $r=\frac{\pi}{\sqrt{k}}$,
one also has to collapse the subset $Y\times\{\frac{\pi}{\sqrt{k}}\}$ to a
point.
The point in $C_{k,r}(Y)$
corresponding to $Y\times\{0\}$ is called the {\em tip} or {\em apex}
of the cone.
The {\em complete} $k$-cone or simply
{\em $k$-cone} $C_k(Y)$ over $Y$ is defined as
$C_{k,\infty}(Y):=\cup_{r>0}C_{k,r}(Y)$ if $k\leq0$
and as
$C_{k,\frac{\pi}{\sqrt{k}}}(Y)$ if $k>0$.
The complete 1-cone over a space is also called its {\em metric
suspension}.

\medskip
We define cone manifolds as certain metric spaces locally isometric to
iterated cones.
To make this precise, we proceed by induction over the dimension.
We first make the convention that the connected
$1$-dimensional cone manifolds of curvature $1$ are circles
of length $\leq2\pi$ or compact intervals of length $\leq\pi$.

\begin{dfn}[Cone manifolds]
A $n$-dimensional {\em conifold}
of curvature $k$, $n\geq2$, is a complete
geodesic metric space locally isometric to the $k$-cone over a
connected $(n-1)$-dimensional conifold of curvature $1$.

A {\em cone manifold} is a conifold which is topologically a manifold.
\end{dfn}

Conifolds of curvature $k=+1$, $k=0$ or $k=-1$ are called
{\em spherical, Euclidean} or {\em hyperbolic}, respectively.

Spelled out in more detail, the definition requires
that for every point $x$ in a $n$-conifold $X$
there exists a radius $\eps>0$ and an isometry from the closed ball
$\ol B_{\eps}(x)$ to the $k$-cone $C_{k,\eps}(\La_xX)$ over a metric space
$\La_xX$ carrying $x$ to the tip of the cone.
Moreover,
$\La_xX$ must be itself a $(n-1)$-conifold of curvature $1$.

The metric space $\La_xX$ is called the {\em space of directions} or
{\em link} of $X$ at $x$.\footnote{
The standard geometric notation would be $\Si_xX$, but we already make
extensive use of the letter $\Si$, namely for the singular locus of an
orbifold.}
It can be defined intrinsically as the space of germs of
geodesic segments in $X$ emanating from $x$ equipped with the angular
metric.
It is implicit in the definition that the links $\La_xX$ are complete
metric spaces.
Since they have curvature $1$,
it follows that they are compact with diameters $\leq\pi$,
see the discussion at the end of this section.

We note that all conifolds of dimension $\leq2$ are manifolds.
The links in 3-dimensional conifolds are,
according to the Gau\3-Bonnet Theorem (extended to singular surfaces),
topologically 2-spheres, 2-discs or projective planes.
If none of the links is a projective plane, then the conifold is a manifold.
The wider concept of conifold will play no role in this paper;
later on, we will only consider cone {\em manifolds} of dimensions $\leq3$.

\begin{ex}[Geometric orbifolds]
A {\em geometric orbifold} of dimension $n$ and curvature $k$
is a complete geodesic metric space
which is locally isometric to the quotient of the model space
$\M_{k}^{n}$
by a finite group of isometries.
\end{ex}

Unlike topological orbifolds,
geometric orbifolds are always {\em global} quotients,
i.e.\ they are (even finite) quotients of manifolds of constant curvature
by discrete group actions.

\medskip
We define the {\em boundary} of a conifold by induction over the
dimension.
The boundary points of a 1-conifold are the endpoints of its interval
components.
The boundary points of a $n$-conifold, $n\geq2$, are the points whose
links have boundary.

A point $x$ in a conifold $X$ is called a {\em smooth interior} point
if $X$ is locally at $x$ isometric to the model space $\M_k^n$ of the same
curvature and dimension as $X$,
or equivalently,
if the link $\La_xX$ is a unit sphere.
If $\La_xX$ is a unit hemisphere, the point $x$ is a {\em smooth
boundary} point.
All other points are called {\em singular}.
We denote by $X^{smooth}$ the subset of smooth points,
and by $\Si_{X}$ its complement, the {\em singular locus}.

Let us go through this in low dimensions.
One-dimensional cone manifolds contain only regular points.
If $S$ is a {\em cone surface}, i.e.\ a cone $2$-manifold,
then $\Si_S$ is a discrete subset.
A singular point is either a {\em corner} of the boundary, if its link is an
interval of length $<\pi$,
or a {\em cone point} in the interior, if its link is a circle of
length $<2\pi$.
In the latter case, the length of the circle is called the
{\em cone angle}.

Consider now a $3$-dimensional cone manifold $X$.
In this case, the singular set $\Si_{X}$ is one-dimensional,
namely a geodesic graph.
We define $\Si_{X}^{(1)}\subseteq \Si_{X}$ as the subset of singular
points $x$ whose link $\La_xX$ is the metric suspension of (complete
$1$-cone over) a circle.
The length of the circle is called the cone angle at $x$.
We call the closure of a component of $\Si_{X}^{(1)}$ a {\em singular
edge}.
The cone angle is constant along edges,
and we can thus speak of the {\em cone angle} of an edge.
The complement $\Si_{X}^{(0)}:=\Si_{X}-\Si_{X}^{(1)}$
is discrete and its points are called {\em singular vertices}.

Notice that a cone surface or a cone 3-manifold without boundary
is a geometric orbifold
iff all cone angles are divisors of $2\pi$.
In particular
the cone angles of a geometric 3-orbifold are $\leq\pi$, and due to
this fact we will be mostly interested in cone manifolds {\em with
cone angles $\leq\pi$}.

\begin{prop}
Conifolds of curvature $k$
are metric spaces with curvature $\geq k$ in the sense of Alexandrov.
\end{prop}

This can be readily seen by induction over the dimension using the
following facts:
Since conifolds are metrically complete by assumption, a local
curvature bound implies a global curvature bound (Toponogov's Theorem);
the $k$-cones over compact intervals of length $\leq\pi$ and circles
of length $\leq2\pi$ are spaces with curvature $\geq k$;
the $k$-cone over a space with curvature $\geq1$ is a space with
curvature $\geq k$.
Note also that spaces with curvature $\geq1$ have diameter $\leq\pi$,
due to the singular version of the Bonnet-Myers theorem, cf.\
\cite[Thm. 3.6]{BGP}.


All our geometric considerations will take place within the framework
of metric spaces with curvature bounded below.
For this theory, we refer the reader to the fundamental paper
\cite{BGP} and the introductory text \cite[Ch. 10]{burago}.

\subsection{Exponential map, cut locus, injectivity radius}
\label{sec:expcutinj}

Consider a connected conifold $X$ of curvature $k$ and dimension $\geq2$.

For a point $p\in X$,
according to our requirement on the local geometry of conifolds,
there exists $\eps>0$ such that
the cone $C_{k,\eps}(\La_pX)$ canonically embeds into $X$,
its tip $O$ being mapped to $p$.
This embedding extends naturally to a map from a larger domain inside the
complete cone $C_k(\La_pX)$ as follows:
Let ${\cal E}(p)\subseteq C_k(\La_pX)$ be the union of all geodesic
segments $\ol{Oy}$, such that there exists a geodesic segment $\ol{px_y}$ in
$X$ with the same initial direction
modulo the natural identification $\La_O(C_k(\La_pX))\cong\La_pX$.
The subset ${\cal E}(p)$ is star-shaped with respect to $O$,
and we define the {\em exponential map} in $p$
\[ \exp_p:{\cal E}(p)\lra X \]
as the map sending each point $y$ to the respective point $x_y$.

The {\em conjugate radius} is defined, purely in terms of the curvature,
as $r_{conj}:=\frac{\pi}{\sqrt{k}}$ if $k>0$ and
$r_{conj}:=\infty$ if $k\leq0$,
i.e.\ $r_{conj}=diam(C_k(\La_pX))$.
The {\em geodesic radius} in a point $p$, $0<r_{geod}(p)\leq r_{conj}$,
is the radius of the largest ball in $C_k(\La_pX)$ around $O$ on which
$\exp_p$ is defined.

Let $x$ be an interior point of a geodesic segment $\si=\ol{pq}$.
Then $\La_xX$ has extremal diameter $\pi$ and, by the Diameter Rigidity
Theorem, is a metric suspension with the directions of $\si$ in $x$ as poles.
The equator of the suspension consists of the directions at $x$ perpendicular
to $\si$.

For any $0<d<\min(d(p,q),r_{conj})$ exists a sufficiently small
$\de>0$ such that the ``thin'' cone $C_{k,d}(B_{\de}(\La_p\si))$
is contained in ${\cal E}(p)$ and embeds via $\exp_p$ locally
isometrically into $X$. Here $\La_p\si\in\La_pX$ denotes the
direction of $\si$ at its endpoint $p$.

If $\si$ has length $<r_{conj}$,
and if $\si'=\ol{pq'}$ is sufficiently Hausdorff close to $\si$,
then there exists an isometrically immersed (2-dimensional) triangle of
constant curvature $k$ with $\si$ and $\si'$ as two of its sides.
It follows that
there don't exist other geodesic segments with the same endpoints as
$\si$ and arbitrarily Hausdorff close to $\si$.

We now focus our attention on {\em minimizing} geodesic segments.
Let $p$ and $q$ be points with $d(p,q)<r_{conj}$.
Our discussion implies that
there are at most finitely many minimizing geodesic segments
$\si_1,\dots,\si_m$ connecting them.

If $x$ is a point sufficiently close to $q$,
then for every $i$ there exists a locally isometrically embedded triangle
$\De_i$ with $x$ as vertex and $\si_i$ as opposite side.
Moreover,
any minimizing segment $\tau=\ol{px}$ is Hausdorff close to one of the
segments $\si_i$ and coincides with the side $\ol{px}$ of the corresponding
triangle $\De_i$.
So, there exists a minimizing segment $\ol{px}$ Hausdorff close to $\si_j$ iff
$\angle_q(\si_j,x)=\min_i\angle_q(\si_i,x)$.

Let ${\cal D}(p)\subseteq{\cal E}(p)$ be the union of all geodesic
segments $\ol{Oy}$ in $\La_pX$ whose images $\ol{px_y}$ under
$\exp_p$ are {\em minimizing} segments. Let $\dot {\cal
D}(p)\subseteq{\cal D}(p)$ be the subset consisting of $O$ and all
interior points of such segments $\ol{Oy}$. Note that $\dot {\cal
D}(p)$ is open
and its closure equals ${\cal D}(p)$.
We have ${\cal D}(p)-\dot {\cal D}(p)=\D {\cal D}(p)$
except in the special case when $k>0$ and $X$ is a metric suspension with tip
$p$.

\begin{dfn}
[Cut locus] The subset $Cut_X(p)=Cut(p):=\exp_p({\cal D}(p)-\dot
{\cal D}(p))\subset X$ is called the {\em cut locus} with respect
to the point $p$.
\end{dfn}

In other words,
$Cut(p)$ is the complement of the union of $p$ and all
minimizing half-open segments $\ga:[0,l)\ra X$ with initial point $\ga(0)=p$.

In the same way one defines the cut locus $Cut(F)$ with respect to a
{\em finite set} $F\subset X$.
Our above discussion implies

\begin{prop}
\label{loccon}
The cut locus is {\em locally conical},
i.e.\ for any point $q$ with $d(p,q)<r_{conj}$ exists $\eps>0$ such that
$Cut_X(p)\cap B_{\eps}(q)= C_{k,\eps}(Cut_{\La_qX}(F))$
where $F\subset\La_qX$
is the finite set of directions of minimizing segments between $p$ and $q$.
\end{prop}

If $k>0$ and $X$ is a metric suspension with tip $p$, then $Cut(p)$ consists
of just one point, namely the antipode of $p$.

In all other cases,
induction over the dimension using Proposition \ref{loccon} yields
that $Cut(p)$ is a possibly empty, locally finite, piecewise totally geodesic
polyhedral complex of codimension one,
and ${\cal D}(p)$ is a locally finite polyhedron in $C_k(\La_pX)$ with
geodesic faces.
The conifold $X$ arises from ${\cal D}(p)$ by identifications on the boundary,
namely by isometric face pairings.

\begin{dfn}
[Dirichlet polyhedron]
${\cal D}(p)\subseteq C_k(\La_pX)$ is called the {\em Dirichlet polyhedron}
with respect to $p$.
\end{dfn}

In dimension $2$,
the Dirichlet polyhedra are polygons.
If $X$ is a cone surface,
then the vertices of ${\cal D}(p)$ correspond to
either smooth interior points of $X$ with $\geq3$ minimizing
segments towards $p$, to boundary points or to cone points.
In the latter cases there may exist only one minimizing segment to $p$.
If this happens for a cone point,
then the angle at the corresponding vertex of ${\cal D}(p)$ equals the cone
angle.
This is the only way, in which {\em concave} vertices of the Dirichlet polygon
can occur:
Every vertex of ${\cal D}(p)$ with angle $>\pi$ corresponds to a cone point
which is connected to $p$ by exactly one minimizing segment.

The discussion in dimension $3$ is analogous.
In particular,
if $X$ is a 3-conifold
then edges of Dirichlet polyhedra with dihedral angles $>\pi$
project via the exponential map to (parts of) singular edges with cone angles
$>\pi$.
Therefore we have the following strong restriction on the geometry of
Dirichlet polyhedra for cone angles $\leq\pi$:

\begin{prop}[Convexity]
In the case of cone angles $\leq\pi$,
the Dirichlet polyhedra are {\em convex}.
\end{prop}

The exponential map is a local isometry near the tip $O$ of $C_k(\La_pX)$.

\begin{dfn}
[Injectivity radius]
The {\em injectivity radius} in $p$, $0<r_{inj}(p)\leq r_{geod}(p)$,
is the radius of the largest open ball in $C_k(\La_pX)$ around $O$
on which $\exp_p$ is an embedding, i.e.\ it is maximal with the property
that all geodesic segments of length $<r_{inj}(p)$ starting in $p$ are
minimizing.
\end{dfn}

Since the cut locus $Cut(p)$ is closed, there exist cut points $q$ at minimal
distance $r_{inj}(p)$ from $p$.
The minimizing segments $\ol{pq}$ must have angles $\geq\frac{\pi}{2}$ with
the cut locus.
Since $diam(\La_qX)\leq\pi$,
there can be at most two minimizing segments $\ol{pq}$.
If there are two, they meet with maximal angle $\pi$ at $q$
and form together a geodesic loop with base point $p$ and midpoint $q$.
If there is a unique minimizing segment $\ol{pq}$
and if $q$ does not belong to the boundary,
then $q$ must lie on a (closed) singular edge with cone angle $\geq\pi$.
Note that this alternative cannot occur for cone angles $<\pi$.

\medskip
The injectivity radius varies continuously with $p$
on the smooth part and along singular edges.
However it converges to zero e.g.\ along sequences of smooth points
approaching the singular locus.
In the singular setting,
the injectivity radius is not the right measure for the simplicity of the
local geometry.
In order to measure up to which scale the local geometry is given by certain
simple models, the following modification turns out to be useful,
at least as long as the cone angles are $\leq\pi$.

\begin{dfn}[Cone injectivity radius]
\label{dfn:coneinjrad}
The {\em cone injectivity radius} $r_{cone-inj}(p)$ in $p$
is the supremum of all $r>0$ such
that the ball $B_r(p)$ is contained in a standard ball,
i.e.\ such that there exist $q\in X$ and $R>0$ with the following property:
$B_r(p)\subseteq B_R(q)$ and $\ol B_R(q)\cong C_{k,R}(\La_qX)$.
\end{dfn}

\subsection{Spherical cone surfaces with cone angles $\leq\pi$}

In this section we will discuss closed cone surfaces $\La$
with curvature $1$ and cone angles $\leq\pi$,
whose underlying topological surface is a $2$-sphere.
They occur as links of 3-dimensional cone manifolds with cone angles $\leq\pi$,
the class of cone manifolds mostly relevant for us in this paper.

\begin{prop}
[Classification]
\label{linkclass}
Let $\La$ be a spherical cone surface with cone angles $\leq\pi$ which is
homeomorphic to the 2-sphere.
Then $\La$ is isometric to either
\BI
\item the unit 2-sphere $S^2$,
\item the metric suspension $S^2(\alpha,\alpha)$ of a circle of length
  $\al\leq\pi$, or to
\item $S^2(\alpha,\beta,\gamma)$,
the double along the boundary of a spherical triangle with angles
  $\frac{\al}{2},\frac{\beta}{2},\frac{\ga}{2}\leq\frac{\pi}{2}$.
\EI
\end{prop}
\begin{proof}
The assertion is clear in the smooth case and we therefore assume that $\La$
has cone points.
Due to Gau\3-Bonnet, there can be at most three cone points.

If $\La$ has only one cone point $c$, then $\La-\{c\}$ is simply connected and
hence can be developed (isometrically immersed) into $S^2$.
A circle of small radius centered at $c$ cannot close up under the developing
map and we obtain a contradiction.
Thus $\La$ must have two or three cone points.

If $\La$ has two cone points, we connect them by a minimizing segment $\si$.
By cutting $\La$ open along $\si$ we obtain a spherical surface which is
topologically a disc and whose boundary consists of two edges of equal
length.
It can be developed into $S^2$ as well,
and it follows that the surface is a spherical bigon, i.e.\ the metric
suspension of an arc.
We obtain the second alternative of our assertion.

If $\La$ has three cone points, we connect any two of them by a minimizing
geodesic segment.
The segments don't intersect and they divide $\La$ into two spherical
triangles.
The triangles are isometric because they have the same side lengths, and we
obtain the third alternative.
\end{proof}

A consequence of the classification is the following description for the local
geometry of a cone 3-manifold with cone angles $\leq\pi$.

\begin{figure}
\begin{center}
 \includegraphics[scale=1]{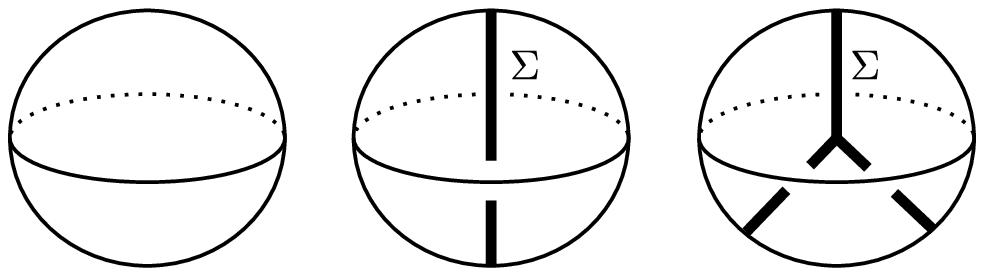}
 \end{center}
    \caption{ }\label{fig:2.1}
\end{figure}

\begin{cor}
    \label{localmodel}
    If $p$ is an interior point in a cone 3-manifold with cone angles
    $\leq\pi$, then a sufficiently small ball
    $B_{\eps}(p)$ is isometric to one of the following (see
    Figure~\ref{fig:2.1}):
\begin{itemize}
\item[--] a ball of radius $\eps$ in a smooth model space $\M_k^3$,
\item[--] a singular ball $C_{k,\eps}(S^2(\al,\al))$ with a singular axis of
  cone angle   $\al$, or
\item[--] a singular ball $C_{k,\eps}(S^2(\al,\beta,\ga))$
with three singular edges emanating from a singular vertex in the center.
\end{itemize}
In particular, the singular locus $\Si_X$ is a {\em trivalent} graph, i.e.\
its vertices have valency at most three.
\end{cor}

In the remainder of this section,
we collect some properties of spherical cone surfaces used later.

\begin{lem}
\label{2dirbigangle}
Let $\La$ be as in Proposition \ref{linkclass}.
Then $\La$ does not contain three points with pairwise distances
$>\frac{2\pi}{3}$.
\end{lem}
\begin{proof}
This is a direct implication of the lower curvature bound $1$ because
the circumference of geodesic triangles has length $\leq2\pi$.
\end{proof}

\begin{dfn}
A {\em turnover} is a cone surface
which is homeomorphic to the 2-sphere
and which has three cone points, all with cone angle $\leq\pi$.
\end{dfn}

Geometrically, a turnover is the double along the boundary of a triangle in a
2-dimensional model space $\M^2_k$ with angles $\leq\frac{\pi}{2}$.

\begin{lem}
\label{fatlink}
(i)
A spherical turnover $\La$ has diameter $\leq\frac{\pi}{2}$.

(ii)
If $\La$ is a spherical turnover with $diam(\La)=\frac{\pi}{2}$,
then at least two of the three cone angles equal $\pi$.
If two points $\xi,\eta\in\La$ have maximal distance $\frac{\pi}{2}$
then at least one of them, say $\xi$, is a cone point,
and $\eta$ lies on the minimizing segment joining the other two cone points,
and these must have cone angles $=\pi$.
\end{lem}
\begin{proof}
(i)
Suppose that $d(\xi,\eta)\geq\frac{\pi}{2}$ and that $\zeta$ is a
cone point $\neq\xi,\eta$. Any geodesic triangle
$\De(\xi,\eta,\zeta)$ has angle $\leq\frac{\pi}{2}$ at $\zeta$.
We denote
$rad(\La,\zeta):=\max d(\zeta,\cdot)$.
Since
$rad(\La,\zeta)\leq\frac{\pi}{2}$, hinge comparison implies that
$d(\xi,\eta)\leq\frac{\pi}{2}$.

(ii) In the case of equality it follows that the cone angle at
$\zeta$ equals $\pi$ and that one of the points $\xi$ and $\eta$,
say $\xi$, has distance $\frac{\pi}{2}$ from $\zeta$. If $\xi$
were not a cone point, then it would lie on the segment connecting
the two cone points $\neq\zeta$ and only $\zeta$ would have
distance $\frac{\pi}{2}$ from $\xi$, contradicting
$d(\xi,\eta)=\pi/2$. Hence $\xi$ must be a cone point, and it
follows that $\eta$ lies on the segment joining $\zeta$ and the
cone point $\neq\zeta,\xi$.
\end{proof}

\begin{lem}
\label{fatturnover}
For $\al<\pi$ there exists $D=D(\al)<\frac{\pi}{2}$ such that:
If $\La$ is a spherical turnover
with at least two cone angles $\leq\al$
then $diam(\La)\leq D(\al)$.
\end{lem}
\begin{proof}
$\La$ is the double of a spherical triangle $\De$ with two angles
$\leq\al/2$ and third angle $\leq\frac{\pi}{2}$. Since the angle
sum of a spherical triangle is $>\pi$, all angles of $\De$ are
$>\frac{\pi-\al}{2}$. Such triangles can Gromov-Hausdorff converge
to a point, but not to a segment. Hence the Gromov-Hausdorff
closure of the space of turnovers as in the lemma is compact and
contains as only additional space the point. It follows that the
diameter assumes a maximum $D(\al)$ on this space of turnovers. By
part (ii) of Lemma \ref{fatlink}, we have $D(\al)<\frac{\pi}{2}$.
\end{proof}

\begin{lem}
\label{thicklink}
For $\al<\pi$ and $0<d\leq\frac{\pi}{2}$
there exists $r=r(\al,d)>0$ such that:
If $\La$ is a spherical turnover
with diameter $\geq d$ and cone angles $\leq\al$,
then it contains an embedded smooth round disc with radius $r$.
\end{lem}
\begin{proof}
The turnover $\La$ is the double of a spherical triangle $\De$
with acute angles $\leq\al/2$ and a lower diameter bound. Since
the angle sum of spherical triangles is $>\pi$, we also have the
positive lower bound $\pi-\al$ for the angles of $\De$. Such
triangles have a lower bound on their inradius, whence the claim.
\end{proof}

\subsection{Compactness for spaces of thick cone manifolds}
\label{sec:compact}

The space of pointed cone 3-manifolds with bounded curvature
is {\em pre}compact in the Gromov-Hausdorff topology
by Gromov's compactness theorem, cf.\ \cite{GLP},
because the volume growth is at most as strong as in the model
space.
Due to general reasons,
the limit spaces in the Gromov-Hausdorff closure are spaces with curvature
bounded below.
We will show that, under appropriate assumptions,
limits of cone 3-manifolds are
still cone 3-manifolds.

\begin{dfn}
\label{dfn:thick}
For $\rho>0$,
a cone manifold $X$ is said to be {\em $\rho$-thick}
(at a point $x$)
if it contains an embedded smooth standard ball of radius $\rho$
(centered at $x$).
Otherwise $X$ is called {\em $\rho$-thin}.
\end{dfn}

For $\kappa,i,a>0$
we denote by ${\cal C}_{\kappa,i,a}$
the space of pointed cone 3-manifolds $(X,p)$
with constant curvature $k\in[-\kappa,\kappa]$,
cone angles $\leq\pi$
and base point $p$ which satisfies $r_{inj}(p)\geq i$ and $area(\La_pX)\geq a$.
Let ${\cal C}_{\kappa,i}:={\cal C}_{\kappa,i,4\pi}$
be the subspace of cone manifolds with smooth base point;
they are $i$-thick at their base points.

\begin{thm}[Compactness for thick cone manifolds with cone angles $\leq\pi$]
\label{thm:thickconecpt}
The spaces
${\cal C}_{\kappa,i}$ and ${\cal C}_{\kappa,i,a}$
are compact in the Gromov-Hausdorff topology.
\end{thm}

The main step in the proof of the theorem is the following result.

\begin{prop}[Controlled decay of injectivity radius]
\label{decayofinj}
For $\kappa\geq0$, $R\geq i>0$ and $a>0$
there exist $r'(\kappa,i,a,R)\geq i'(\kappa,i,a,R)>0$
such that the following holds:

Let $X$ be a closed cone 3-manifold with curvature $k\in[-\kappa,\kappa]$
and cone angles $\leq\pi$. Let $p\in X$ be a point with $r_{inj}(p)\geq i$
and $area(\La_pX)\geq a$. Then for every point $x\in B_R(p)$ the ball
$B_{i'}(x)$ is contained in a standard ball with radius $\leq r'$. In
particular, we have $r_{cone-inj}\geq i'$ on $B_R(p)$.
\end{prop}

By a {\em standard ball} we mean the $k$-cone over a spherical cone surface
homeomorphic to the 2-sphere, cf.\ Definition \ref{dfn:coneinjrad}.

\begin{proof}
{\em Step 0.}
It follows from the classification of links, cf.\ Proposition
\ref{linkclass}, that $\La_pX$ contains a smooth standard disc with radius
bounded below in terms of $a$,
and hence the ball $B_i(p)$ contains an embedded {\em smooth}
standard ball with a lower bound on its radius in terms of $\kappa,i$ and
$a$.
We may therefore assume without loss of generality that $p$ is a smooth
point.

{\em Step 1.}
We have a lower bound
$vol(B_R(x)-B_{i/2}(x))\geq v(\kappa,i)>0$
because $B_R(x)-B_{i/2}(x)$
contains a smooth standard ball of radius $\geq i/4$.
Let $A_x\subseteq\La_xX$ denote the subset of
initial directions of minimizing geodesic segments with length $\geq i/2$.
The lower bound for
the volume of the annulus $B_R(x)-B_{i/2}(x)$
implies a lower bound $area(A_x)\geq a_1(\kappa,i,R)>0$.

{\em Step 2.} By triangle comparison, there exists for $\eps>0$ a
number $l=l(\kappa,i,\eps)>0$ such that: Any geodesic loop of
length $\leq 2l$ based in $x$ has angle $\geq\frac{\pi}{2}-\eps$
with all directions in $A_x$. The same holds for the angles of
$A_x$ with segments of length $\leq l$ starting in $x$ and
perpendicular to the singular locus $\Si_{X}$. Thus, if
$r_{inj}(x)\leq l$, then minimizing segments from $x$ to closest
cut points must have angles $\geq\frac{\pi}{2}-\eps$ with all
directions in $A_x$, cf.\ our discussion of the cut locus in
Section \ref{sec:expcutinj}. We use this observation to obtain
lower bounds for the injectivity radius.

\begin{lem}
\label{oppdir}
For $a'>0$
there exists $\eps=\eps(a')>0$
such that:

Let $\La$ be a spherical cone surface homeomorphic to the 2-sphere and
with cone angles $\leq\pi$.
Let $A\subset\La$ be a subset with $area(A)\geq a'$.
Then $\La=N_{\frac{\pi}{2}-\eps}(A)$ if $\La$ is a turnover.
If $\La$ has 0 or 2 cone points,
then there exists a point $\eta$
such that
$\La - N_{\frac{\pi}{2}-\eps}(A) \subset B_{\frac{\pi}{2}-\eps}(\eta)$.
In the case that $\La$ has two cone points,
$\eta$ can be chosen as a cone point.
\end{lem}
\begin{proof}
When $\La$ is a turnover,
the description in Lemma \ref{fatlink}
of segments of maximal length $\frac{\pi}{2}$ implies:
Points in $\La$ with radius (Hausdorff distance from $\La$) close to
$\frac{\pi}{2}$ must be close to one of the three minimizing segments
connecting cone points, i.e.\ lie in a region of small area.
Hence $A$ contains points with radius $<\frac{\pi}{2}-\eps$ for sufficiently
small $\eps>0$ depending on $area(A)$.

If $\La$ has 0 or 2 cone points then it is isometric to the unit sphere $S^2$
or the metric suspension of a circle with length $\leq\pi$,
cf.\ the classification in Proposition \ref{linkclass}.
In both cases the assertion is easily verified.
\end{proof}

We choose $\eps:=\eps(a_1)$
with $a_1=a_1(\kappa,i,R)$ as in step 1,
and accordingly $l=l(\kappa,i,\eps)=l(\kappa,i,R)$.

{\em Step 3.}
For a singular vertex $x$ Lemma \ref{oppdir} implies that
$r_{inj}(x)\geq i_1=i_1(\kappa,i,R):=l(\kappa,i,R)>0$.

{\em Step 4.}
Assume that $x$ is a singular point
with $r_{inj}(x)\leq i_1=l$
at distance $\geq i_1/4$ from all singular vertices,
and choose the singular direction $\eta_x\in\La_xX$ according to
Lemma \ref{oppdir}.
By the assumption on the injectivity radius,
there exists a geodesic loop $\la$ of length $\leq 2l$ based at $x$
or a segment $\ol{xy}$ of length $\leq l$
meeting $\Si_{X}$ orthogonally at at a point $y$.
Either of them has angles $\geq\frac{\pi}{2}-\eps$ with the directions in
$A_x$ and therefore angles $\leq\frac{\pi}{2}-\eps$ with the direction
$\eta_x$.

In the case of a loop,
consider the geodesic variation of $\la$ moving its base point
with unit speed in the direction $\eta_x$.
Since both ends of the loop have angle $\leq\frac{\pi}{2}-\eps$ with $\eta_x$,
the first variation formula implies that
the length of $\la$ decreases at a rate $\leq-2\sin\eps$.
Similarly, in the case of a segment, we obtain that
$r_{inj}$ decreases at a rate $\leq-\sin\eps$.
It follows that
$r_{inj}(x)\geq \frac{i_1}{4}\cdot \sin\eps=:i_2=i_2(\kappa,i,R)$.

{\em Step 5.}
Suppose now that $x$ is a smooth point
with $r_{inj}(x)\leq i_2$ at distance $\geq i_2/4$ from $\Si_{X}$.
We choose the direction $\eta_x\in\La_xX$ according to Lemma \ref{oppdir}.
As in Step 4,
we see that $r_{inj}$ decays in the direction $\eta_x$
with rate $\leq-\sin\eps$.
It follows that
$r_{inj}(x)\geq \frac{i_2}{4}\cdot \sin\eps=:i_3=i_3(\kappa,i,R)$.

{\em Conclusion.}
The assertion holds for $r':=i_1$ and $i':=i_3$.
\end{proof}

{\em Proof of Theorem \ref{thm:thickconecpt}.}
Let $(Y,q)$ be an Alexandrov space in the Gromov-Hausdorff closure of
${\cal C}_{\kappa,i,a}$.
It is the Gromov-Hausdorff limit of a sequence of pointed cone manifolds
$(X_n,p_n)\in{\cal C}_{\kappa,i,a}$.
For a point $y\in Y$,
we pick points $x_n\in X_n$ converging to $y$.
The metric ball $B_{\rho}(y)\subset Y$ is then the
Gromov-Hausdorff limit of the balls $B_{\rho}(x_n)$ in the approximating cone
manifolds $X_n$.

Proposition~\ref{decayofinj} yields numbers $r'\geq i'>0$
such that each ball $B_{i'}(x_n)$ is contained in a standard ball
$B_{r'_n}(x'_n)$
with radius bounded above by $r'_n\leq r'$.
Moreover, the lower bound on the volumes of the balls $B_i(p_n)$
yields a uniform estimate $area(\La_{x'_n}X_n)\geq a'(\kappa,i,a,d(q,y))>0$.

It is clear from the classification of links in Proposition \ref{linkclass}
that the space ${\cal C}^2_{a'}$
of spherical cone surfaces homeomorphic to the 2-sphere
with cone angles $\leq\pi$ and area $\geq a'$ is Gromov-Hausdorff compact.
Thus, after passing to a subsequence, we have that the links $\La_{x'_n}X_n$
converge to a cone surface $\La\in{\cal C}^2_{a'}$.
Moreover,
$r'_n\ra r'_{\infty}\leq r'$ and $k_n\ra k_{\infty}$
where $k_n$ denotes the curvature of $X_n$.
It follows that
$B_{r'_n}(x'_n)\cong C_{k_n,r'_n}(\La_{x'_n}X_n)\ra
C_{k_{\infty},r'_{\infty}}(\La)$.
This means that $Y$ is a cone manifold locally at $y$.
It is then clear that $Y\in {\cal C}_{\kappa,i,a}$.
\qed

\medskip
In our context, Gromov-Hausdorff convergence implies a
stronger type of convergence,
namely a version of bilipschitz convergence for cone manifolds.
Recall that, for $\eps>0$, a map $f:X \to Y$ between metric spaces is
called a $(1+\varepsilon)$-bilipschitz embedding
if
\[ (1+\eps)^{-1}\cdot d(x_1,x_2)<d(f(x_1),f(x_2))<(1+\eps)\cdot d(x_1,x_2)\]
holds for all points $x_1,x_2\in X$.

\begin{dfn}
[Geometric convergence]
A sequence of pointed cone 3-manifolds
$(X_n,x_n)$
{\em converges geometrically} to a pointed
cone 3-manifold $(X_{\infty},x_{\infty})$
if for every $R>0$ and $\varepsilon>0$
there exists $n(R,\eps)\in\N$ such that
for all $n\geq n(R,\eps)$
there is a $(1+\varepsilon)$-bilipschitz embedding
$f_n:B_R(x_{\infty})\to X_n$ satisfying:
\begin{itemize}
\item[(i)] $d(f_n(x_{\infty}),x_n)<\varepsilon$,
\item[(ii)] $B_{(1-\eps)\cdot R}(x_n)\subset f_n(B_R(x_{\infty}))$, and
\item[(iii)] $f_n(B_R(x_{\infty})\cap\Sigma_{\infty})=
f_n(B_R(x_{\infty})) \cap\Sigma_n$.
\end{itemize} \end{dfn}

Note that the definition also implies the inclusion
$f_n(B_R(x_{\infty}))\subset B_{R(1+\varepsilon)+\varepsilon}(x_n)$.

A standard argument (cf. \cite[ch.\ 3.3]{BoP}) using the strong local structure
of cone 3-manifolds and the controlled decay of injectivity radius
(Proposition \ref{decayofinj}) shows that within the spaces
${\cal C}_{\kappa,i}$ and ${\cal C}_{\kappa,i,a}$
the Gromov-Hausdorff topology and the pointed
bilipschitz topology are equivalent.
We therefore deduce from Theorem \ref{thm:thickconecpt}:

\begin{cor}
\label{thm:compactness}
Let $(X_n)$ be a
sequence of cone 3-manifolds with curvatures $k_n\in [-\kappa,\kappa]$,
cone angles $\leq\pi$, and possibly with totally geodesic boundary.
Suppose that, for some $\rho>0$,
each $X_n$ is $\rho$-thick at a point $x_n \in X_n$.

Then, after passing to a subsequence,
the pointed cone 3-manifolds $(X_n,x_n)$ converge geometrically
to a pointed cone 3-manifold $(X_{\infty},x_{\infty})$, with
curvature $k_{\infty} = \lim\limits_{n \to \infty} k_n$.
\end{cor}

Note that
the case with totally geodesic boundary follows from the closed case by
doubling along the boundary.

\section{Noncompact Euclidean cone 3-manifolds}
\label{sec:euc}

A heuristic guideline to describe the geometry of the thin part of
cone $3$-manifolds, (i.e.\ the possibilities for the local
geometry on a {\em uniform} small scale) is that global results
for noncompact Euclidean cone manifolds correspond to local
results for cone manifolds of bounded curvature. For instance, in
the smooth case, the fact that there is a short list of noncompact
Euclidean manifolds reflects the Margulis lemma for complete
Riemannian manifolds of bounded curvature.

We show in this section that there is still a short list of
noncompact Euclidean cone 3-manifolds with cone angles $\leq\pi$.
The corresponding local results for cone manifolds with bounded
curvature will be discussed in Section \ref{sec:conelowdiambound}.

\begin{thm}
[Classification]
\label{thm:euclid}
Every noncompact Euclidean cone 3-manifold $E$ with cone
angles $\leq\pi$ belongs to the following list:

\begin{itemize}
\item
smooth flat 3-manifolds,
i.e.\ line bundles over the 2-torus or the Klein bottle,
and plane bundles over the circle;
\item
complete Euclidean cones $C_0(\La)$ over spherical cone surfaces $\La$
with cone angles $\leq\pi$ and which are homeomorphic to $S^2$;
\item
bundles over a circle or a compact interval
with fiber a smooth Euclidean plane or a singular plane $\M^2_0(\theta)$
with $\theta\leq\pi$;
\item
$\R$ times a closed flat cone surface with cone angles $\leq\pi$;
bundles over a ray with fiber a closed flat cone surface with cone angles
$\leq\pi$.
\end{itemize}
\end{thm}

By {\em bundles} we mean metrically locally trivial bundles.
Line bundles refer to bundles with fiber $\cong\R$.
In the case of bundles over a ray or a compact interval,
the fibers over the endpoints are singular with index two.

We give a short direct proof of the classification without using general
results for nonnegatively curved manifolds such as the Soul Theorem or the
Splitting Theorem, although the ideas are of course related.
The existence of a soul in our special situation is actually a direct
consequence of the list given in Thm.~\ref{thm:euclid}.
Recall that a {\em soul}
is a totally convex compact submanifold of dimension
$<3$ with boundary either empty or consisting of singular edges with cone
angle $\pi$.

\begin{cor}\label{cor:soul} Every noncompact Euclidean cone 3-manifold
 with cone angles $\leq\pi$ is a metrically locally trivial bundle over a soul with fiber
 a complete cone, or a quotient of such a bundle by an
isometric involution.

In particular, the soul is a point iff $E$ is a cone.
\end{cor}

Before giving the proof of Theorem \ref{thm:euclid}
we establish some preliminary lemmas.
Since $E$ is noncompact,
there are globally minimizing rays emanating from every point $x\in E$.
We denote by $R_x\subseteq\La_xE$ the closed set of initial directions of rays
starting in $x$.

\begin{lem}
\label{setofraysconvex}
(i)
$R_x$ is convex,
i.e.\ with any two directions $\xi$ and $\eta$, possibly coinciding,
$R_x$ contains all arcs $\ol{\xi\eta}$ of length $<\pi$.

(ii)
If $x\in\Si_E$,
every cone point of $\La_xE$ at distance $<\frac{\pi}{2}$ from $R_x$
belongs to $R_x$.
\end{lem}
\begin{proof}
(i) The convexity of the Dirichlet polyhedron ${\cal
D}(x)\subseteq C_0(\La_xE)$, cf.\ Section~\ref{sec:expcutinj},
implies that $R_x$ is convex.

(ii)
Suppose that $\xi\in\La_xE$ is a cone point
and $\eta$ is a point in $R_x$ with $d(\xi,\eta)<\frac{\pi}{2}$.

We consider first the case when the cone angle at $\xi$ is $<\pi$.
If $\La_xE$ is the metric suspension of a circle
then there exists a loop of length $<\pi$ based at $\eta$
and surrounding $\xi$.
It follows that $\xi$ is contained in the convex hull of $\eta$
and hence $\xi\in R_x$.
If $\La_xE$ is a spherical turnover,
we cut $\La_xE$ open along $Cut(\xi)$
and obtain a convex spherical polygon with $\xi$ as cone point.
Inside the polygon we find a loop as before.

We are left with the case
that the cone angle at $\xi$ equals $\pi$.
Let $\rho_{\xi}\subset C_0(\La_xE)$ be the singular ray in direction $\xi$.
Observe that,
if $z$ is a point on $\rho_{\xi}$ different from its initial point $x$,
$\ol{yz}$ is a segment perpendicular to $\rho_{\xi}$
and $B$ is a (small) ball around $y$,
then the convex hull of $B$ in $C_0(\La_xE)$ contains $z$.
Now the ray $\rho_{\eta}$ is contained in ${\cal D}(x)$.
Since ${\cal D}(x)$ is convex and has non-empty interior,
arbitrarily close to every point of $\rho_{\eta}$ we find interior points of
${\cal D}(x)$.
Our observation therefore implies
that $\rho_{\xi}\subset{\cal D}(x)$ and $\xi\in R_x$.
\end{proof}

Let $x$ be a point with $r_{inj}(x)<\infty$,
i.e.\ $Cut(x)\neq\emptyset$ and $R_x$ is a proper subset of $\La_xE$.
We then have as further restriction on $R_x$
that there exists a direction of angle $\geq\frac{\pi}{2}$ with $R_x$.
This follows from the next result by
looking at shortest segments to the cut locus:

\begin{lem}
\label{rightangle}
Suppose that $\zeta\in\La_xE$ is the
initial direction of a geodesic loop based at $x$
or of a segment $\ol{xy}$ perpendicular
to $\Sigma_{E}$ at $y$.
Then $\angle_x(\zeta,R_x)\geq\frac{\pi}{2}$.
\end{lem}
\begin{proof}
Let $r:[0,\infty)\ra E$ be a ray starting in $x$.
In the case of a loop $\la$,
the assertion follows by applying angle comparison to the isosceles geodesic
triangle with $\la$ as one of its side and twice the segment $r|_{[0,t]}$ as
the other two sides, and by letting $t\ra\infty$.
Comparison is applied to the angles adjacent to the non-minimizing side $\la$.

In the second case, the argument is similar.
We consider instead the geodesic triangle with sides $\ol{xy}$, $r|_{[0,t]}$
and a minimizing segment $\ol{yr(t)}$ as third side,
and use that every direction at the singular point $y$ has angle
$\leq\frac{\pi}{2}$ with $\ol{xy}$.
\end{proof}

Any of the Lemmas \ref{setofraysconvex}
or \ref{rightangle} implies:

\begin{lem}
\label{vertexthencone}
If $v$ is a singular vertex with $diam(\La_vE)<\frac{\pi}{2}$
then $R_v=\La_vE$ and
$\exp_v$ is a global isometry, i.e.\ $E\cong C_0(\La_vE)$.
\end{lem}

Regarding non-vertex singular points,
Lemma \ref{setofraysconvex} implies:

\begin{lem}
\label{raysatsingedge} Let $x\in\Si_E^{(1)}$. Then either there is
a singular ray initiating in $x$, or all rays emanating in $x$ are
perpendicular to $\si$, where $\sigma$ is the singular edge of
$\Si_E^{(1)}$ containing $x$. In the latter case, if the cone
angle at $\si$ is $<\pi$, then every direction in $x$
perpendicular to $\si$ is the initial direction of a ray.
\end{lem}

{\em Proof of Theorem~\ref{thm:euclid}.}
The smooth case is well-known,
and we assume that the singular locus $\Sigma_E$ is non-empty.

{\em Part 1: The case when cone angles are $<\pi$.}
If $E$ contains a singular vertex,
then $E$ is a cone by Lemmas \ref{fatturnover} and \ref{vertexthencone}.
If $E$ contains a closed singular geodesic, then
Lemma \ref{raysatsingedge} implies that
the exponential map is an isometry from the normal bundle of $\si$ onto $E$,
i.e.\ $E$ is a metrically locally trivial bundle over $\si$ with fiber a plane
with cone point.
We are left with the case that $\Si_{E}$
consists of lines, i.e.\ of complete noncompact geodesics.
We assume that $E$ is not a cone, i.e.\ $r_{inj}<\infty$ everywhere.

Let $\si$ be a singular edge with cone angle $\theta$.
Assume that there exists a ray in $E$ perpendicular to $\si$ in a point $x$.
The singular model space $C_0(\La_xE)$ is isometric to the product
$\M^2_0(\theta)\times\R$.
Note that
$\M^2_0(\theta)$ contains no unbounded proper convex subset
because $\theta<\pi$.
It follows that
${\cal D}(x)$ splits metrically as the product of
$\M^2_0(\theta)$ with a closed connected subset $I$ of $\R$.
Since $E$ is not a cone, $I$ is a proper subset of $\R$
and $\D{\cal D}(x)$ consists of one or two singular planes
$\cong\M^2_0(\theta)$.
Under our assumption that cone angles are $<\pi$,
the points in $\D{\cal D}(x)$ away from the singular axis project to
{\em smooth} cut points.
It follows that $\si$ closes up, contradiction.

Hence there are no rays in $E$ perpendicular to $\si$.
Lemma \ref{raysatsingedge} leaves the possibility
that from each point $x\in\si$ emanates at least one {\em singular} ray.
Let us denote by $A,B\subseteq\si$ the sets of initial points of singular rays
directed to the respective ends of $\si$.
Both subsets $A$ and $B$ are closed, connected and unbounded.
So either they have non-empty
intersection or one of them, say $A$, is empty and $B=\si$.
In the latter case, $\si$ would be globally minimizing and we obtain a
contradiction with $A=\emptyset$.
Only the first case is possible,
i.e.\ there exists a point $x$ on $\si$ which divides $\si$ into
two rays.

Then ${\cal D}(x)$ contains the entire singular axis of $\M^3_0(\theta)$
and, by convexity, it splits as
${\cal D}(x)\cong\R\times C_x$
where $C_x\subset M^2_0(\theta)$ denotes the cross section through $x$.
Since $E$ is not a cone,
$C_x$ is a {\em proper} convex subset.
It follows that $C_x$ is compact
and hence a finite-sided polygon with one cone point.
Accordingly,
$\D{\cal D}(x)$ consists of finitely many strips of finite width.

Away from the edges the identifications on $\D{\cal D}(x)$ are given by an
involutive isometry $\iota$,
and on the edges by its continuous extension.
It must preserve the direction parallel to the singular axis of
$M^3_0(\theta)$.
Moreover,
$\iota$ preserves distance from $x$.
It follows that $\iota$ maps $\D C_x$ onto itself and
$C_x$ projects to an embedded
totally geodesic closed surface $S\subset E$ with at least one cone point.
Due to Gau\3-Bonnet,
$S$ must be a turnover
and is in particular two-sided.
It follows that $E\cong\R\times S$.

{\em Part 2: The general case of cone angles $\leq\pi$.}
We expand the above analysis.
We assume again that $E$ is not a cone, i.e.\ $r_{inj}<\infty$ everywhere.

For $x\in E$,
let us denote by
$\dot\D{\cal D}(x)$
the smooth part of the boundary of the Dirichlet polyhedron,
i.e.\ the complement of the edges.
The identifications on $\D{\cal D}$ are the continuous extension of an
involutive self-isometry $\iota$ of $\dot\D{\cal D}(x)$.
Unlike in the case of cone angles $<\pi$, $\iota$ may now have fixed points;
the fixed point set $Fix(\iota)$
is a union of segments and projects to the interior points on singular edges
with cone angle $\pi$
which are connected to $x$ by exactly one minimizing segment.

{\em Step 1.}
Let $x$ be an interior point of a singular edge $\si$ with cone angle
$\theta\leq\pi$
and suppose that $x$ is not the initial point of a {\em singular} ray.
Then, starting at $x$, $\si$ remains in
both directions minimizing only for finite time, i.e.\ ${\cal D} (x)$
intersects the singular axis of
$C_0(\La_xE)\cong \M^3_0(\theta)$ in a compact subsegment $I$.
By convexity, we have ${\cal D}(x)\subseteq I\times\M^2_0(\theta)$,
compare the proof of part (ii) of Lemma \ref{setofraysconvex}.
The cross section $C_x\subseteq\M^2_0(\theta)$
of ${\cal D}(x)$ perpendicular to $I$ through $x$
is an unbounded convex subset.

{\em Step 1a.} If $C_x=\M^2_0(\theta)$, then ${\cal D} (x)\cong
I\times \M^2_0(\theta)$ and ${\cal D}(x)$ consists of two copies
of $\M^2_0(\theta)$. The involution $\iota$ on $\dot\D{\cal D}(x)$
either exchanges the boundary planes or it is a reflection on each
of them. By a reflection on the singular plane $\M^2_0(\theta)$ we
mean an involutive isometry whose fixed point set is the union of
two rays emanating from the cone point into ``opposite''
directions with angle $\frac{\theta}{2}$. Thus $E$ is a bundle
with fiber $\cong\M^2_0(\theta)$ over a circle or a compact
interval; in the latter case the fibers over the endpoints of the
interval are singular with index two, meaning that they are
index-two branched subcovers of the generic fiber.

{\em Step 1b.} If $C_x$ is a proper subset of $\M^2_0(\theta)$,
then $\theta=\pi$ because $C_x$ is unbounded. There is a unique
ray $r\subset {\cal D}(x)$ with initial point $x$. Let $H$ be the
half plane in $\M^3_0(\pi)$ bounded by the singular axis and
containing $r$. Cutting ${\cal D}(x)$ open along $H$ yields a
convex polyhedron ${\cal D}'$ which splits as ${\cal
D}'\cong\R\times P$ where $P$ denotes the cross section containing
$I$. The cross section $P$ is a compact convex polygon with $I$ as
one of its sides and angles $\leq\frac{\pi}{2}$ at both endpoints
of $I$.

The cone manifold $E$ arises from ${\cal D}'$ by identifications
on the boundary. As before, away from the edges they are given by
an isometric involution $\iota'$ with one-dimensional fixed point
set. The involution $\iota'$ carries lines to lines and, since it
preserves distance from $x$, it also preserves $\D P$. The fixed
point set of $\iota'$ consists of midlines of strips in ${\cal
D}'$ and of edges of $P$, in our situation including $I$. After
performing the identifications, $P$ becomes a compact totally
geodesic cone surface $S\subset E$. The boundary $\D S$ is a union
of singular edges with cone angle $\pi$. Every corner of $\D S$ is
the initial point of a singular ray perpendicular to $S$, and the
angle at the corner equals half the cone angle of the ray and
hence is $\leq\frac{\pi}{2}$. We obtain that $E$ is a line bundle
over $S$ with singular fibers (rays) over the boundary.

The cone manifold $E$ can be described as a bundle in a different
way. Let us denote by $P_t$ the cross section $\{t\}\times P$ of
${\cal D}'$ where we identify $P$ with $P_0$. Then, for $t>0$, the
union of the two cross sections $P_{\pm t}$ projects to a totally
geodesic closed cone surface $S_t\subset E$. All the surfaces
$S_t$ are canonically isometric, say to a surface $\hat S$. We see
that $E$ fibers over $[0,\infty)$ with fiber $\hat S$; the
singular fiber over $0$ is isometric to $S$ and obtained from
$\hat S$ by dividing out a reflection.

{\em Step 2.}
In the following we can assume
that in each singular point initiates a {\em singular} ray.
It is a consequence
that all singular edges emanating from singular vertices are rays.
If there exists a singular vertex $v$,
Lemma \ref{setofraysconvex} implies
that $R_v=\La_vE$ and $E$ is a cone,
contrary to our assumption.
Hence $E$ contains no singular vertices
and $\Si_{E}$ is a union of lines.

As in part 1 it follows that each singular line $\si$ contains a
point $x$ dividing it into two rays, and ${\cal D} (x)\cong
\R\times C_x$ where $C_x$ is the cross section of ${\cal D} (x)$
through $x$. The section $C_x$ is a proper convex subset of
$\M^2_0(\theta)$ where $\theta\leq\pi$ is the cone angle at $\si$.
It is proper because $E$ is not a cone.

{\em Step 2a.}
If $C_x$ is bounded,
then $\D{\cal D}(x)$ is a finite union of strips of finite width.
We argue as in step 1b and obtain that $E$ splits off an $\R$-factor or fibers
over a ray.

{\em Step 2b.}
If $C_x$ is unbounded, then $\theta=\pi$ because $C_x$ is a proper
subset of $\M^2_0(\theta)$.
Moreover, $C_x$ is a Euclidean surface with one cone point of angle $\pi$ and
one boundary line;
it can be constructed from a flat strip by identifying one boundary line to
itself by a reflection.
Hence $\D{\cal D}(x)$ is a smooth Euclidean plane.
The involution $\iota$ preserves $d(x,\cdot)$ and therefore fixes
the unique point on $\D{\cal D}(x)$ closest to $x$.
It follows that $\iota$ is a reflection at a line through $x$,
and $E$ is a plane bundle over a compact interval.
Over each endpoint of the interval
there is a singular fiber isometric to a half plane and bounded by a singular
line with cone angle $\pi$.

\medskip
The proof of Theorem \ref{thm:euclid} is now complete.
\qed

\section{The local geometry of cone 3-manifolds with lower diameter bound}
\label{sec:conelowdiambound}

\subsection{Umbilic tubes}
\label{sec:tubes}

We start by describing certain simple cone manifolds which serve as local
models and building blocks for the thin part of arbitrary cone 3-manifolds.

\medskip
The smooth 3-dimensional model space $\M^3_k$ of constant curvature
$k$ can be viewed as the complete $k$-cone over the unit 2-sphere.
More generally, we define for a spherical cone surface $\La$
the {\em singular} model space $\M_k^3(\La)$ as the complete $k$-cone
$C_k(\La)$.
For the metric suspension
$\La(\al):=C_1(S^1(\al))$
of the circle $S^1(\al)$ with length $\al<2\pi$
we obtain the model space $\M^3_k(\al):=\M_k^3(\La(\al))$ of curvature $k$ with
a singular axis of cone angle $\al$.
The singular model spaces $\M_k^3(\La)$ serve as local models for cone
3-manifolds,
cf.\ Corollary \ref{localmodel}.

\medskip
Recall that an embedded connected surface $S$ in a model space $\M^3_k$ is
called {\em umbilic}
if in each point both principal curvatures are equal.
It follows that the principal curvatures in all points have the same value
which we denote by $pc(S)$.
The local extrinsic geometry of the surface is determined by its principal
curvature.
Its intrinsic Gau\3 curvature is
given by $k_S=k+pc(S)^2$.
We call $S$ {\em spherical}, {\em horospherical}, respectively
{\em hyperspherical}, depending on whether $k_S>0$, $k_S=0$ or $k_S<0$.

The model spaces $\M^3_k$ admit the following {\em umbilic foliations},
i.e.\ foliations by umbilic surfaces:
\begin{itemize}
\item
For all $k$ the {\em spherical} foliation
by distance spheres around a fixed point;
\item
for $k\leq0$ the {\em parabolic foliation}
by parallel planes if $k=0$,
respectively by horospheres centered at a fixed point at infinity if $k<0$;
\item
for $k<0$ the {\em hyperbolic foliation}
by equidistant surfaces from a fixed totally geodesic plane.
\end{itemize}
The leaves of these foliations are spherical, horospherical,
respectively hyperspherical.

We proceed to construct certain singular spaces with umbilic foliations.
Fix a cone surface $S$ with curvature $k_S\geq k$,
and let $L$ be a leaf of an umbilic foliation of $\M_k^3$ with curvature
$k_S$.
The type of the foliation depends on the sign of $k_S$.
We can develop the universal cover
$\widetilde{S^{smooth}}$ along $L$,
i.e.\ there exists an isometric immersion
$dev:\widetilde{S^{smooth}}\ra L$.
Let $N$ be a unit normal vector field along $L$
and consider the metric obtained from
pulling back the Riemannian metric of model space via
\[
\begin{array}{rcl}
\widetilde{S^{smooth}}\times\R&\lra& \M^3_k\\
 (x,t)&\mapsto &\exp(tN(dev(x))).
 \end{array}
\]
We choose the maximal open interval $I$ containing $0$
such that the Riemannian metric on $\widetilde{S^{smooth}}\times I$
is {\em non-singular}.
This metric has constant curvature $k$ and
descends to $S^{smooth}\times I$.

\begin{dfn}
We call the cone 3-manifold resulting
from metric completion of $S^{smooth}\times I$
the {\em complete $k$-tube over $S$}
and denote it by
$Tube_k(S)$.
\end{dfn}

The tubes have natural foliations
by umbilic surfaces equidistant from $S$;
the leaves are homothetic to $S$.
To each cone point of $S$ corresponds a singular edge of $Tube_k(S)$.
If $k_S>0$,
then $Tube_k(S)$ is just the complete $k$-cone over $k_S^{-1/2}\cdot S$,
i.e.\ the surface $S$ rescaled by the factor $k_S^{-1/2}$.
If $k_S\leq0$ (and hence $k\leq0$) then $I=\R$.
\begin{dfn}
We call $Tube_k(S)$ the
{\em complete $k$-cusp} over $S$ if $k<k_S=0$,
the {\em complete $k$-neck} if $k\leq k_S<0$,
and the {\em complete (Euclidean) cylinder} if $k=k_S=0$.
\end{dfn}

By an {\em umbilic tube} we mean a closed connected subset of a complete tube
which is a union of leaves of the natural umbilic foliation.
We will use the following terminology for different types of umbilic tubes:
A {\em standard ball}
is a truncated cone over a spherical cone surface which is homeomorphic to
the 2-sphere.
A {\em cusp} is a convex umbilic tube inside a complete cusp which is bounded
by one umbilic leaf.
A {\em neck} is a convex umbilic tube inside a complete neck
bounded by two umbilic leaves;
a neck has a totally geodesic {\em central leaf}.
A {\em cylinder} is an umbilic tube inside
a complete Euclidean cylinder
bounded by at most two totally geodesic leaves.

\subsection{Statement of the main geometric results}
\label{sec:thin}

The main result of this chapter
is the following description of the
thin part of cone 3-manifolds with lower diameter bound and cone angles
bounded away from $\pi$.
To simplify the exposition,
we will also assume a lower bound on cone angles.

\begin{thm}
[Thin part]
\label{descthin}
For $\kappa,D_0>0$ and $0<\beta<\al<\pi$
there exist constants $i=i(\kappa,\al,D_0,\beta)>0$,
$P=P(\kappa,\al,D_0)>0$
and
$\rho=\rho(\kappa,\al,D_0,\beta)>0$
such that:

Let $X$ be an orientable cone 3-manifold without boundary
which has curvature $k\in[-\kappa,0)$, cone angles
$\in[\beta,\al]$ and $diam(X)\geq D_0>0$.
Then $X$ contains a, possibly empty,
disjoint union $X^{thin}$ of submanifolds
which belong to the following list:
\begin{itemize}
\item
smooth Margulis tubes: tubular neighborhoods of closed geodesics
and smooth cusps of rank one or two,
\item
tubular neighborhoods of closed singular geodesics,
\item
umbilic tubes with one or two boundary components which are
strictly convex.
\end{itemize}
Furthermore,
the boundary of each component of $X^{thin}$ is
non-empty, strictly convex with principal curvatures $\leq P$,
and each of its (at most two) components is thick in the sense that it
contains a smooth point with injectivity radius $\geq\rho$ (measured in $X$);
each component of $X^{thin}$ contains an embedded smooth standard ball of
radius $\rho$;
all singular vertices are contained in $X^{thin}$,
and on $X-X^{thin}$ holds $r_{cone-inj}\geq i$.
\end{thm}

The proof will be given in Section \ref{sec:thinpart}.

We call $X^{thin}$ the {\em thin part} of $X$
and its components {\em thin submanifolds} or {\em Margulis tubes}.
Notice that some components of $X^{thin}$ may be
balls around singular vertices with thick links;
one may argue whether such components should be called thin as well.

We deduce two important consequences of Theorem \ref{descthin}
which we will use in the proof of the main theorem.

\begin{cor}
[Thickness]
\label{thick}
There exists
$r=r(\kappa,\al,D_0,\beta)>0$ such that:
If $X$ is as in Theorem \ref{descthin} then
$X$ is $r$-thick, i.e.\ contains an embedded smooth standard ball of radius
$r$.
\end{cor}
\begin{proof}
If $X^{thin}\neq\emptyset$, we find a thick smooth point on $\D X^{thin}$.
If $X^{thin}=\emptyset$, there are no singular vertices and the lower bounds on
$r_{cone-inj}$ and the cone angles imply thickness as well.
\end{proof}

\begin{cor}
[Finiteness]
\label{finiteness}
Let $X$ be as in Theorem \ref{descthin} and
suppose in addition that
$vol(X)<\infty$.
Then $X$ has finitely many ends and all of them are
(smooth or singular) cusps with compact cross sections. In other
words, $X$ has a compact core with horospherical boundary.
\end{cor}
\begin{proof}
According to Theorem \ref{descthin}
each thin submanifold contributes a definite quantum to the volume of $X$.
Thus $X^{thin}$ can have only finitely many components.
Finiteness of volume implies moreover
that thin submanifolds are compact or cusps with compact cross sections.

Consider a (globally minimizing) ray $r:[0,\infty)\ra X$.
There is a uniform lower bound on the volume of balls with radius $i$ and
centers outside $X^{thin}$,
where $i$ is the constant in Theorem \ref{descthin}.
Hence,
by volume reasons,
$r$ enters $X^{thin}$ after finite time.
A thin submanifold containing a ray is noncompact
and must therefore be a cusp.
We conclude that the complement of all cusp components of $X^{thin}$ is
compact,
because otherwise it would contain a ray which would
end up in yet another cusp, a contradiction.
\end{proof}

\subsection{A local Margulis lemma for incomplete manifolds}
\label{seclocmarglem}

The results in this section will be applied to the smooth part of cone
manifolds.

Let $M$ be an {\em incomplete} 3-manifold of constant negative
sectional curvature $k\in[-\kappa,0)$. Our discussion could be
carried out for arbitrary curvature sign. However, we restrict to
negative curvature for simplicity and because this is the only
case needed later.

We recall that the {\em developing map} is a local Riemannian
isometry $dev\! :\tilde M\ra \M^3_k$. It is unique up to
postcomposition with an isometry and induces the {\em holonomy}
homomorphism
$hol\!:Isom(\tilde M)\ra Isom(\M^3_k)$. The action
$\Ga:=\pi_1(M)\curvearrowright\tilde M$ of the fundamental group
by deck transformations on the universal cover transfers, via
composition with $hol$, to {\em holonomy action}
$\Ga\curvearrowright \M^3_k$. Whereas the deck action is properly
discontinuous and free, the holonomy action is in general
non-discrete.

Even though $\tilde M$ may have complicated geometry, the next
result shows that complete distance balls in $\tilde M$ are
standard; recall the definitions of various radii from Section
\ref{sec:expcutinj}.

\begin{lem}
\label{balls}
Let $\tilde x\in\tilde M$ be a lift of $x\in M$.
Then $r_{inj}(\tilde x)=r_{geod}(x)$.
\end{lem}
\begin{proof}
We have $r_{inj}(\tilde x)\leq r_{geod}(\tilde x)=r_{geod}(x)$.
The immersion $B_{r_{geod}(x)}(\tilde x)\looparrowright
\M^3_k$ into model space given by the developing map must be an
isometry onto a round ball.
Therefore also $r_{inj}(\tilde x)\geq r_{geod}(x)$.
\end{proof}

We will use these standard balls in $\tilde M$ to localize the
usual arguments in the Margulis lemma for complete manifolds of
bounded curvature.

For $\de'>0$ and for a point $y\in\tilde M$ with
$r_{geod}(y)>\de'$, we define $\Ga_y(\de')\subset\Ga$ as the
subgroup generated by all elements $\ga$ with $d(\ga \,y,y)<\de'$.
It is non-trivial if the corresponding point in $M$ has small
injectivity radius. For $r,\de>0$ and points $y\in\tilde M$ with
$r_{geod}(y)>2r+\de$ let us moreover define $A_y(r,\de)\subset\Ga$
as the subgroup generated by all elements which have displacement
$<\de$ everywhere on the closed ball $\bar B_r(y)$. The definition
is made so that, if $\de$ is small compared to $r$, then the
generators of $A_y(r,\de)$ have small rotational part. The groups
$\Ga_y$ and $A_y$ are locally semi-constant: for $z$ sufficiently
close to $y$ holds $\Ga_z\supseteq \Ga_y$ and $A_z\supseteq A_y$.
A pigeonhole argument shows that for sufficiently small
$\de'=\de'(\kappa,r,\de)>0$ holds: $A_y(r,\de)$ is non-trivial if
$\Ga_y(\de')$ is.

The standard commutator estimate yields:

\begin{prop}
\label{localmargulis} For $R>0$ there exist constants $r>>\de>0$
also depending on $\kappa$
such that for every point $y\in\tilde M$ with $r_{geod}(y)>R$ the
group $A_y(r,\de)$ is abelian.
\end{prop}
\begin{rem}
In the more general situation of variable curvature
one obtains that the groups are nilpotent.
We are using that
nilpotent subgroups of $Isom(\M^3_k)$ are abelian.
\end{rem}

We fix $R,r,\de>0$ so that \ref{localmargulis} holds.
We define the {\em thin part} $\tilde M^{thin}$ of $\tilde M$ as the
open subset of points $y$ with $r_{geod}(y)>R$ and non-trivial
$A_y(r,\de)$, and the thin part $M^{thin}$ as its projection to
$M$.

There is a natural codimension-one locally homogeneous Riemannian
{\em foliation} on the thin part.
This can be seen as follows.

Consider a point $y\in M^{thin}$, i.e.\ $r_{geod}(y)>R$ and
$A_y=A_y(r,\de)$ is non-trivial. Let $A'_y\subset Isom(\M^3_k)$ be
the image of $A_y$ under the holonomy homomorphism $hol$. Observe
that by Lemma \ref{balls}, $hol$ is injective on isometries $\phi$
with $d(\phi\, y,y)<r_{geod}(y)$, and hence $A'_y$ is still
non-trivial. The group $A_y$ is generated by small elements which
in particular preserve orientation. The classification of
isometries of $\M^3_k$ implies that $A'_y$ either preserves a
unique geodesic (axis) or the horospheres centered at a unique
point at infinity. In both cases there is a natural choice of a
connected abelian subgroup $H_y\subset Isom(\M^3_k)$ containing
$A'_y$, namely the identity component of the stabilizer of the
axis, respectively, the group of translations along the
horospheres. Moreover, there is a corresponding $A'_y$-invariant
locally homogeneous Riemannian foliation ${\cal F}_y$ of $\M^3_k$,
namely by $H_y$-orbits. The leaves are equidistant surfaces of the
axis or they are horospheres. ${\cal F}_y$ pulls back by the
developing map to a foliation of $\tilde M^{thin}$ near $y$. The
local semi-constancy of $A_y$ implies that these locally defined
foliations fit together to form a natural $\Ga$-invariant
foliation $\tilde{\cal F}$ of $\tilde M^{thin}$. There may be
one-dimensional singular leaves, namely geodesic segments in
$\tilde M$ fixed by small deck transformations; for instance,
complete $A_y$-invariant geodesics project to short closed
geodesics in $M$. $\tilde{\cal F}$ descends to a foliation ${\cal
F}$ of $M^{thin}$. Note that the regular (two-dimensional) leaves
of ${\cal F}$ are intrinsically flat and extrinsically strictly
convex.

\subsection{Near singular vertices and short closed singular geodesics}

In this section,
we make the following general assumption:

\begin{ass}
\label{ass1} We assume that $X$ is a cone 3-manifold of constant curvature
$k\in[-\kappa,\kappa]$ with cone angles $\leq\pi$ and $diam(X)\geq D_0>0$.
\end{ass}

The following result parallels Lemma~\ref{vertexthencone}:

\begin{lem}
[Thick vertex]
\label{lem:thickvertex}
For $0<d<\frac{\pi}{2}$
there exists $i=i(\kappa,D_0,d)>0$ such that:
If $v$ 
is a singular vertex with $diam(\La_vX)\leq d$,
then $r_{inj}(v)\geq i$.
\end{lem}
\begin{proof}
Since $diam(X)\geq D_0$, there exists a point $y$ with $d(y,v)\geq
D_0/2$. Let $x$ be a point in $Cut(v)$ closest to $v$. Either $x$
is the midpoint of a geodesic loop $l$ of length $2r_{inj}(v)$
based at $v$, or $x$ belongs to a singular edge with cone angle
$\pi$ and there is a (unique) minimizing geodesic segment
$s=\ol{vx}$ of length $r_{inj}(v)$ which is perpendicular to the
singular locus at $x$, cf.\ our discussion of the cut locus in
Section \ref{sec:expcutinj}. In both cases, we have a geodesic
triangle $\De(v,y,x)$ with $\angle_x(v,y)\leq\frac{\pi}{2}$. By
our assumption on the diameter of $\La_vX$ holds moreover
$\angle_v(y,x)\leq d$. Triangle comparison yields a positive lower
bound $i(\kappa,D_0,d)$ for $r_{inj}(v)=d(v,x)$.
\end{proof}

\begin{rem}
Lemma~\ref{lem:thickvertex} allows to apply the compactness
results of Section \ref{sec:compact} in many situations, for
instance to cone manifolds $X$ with a singular vertex $v$ where
the cone angles of at least two adjacent singular edges are
$\leq\pi-\eps<\pi$, since in this case $diam(\La_vX)\leq d(\eps)
<\frac{\pi}{2}$, cf.\ Lemma \ref{fatturnover}.
\end{rem}

\begin{dfn}
The {\em normal injectivity radius}
of a closed (smooth or singular) geodesic $\ga$
is the
maximal radius $r_{inj}(\ga)\in(0,\infty]$
up to which the exponential map on the normal bundle of $\ga$
is defined and is an {\em embedding},
i.e.\ for every direction $\xi$ perpendicular to $\ga$
and for every $0<l<r_{inj}(\ga)$
there exists a geodesic segment of length $l$ with initial direction $\xi$
which minimizes distance from $\ga$.
\end{dfn}

Parallel to Lemma~\ref{raysatsingedge} we have:

\begin{lem}
[Normal injectivity radius at short singular circles]
\label{injatclosedsinggeo}
For $0<\al<\pi$ there exist
$l=l(\kappa,D_0,\al)>0$ and $n=n(\kappa,D_0,\al)>0$ such that:

A singular closed geodesic $\si$ with length $\leq l$
and cone angle $\leq\al<\pi$
has normal injectivity radius $\geq n$.
\end{lem}
\begin{proof}
We will choose $l$ smaller than $\frac{D_0}{3}$ and hence can
pick a point $y$ at distance $d(y,\si)\geq\frac{D_0}{3}$ from $\si$.

Consider a minimizing segment
$\tau=\ol{wy}$ from a point $w\in\si$ to $y$.
We apply comparison to
the geodesic triangle with sides $\tau,\si,\tau$.
This can be done although the side $\si$ is of course not minimizing.
We obtain (for both angles between $\si$ and $\tau$ at $w$):
\begin{equation}
\label{anglewithlongseg}
\angle_w(\si,\tau) \geq\frac{\pi}{2}-\eps(\kappa,D_0,l)
\end{equation}
with
$\eps=\eps(\kappa,D_0,l)>0$ and $\eps\ra0$ as $l\ra0$.

We proceed as in the proof of Lemma \ref{lem:thickvertex}.
Let $x$ be a point in $Cut(\si)$ closest to $\si$.
Either $x$ is the midpoint of a segment of length $2r_{inj}(\si)$
which is perpendicular to $\si$ at both endpoints,
or $x$ belongs to a singular edge with cone angle $\pi$
and there is a minimizing segment $\ol{wx}$ of length $r_{inj}(\si)$
perpendicular to $\si $ at $w$ and to $\Si_X$ at $x$.
In both cases there exists a point $w\in\si$
and a geodesic triangle $\De=\De(w,y,x)$
with the properties:
(i) $d(w,y)\geq D_0/3$;
(ii) $d(w,x)=r_{inj}(\si)$;
(iii) $\angle_x(w,y)\leq\frac{\pi}{2}$; and
(iv) the side $\ol{wx}$ is perpendicular to $\si$.

We use property (iv)
to bound the angle of $\De$ at $w$ from above:
The link $\La_wX$ at $w$ is the metric suspension of a circle of length
$\leq\al$,
and hence
(\ref{anglewithlongseg}) implies
$\angle_w(y,x)\leq \al/2+\eps$.
By choosing $l=l(\kappa,D_0,\al)>0$ sufficiently small,
we can assure, for instance, that
(v) $\angle_w(y,x)\leq (\al+\pi)/4<\frac{\pi}{2}$.
Triangle comparison using the properties (i)-(v) yields
a positive lower bound $n(\kappa,D_0,\al)$ for $r_{inj}(\si)$.
\end{proof}

\subsection{Near embedded umbilic surfaces}
\label{secnecks}

In this section,
$X$ denotes an orientable cone 3-manifold without boundary
which has curvature $k\in[-\kappa,\kappa]$
and cone angles $\leq\pi$.
We do not need to assume a lower diameter bound.

\begin{dfn}
Suppose that $S\subset X$ is an embedded compact connected surface
such that $(S-\D S)\cap \Si_{X}$ is discrete and $\D S$ is a union of
singular edges with cone angle $\pi$.
We call the surface $S$ {\em umbilic} if
$S^{smooth}:=S-\D S-\Si_{X}$ is umbilic.
\end{dfn}
It follows that $S-\D S$ meets the singular locus
orthogonally in non-vertex singular points.
Moreover $\D S$ can be non-empty only in the totally geodesic case.

Nearby equidistant surfaces of umbilic surfaces are also umbilic.
We say that two compact connected embedded umbilic surfaces
in a cone manifold are {\em parallel}
if their union bounds an embedded umbilic tube.

\medskip
In the first part of our discussion,
we make the following assumption.
Results in the general case will be deduced afterwards.

\begin{ass}
\label{ass2}
Suppose that $S$ is {\em separating} and {\em not totally geodesic}.
\end{ass}

Since $S$ is not totally geodesic, it is two-sided.
It has a {\em convex} and a {\em concave side} defined as follows:
We say that a locally defined unit normal vector field $N$ along $S$
points to the convex side
if the principal curvature of $S$ with respect to $N$ is positive,
i.e.\ if the shape operator $DN$,
defined on tangent spaces to $S$ at smooth points,
is a positive multiple of the identity.
We call the other side of $S$ concave.

Analogously to the cut locus with respect to a point, compare our
discussion in Section \ref{sec:expcutinj}, one can define the {\em
cut locus $Cut(S)$ with respect to the umbilic surface $S$}. Let
$U(S)$ be the union of $S$ and all half-open geodesic segments
$\ga\!:[0,l)\ra X$ emanating from $S$ in orthogonal direction,
$\ga(0)\in S$, and minimizing the distance to $S$. It is an open
subset of $X$. We call the metric completion ${\cal D}(S)$ of
$U(S)$
the 
{\em Dirichlet domain} relative to $S$.
It canonically embeds into
$Tube_k(S)$
and there is a natural quotient map
\begin{equation}
\label{projdirix} \phi:{\cal D} (S)\lra X
\end{equation}
The {\em cut locus} $Cut(S)$ is defined as the complement $X-U(S)$.
Since $S$ separates $X$,
each connected component of the cut locus is either a locally finite
totally geodesic 2-complex or a point corresponding to a tip of
$Tube_k(S)$ contained in ${\cal D}(S)$;
with every tip, ${\cal D}(S)$ contains the entire component of $Tube_k(S)-S$.

The upper bound $\pi$ on cone angles implies that ${\cal D} (S)$ is
{\em convex}.

We will denote by
$X_{conv}(S)$, $Cut_{conv}(S)$, ${\cal D}_{conv}(S)$
and $\D_{conv}{\cal D}(S)$
the portions of $X$, the cut locus, Dirichlet domain and its
boundary on the convex side of $S$,
and similarly by
$X_{conc}(S)$, $Cut_{conc}(S)$, ${\cal D}_{conc}(S)$ and
$\D_{conc}{\cal D}(S)$ the portions on the concave side.

\medskip
The next two lemmas concern the component $X_{conc}(S)$ on the concave side.

\begin{lem}
\label{stronglycurvedboundarythenball}
If $S$ is spherical or horospherical,
then it bounds a standard ball respectively a cusp embedded in $X$.
\end{lem}
\begin{proof} The Dirichlet domain ${\cal D} (S)$ is convex and therefore contains the convex hull
of $S$ in $Tube_k(S)$. Since $S$ is not hyperspherical, the convex
hull fills out the whole component of $Tube_k(S)$ on the concave side
of $S$. This is a standard ball or cusp, according to whether $S$ is
spherical or horospherical, and it embeds into $X$
via the map (\ref{projdirix}).
\end{proof}

The umbilic surface $S$ can be hyperspherical only if $k<0$. In
this case we define $\rho=\rho(k,pc(S))$ as the distance from $S$
to the totally geodesic central leaf $L_{central}$ in $Tube_k(S)$.
We denote by $T$ the umbilic tube between $S$ and $L_{central}$.

\begin{lem}
\label{cutawayfromconvexbd}
(i)
If $S$ is hyperspherical
then $d(S,Cut_{conc}(S))\geq\rho$
and the map (\ref{projdirix})
is an embedding on $T-L_{central}$. It is an
embedding on $T$ if $d(S,Cut_{conc}(S))>\rho$.

(ii)
Rigidity:
If $d(S,Cut_{conc}(S))=\rho$,
then $\D_{conc} {\cal D} (S)=L_{central}$.
The map (\ref{projdirix}) restricts on $L_{central}$
to a 2-fold ramified covering over $Cut_{conc}(S)$.
The corresponding identifications on $L_{central}$
are given by
an orientation reversing isometric involution $\tau$;
its fixed point set
is either empty or a piecewise geodesic one-manifold
and maps homeomorphically onto the boundary of $Cut_{conc}(S)$
which is a union of singular edges with cone angle $\pi$.
\end{lem}
\begin{proof}
(i)
$T$ is the closed convex hull of $S$ in $Tube_k(S)$ and therefore belongs
to ${\cal D} (S)$. This implies the first part of the assertion.

(ii)
Note that as soon as ${\cal D} (S)$
contains a neighborhood of a point of $L_{central}$, then it
contains a neighborhood of the entire leaf $L_{central}$ and thus
$d(S,Cut_{conc}(S))>\rho$.
We are using here that $S$ is connected.
Therefore, if $d(S,Cut_{conc}(S))=\rho$, then
$L_{central}=\D_{conc} {\cal D} (S)$. Thus $X_{conc}(S)$  arises from
$T$ by boundary identifications on $L_{central}$, and
$L_{central}$ maps via (\ref{projdirix}) onto $Cut_{conc}(S)$.
It is clear that the identifications on $L_{central}$ arise from an isometric
involution $\tau$.
It must be orientation reversing because $X$ is orientable by assumption.
\end{proof}

\medskip\noindent
Now we investigate
the cut locus on the convex side of $S$.

Let $S_{k,P}$ be a complete umbilic surface with principal curvature
$pc(S_{k,P})=P>0$ in the smooth model space $\M^3_k$, and let $y$ be a
point on the convex side at distance $h>0$ from $S_{k,P}$. Consider the
convex hull $C$ of $S_{k,P}$ and $y$. It is rotationally symmetric, and we
define $\psi=\psi(k,P,h)\in(0,\frac{\pi}{2}]\cup\{\pi\}$ as its opening
angle, i.e.\ we set $\psi:=\pi$ if $y$ is an interior point of $C$ --
which can only happen if $k>0$ -- and define $\psi$ as the radius of the
disc $\L_yC$  otherwise.

We are interested in lower bounds for $\psi$.
Since the function $\psi(k,P,h)$ is not monotonic in all variables, we set
$\hat\psi(\kappa,P,h):=
\inf_{-\kappa\leq k\leq\kappa,0<P'\leq P,0<h'\leq h}\psi(k,P',h')$.
Then we have for all $\kappa,P>0$:
\begin{equation}
\label{openingangle}
\lim_{h\ra0}\hat\psi(\kappa,P,h)=\frac{\pi}{2}.
\end{equation}

\begin{lem}
\label{anglecutlocusminseg}
If $pc(S)\leq P$
and if $x\in Cut(S)$ with $d(x,S)\leq h$,
then the angle at $x$ between $Cut(S)$ and any shortest segment
from $x$ to $S$ is $\geq\hat\psi=\hat\psi(\kappa,P,h)$.
In particular,
the angle at $x$ between any two shortest segments to $S$ is
$\geq2\hat\psi$.
\end{lem}
\begin{proof} A shortest segment from $x$ to $S$ corresponds to a point
$\bar x\in\D {\cal D} (S)$. Let $\xi$ denote the direction at
$\bar x$ of the perpendicular to $S$. The Dirichlet domain ${\cal
D} (S)$ contains the convex hull of $S$ and $\bar x$ which in turn
contains, locally at $\bar x$, the cone over the disc of radius
$\hat\psi$ around $\xi$ in $\La_{\bar x}Tube_k(S)$. This shows the
first assertion, and the second is a direct consequence.
\end{proof}

We use Lemma~\ref{anglecutlocusminseg} to bound the number of
shortest segments from a point $x$ to $S$
and to rule out branching of the cut locus sufficiently close to $S$.
We obtain the following description of the geometry of $Cut(S)$ near $S$:

\begin{lem}
\label{cutnearsigma}
There exists $h=h(\kappa,P)>0$ with the following
property:
If $pc(S)\leq P$
and if $x\in Cut(S)$ with $d(x,S)<h$,
then there are at most two
shortest segments from $x$ to $S$.

If there are exactly two shortest segments $\tau_1$ and $\tau_2$, then
$Cut(S)$ is totally geodesic near $x$. If in addition $x$ is
singular, then $\tau_1\cup\tau_2$ forms a singular segment
orthogonal to $S$ at both endpoints and with $x$ as midpoint.

If there is only one shortest segment $\tau$, then either $\tau$ is smooth
and $x$ is an interior point of a singular edge $\si$ with cone angle
$\geq2\hat\psi$, or $x$ is a singular vertex, $\tau$ a singular segment,
and the other two singular segments $\si_1$ and $\si_2$ emanating from $x$
have cone angles $\geq2\hat\psi$. In the first case, $Cut(S)$ is near $x$
a totally geodesic half disc bounded by $\si$; in the second case it is a
sector, that is, the $k$-cone over an arc of length
$\leq\frac{\alpha(\tau)}{2}\leq\frac{\pi}{2}$
bounded by $\si_1$ and $\si_2$, where
$\alpha(\tau)$ denotes the cone angle at $\tau$.
\end{lem}

\begin{proof}
Using (\ref{openingangle}),
we choose $h>0$ sufficiently small so that
$\hat\psi(\kappa,P,h)>\frac{\pi}{3}$.
According to Lemma \ref{anglecutlocusminseg}
any two shortest segments from $x$ to $S$ have angle $>\frac{2\pi}{3}$, and
hence by Lemma~\ref{2dirbigangle} there can be at most two of them.

Regarding the second part,
the assertion is clear for smooth points $x$.
Suppose therefore that $x$ is singular
and that there are two shortest
segments $\tau_1$ and $\tau_2$ from $x$ to $S$.
Since $diam(\La_{x}X)>\frac{\pi}{2}$,
$x$ cannot be a singular vertex, cf.\ Lemma \ref{fatlink}.
Hence $x$ lies on
a singular edge $\si$ and divides it into singular segments
$\si_1$ and $\si_2$.

Note that if the metric suspension of a circle of length $\leq\pi$
contains two points with distance $>\frac{2\pi}{3}$,
then each pole of the suspension lies within distance $<\frac{\pi}{3}$
of one of the points.
Thus, after reindexing if necessary,
we have $\angle_x(\si_i,\tau_i)<\frac{\pi}{3}$.
By Lemma~\ref{anglecutlocusminseg},
the $\si_i$ cannot belong to $Cut(S)$ near $x$.
Hence $\tau_i\subset\si_i$.

Suppose now that there is just one shortest segment $\tau$
from $x$ to $S$.
If $x$ is an interior point of a singular edge $\si$ with cone angle $\beta$
then, near $x$, the cut locus is a totally geodesic half disc
bounded by $\si$.
The angle between $\tau$ and $Cut(\Si)$ at $x$
is hence $\leq\frac{\beta}{2}$,
and
Lemma~\ref{anglecutlocusminseg}
implies $\beta\geq2\hat\psi$.

We are left with the case that $x$ is a singular vertex. By
Lemma~\ref{anglecutlocusminseg}, the link $\La_xX$ has injectivity radius
$>\frac{\pi}{3}$ at the direction tangent to $\tau$, and an area
estimate implies that $\tau$ must be singular.
(A spherical turnover with cone angles $\leq\pi$ has area
$\leq\frac{1}{4}area(\S^2)$, which equals the area of a smooth
spherical disc with radius $\frac{\pi}{3}$.
Hence the direction
of $\tau$ at $x$ cannot be a smooth point of $\La_xX$.) Our
previous argument shows that the cone angles at singular points
near $x$ and not on $\tau$ are $\geq2\hat\psi$. The rest follows.
\end{proof}

\begin{cor}
\label{closecuttotgeod2}
There exists $h=h(\kappa,P)>0$ such that:
If $pc(S)\leq P$
then, up to distance $h$ from $S$,
$Cut(S)$ is a totally geodesic surface,
possibly with boundary.
\end{cor}

Our next motive is that, due to the convexity of ${\cal D} (S)$, $Cut(S)$
cannot bend away from $S$ too fast. If $S$ has small diameter, or bounded
diameter and small principal curvature, this will force the cut locus to
close up as soon as it approaches $S$ sufficiently.

\begin{lem}
\label{cutstaysclose}
For $h>0$
there exist $d_1=d_1(\kappa,P,h)>0$ and $\tilde h=\tilde h(\kappa,P,h)>0$
such that holds:

If $0<pc(S)\leq P$,
$diam(S)\leq d_1$ and $d(S,Cut_{conv}(S))\leq\tilde h$,
then every segment emanating from $S$
in the perpendicular direction to the convex side
hits the cut locus within distance $<h$.
Moreover,
$Cut_{conv}(S)$ is a compact totally geodesic surface, possibly with boundary,
which is entirely contained in the tubular neighborhood
$N_h(S)$ of radius $h$ around $S$.
\end{lem}
\begin{proof} Suppose that there exists a unit speed segment
$\tau:[0,h]\ra {\cal D} (S)$ of length $h$ emanating in the
perpendicular direction to the convex side of $S$. Moreover,
consider another such segment $\tilde\tau:[0,l]\ra {\cal D} (S)$
of length $l$ which connects $S$ to the nearest point of
$\D_{conv}{\cal D} (S)$. We then have
$d(\tau(0),\tilde\tau(0))\leq diam(S)\leq d_1$, and, due to the
convexity of ${\cal D}(S)$,
$\angle_{\tilde\tau(l)}(\tilde\tau(0),\tau(h))\leq\frac{\pi}{2}$.
The segments $\tau$ and $\tilde\tau$ are opposite sides of a
(two-dimensional) quadrangle $Q$ of constant curvature $k$
embedded in $Tube_k(S)$; the side connecting $\tau(0)$ and
$\tilde\tau(0)$ is concave with curvature $pc(S)\leq P$.
Elementary geometry in the models $\M^2_k$ implies: If
$d_1=d_1(\kappa,P,h)$ is chosen small enough, then $l$ can be
bounded below by a positive constant $\tilde h(\kappa,P,h)$.

The second part of the assertion follows from Corollary
\ref{closecuttotgeod2}.
Namely, we replace $h$ by $h':=\min(h,h(\kappa,P))$,
where $h(\kappa,P)$ is taken from Corollary \ref{closecuttotgeod2},
and adjust $d_1$ and $\tilde h$ accordingly.
\end{proof}

We need the following variant of Lemma~\ref{cutstaysclose} for umbilic
surfaces with small principal curvatures instead of small diameters:

\begin{lem}
\label{cutstaysclose2}
For $d_0>0$ and $h>0$
there exist
$P_0=P_0(\kappa,d_0,h)>0$ and $\tilde h=\tilde h(\kappa,d_0,h)>0$
such that:

If $0<pc(S)\leq P_0$,
$diam(S)\leq d_0$ and $d(S,Cut_{conv}(S))\leq\tilde h$,
then every segment emanating from $S$
in the perpendicular direction to the convex side
hits the cut locus within distance $<h$.
Moreover,
$Cut_{conv}(S)$ is a compact totally geodesic surface, possibly with boundary,
which is entirely contained in the tubular neighborhood
$N_h(S)$.
\end{lem}
\begin{proof}
The first part of the assertion is proven as for Lemma \ref{cutstaysclose}.
Note that
$P_0$ and $\tilde h$ may be chosen monotonically decreasing as $h$ decreases.
Thus,
to obtain the second part,
we may decrease $h$, if necessary, below the value $h(\kappa,P_0)$ from
Corollary \ref{closecuttotgeod2},
and then decrease $P_0$ and $\tilde h$ accordingly.
\end{proof}

%
%

We now suppose in addition that
$X$ has cone angles $\leq\al<\pi$.
The results discussed above then simplify.

If $S$ is hyperspherical, and hence $k<0$,
we obtain on the concave side:

\begin{add}
[to Lemma~\ref{cutawayfromconvexbd}]
\label{addcutawayfromconvexbd}
(i)
In the rigidity case of \ref{cutawayfromconvexbd}, $Cut_{conc}(S)$ is a
closed non-orientable totally geodesic surface,
and the natural map $S\cong L_{central}\ra Cut_{conc}(S)$ is a regular
two-fold covering.

(ii)
There exists $d_0=d_0(\kappa,\al)>0$ such that:
If $diam(S)\leq d_0$,
then the rigidity case in \ref{cutawayfromconvexbd} cannot occur,
i.e.\ the tube $T$ embeds.
\end{add}

\begin{proof}
(i)
The orientation reversing involution $\tau$ on $L_{central}$ cannot have
fixed points
because there are no singular edges with cone angle $\pi$.

(ii)
We have $diam(Cut_{conc}(S))<diam(S)$.
We apply the Gau\3-Bonnet Theorem to $Cut_{conc}(S)$ and note that,
for $d_0$ sufficiently small,
the contribution of its smooth part to the curvature integral is
a small negative number,
say $\in(\al-\pi,0)$,
and the contribution of each cone point belongs to the interval
$[2\pi-\al,2\pi)$.
Since $Cut_{conc}(S)$ is non-orientable, it must be a projective plane
and have curvature integral $2\pi$.
But with one cone point, the curvature integral would amount to $<2\pi$,
and with at least two cone points to
$>(\al-\pi)+2(2\pi-\al)>2\pi$.
We get a contradiction.
\end{proof}

On the convex side of $S$,
Lemma \ref{cutnearsigma} implies for the cut locus near $S$:

\begin{lem}
\label{closecuttotgeod}
There exists $h=h(\kappa,P,\al)>0$ such that:
If $pc(S)\leq P$
then, up to distance $h$ from $S$,
$Cut(S)$ is totally geodesic without boundary.
\end{lem}
\begin{proof}
We choose $h$ sufficiently small so that $\hat\psi(\kappa,P,h)>\al/2$,
compare (\ref{openingangle}).
This rules out in
\ref{cutnearsigma} the possibility of cut points near $S$ with a
unique minimizing segment to $S$. \end{proof}

From Lemma~\ref{cutstaysclose} on the closing up of the cut locus
near umbilic surfaces of small diameter we deduce:

\begin{lem}
\label{cutstaysaway}
There exist $d_2=d_2(\kappa,P,\al)>0$
and $\tilde h=\tilde h(\kappa,P,\al)>0$
such that:
If $pc(S)\leq P$ and
$diam(S)\leq d_2$,
then $d(S,Cut_{conv}(S))>\tilde h$.
\end{lem}

\begin{proof}
We use the constant $h=h(\kappa,P,\al)$ from Lemma \ref{closecuttotgeod}
and accordingly the constants
$d_1=d_1(\kappa,P,h)=d_1(\kappa,P,\al)$ and
$\tilde h=\tilde h(\kappa,P,h)=\tilde h(\kappa,P,\al)$
from Lemma~\ref{cutstaysclose}.

Suppose that $diam(S)\leq d_2$ and $d(S,Cut_{conv}(S))\leq\tilde h$.
Lemmas~\ref{closecuttotgeod} and \ref{cutstaysclose} imply for
$d_2\leq d_1$ that $Cut_{conv}(S)$ is a closed totally geodesic
surface contained in $N_h(S)$.
Then $\D_{conv}{\cal D} (S)$ is a closed totally geodesic
surface as well,
and it follows that $k>0$.

Since $Cut_{conv}(S)$ is non-orientable,
the Gau\3-Bonnet Theorem and the upper cone angle bound $\pi$ imply that it is
a projective plane with at most one cone point.
Hence it is an index two subcover of the complete $k$-cone of a circle of
length $\leq\al$ and has diameter
$\frac{\pi}{2\sqrt{k}}\geq\frac{\pi}{2\sqrt{\kappa}}$.
On the other hand, $diam(Cut_{conv}(S))<diam(S)+2h\leq d_2+2h$.
This yields a contradiction
if $d_2$ is chosen small enough.
\end{proof}

\medskip
We drop now our assumption \ref{ass2}
from the beginning of the section
that $S$ separates $X$ and is not totally geodesic.
We are interested in the situation when $S$ has
small diameter
and controlled principal curvature.
Our above discussion leads to
the following description of the geometry near such surfaces,
which is the main result of this section:

\begin{prop}
[Neighborhoods of umbilic surfaces with small diameter]
\label{thinpartnearsmallumbilic}
For $\kappa,P>0$ and $\al<\pi$ exists
$d=d(\kappa,P,\al)>0$ such that:

Let $X$ be an orientable cone 3-manifold without boundary which has curvature
$k\in[-\kappa,0)$ and cone angles $\leq\al$.
Suppose that $S\subset X$
is a (not necessarily separating) umbilic surface
with $0\leq pc(S)\leq P$ and $diam(S)<d$.
Then $S$ is an umbilic leaf in an embedded umbilic tube $T\subset X$
with convex boundary
and the property that
each of its at most two boundary components has diameter $d$.
\end{prop}
\begin{rem}
\label{twosided}
Note that for $d$ sufficiently small,
$S$ is either horospherical or a turnover.
This follows by applying Gau\3-Bonnet as in the proof of Addendum
\ref{addcutawayfromconvexbd}.
In particular, $S$ is always two-sided.
\end{rem}
\begin{proof}
{\em Step 1.} Suppose that $S$ \emph{separates} and is
\emph{not totally geodesic}.

We choose $d$ smaller than the constant $d_0(\kappa,\al)$
in \ref{addcutawayfromconvexbd}.
By combining \ref{stronglycurvedboundarythenball},
\ref{cutawayfromconvexbd} and \ref{addcutawayfromconvexbd},
there exists an embedded umbilic tube $T_0\subset X$ with $S$ as boundary
component and the following properties:
$T_0$ is a ball if $S$ is spherical,
a cusp if $S$ is horospherical,
and a neck if $S$ is hyperspherical.
$T_0$ has strictly convex boundary with at most two components.
Their principal curvatures are $<\sqrt{\kappa}$ if $S$ is hyperspherical
and $\leq pc(S)$ otherwise;
hence in all cases they are
$\leq P'=P'(\kappa,P):=\max(P,\sqrt{\kappa})$.
We may also assume that their diameters are $<d$.

We decrease $d$ below the constant $d_2(\kappa,P',\al)$ from
\ref{cutstaysaway}.
By applying Lemma \ref{cutstaysaway}
to the boundary components of $T_0$
and repeating this procedure finitely many times,
we obtain that $T_0$ can be enlarged to an embedded tube $T$
whose boundary components have
diameter $\geq d$ and principal curvature $<P'$.

{\em Step 2.}
Suppose now that $S$ {\em does not separate}
but still is {\em not totally geodesic}.

Consider the cyclic covering $p:\hat X\ra X$ associated
to the homomorphism $\pi_1(X)\ra\Z$ given by the oriented
intersection number with $S$.
Any connected component $\hat S$ of
$p^{-1}(S)$ is an umbilic surface isometric to $S$
which separates $\hat X$,
and our previous discussion applies.
First of all,
neither
component of $\hat X$ obtained by cutting along $\hat{S}$ is a ball or cusp,
and thus $S$ is hyperspherical.
Furthermore, by step 1,
$\hat X$ contains an embedded neck $\hat T$
with $\hat{S}$ as an umbilic leaf
and boundary components of diameter $d$.

\begin{sublem}
\label{smallumbilicsdontintersect}
There exists $d_3=d_3(\kappa,P',\al)>0$ such that:
Any two separating umbilic surfaces $S_1,S_2\subset \hat X$
with $pc(S_i)\leq P'$
and $diam(S_i)\leq d_3$ 
are disjoint or coincide.
\end{sublem}

\begin{proof}
We choose $d_3\leq \min(d_2,\tilde h)$
with the constants from \ref{cutstaysaway}.
Suppose that $S_1$ and $S_2$ are not disjoint.
Then $S_1$ is contained in
$N_{d_3}(S_2)=:Z$
which, by step 1, is an umbilic tube.

The tube $Z$, or more precisely the universal cover of $Z^{smooth}$,
develops into a layer of width $2d_3$ in model space $\M^3_k$ bounded by
two leaves $L_1$ an $L_2$ of an umbilic foliation ${\cal F}_{model}$.
The universal cover of $S_1^{smooth}$ develops along a complete umbilic
surface $U$ and leaves out at most a discrete set.
It follows that $U$ is contained in the layer.
If $d_3$ is sufficiently small, $U$ cannot bound a ball contained in the
layer, because $U$ has principal curvature $\leq P'$.
Thus $U$ separates the $L_i$.
In the case when the foliation
${\cal F}_{model}$ is not spherical, this already means that $U$ must be
one of its leaves, i.e.\ is parallel to the $L_i$.

If ${\cal F}_{model}$ is spherical, then $L_1,L_2$ and also $U$
are round spheres, and we need one more observation to see that
$U$ is concentric with the $L_i$. We consider the function
$f=d(L_1,\cdot)$ on model space. Since the development of the
universal cover of $S_1^{smooth}$ into $U$ is equivariant with
respect to its deck group, the restriction of $f$ to $U$ must have
a minimum and maximum point within distance $\leq d_3$. This
forces $U$ to be concentric with the $L_i$ if $d_3$ is small
enough. It then follows that $S_1$ and $S_2$ are parallel and thus
coincide. \end{proof}

We decrease $d$ further so that
$d\leq d_3$. 
All umbilic leaves of $\hat T$
have diameter $\leq d$ and principal curvature $<\sqrt{\kappa}\leq P'$.
Sublemma \ref{smallumbilicsdontintersect} therefore implies,
that any two translates of $\hat T$
by a non-trivial deck transformation of $\hat X\ra X$ are disjoint.
It follows that $\hat T$ projects to an embedded neck in $X$ around $S$,
and we are done in this case, too.

{\em Step 3.} Finally assume that $S$ is {\em totally geodesic}.

If $d$ is sufficiently small, our assumptions apply that $S$ is
two-sided, cf.\ Remark~\ref{twosided}, and we deduce the claim by
applying the above discussion to nearby equidistant surfaces of
$S$.
%

This completes the proof of Proposition \ref{thinpartnearsmallumbilic}.
\end{proof}

\subsection{Finding umbilic turnovers}
\label{secfindumbilic}

As in Section \ref{secnecks}, let $X$ denote an orientable cone
3-manifold without boundary which has curvature
$k\in[-\kappa,\kappa]$ and cone angles $\leq\pi$.

We are interested in conditions which imply the existence of umbilic
turnovers with small diameter and controlled principal curvature.
We will find them as cross sections to minimizing singular segments with cone
angle bounded away from $\pi$ in regions of small injectivity radius,
cf.\ our main result Proposition~\ref{umbiliccrosssec}.

We start with some observations
about the geometry near the middle of minimizing segments in $X$
which express aspects of an almost product structure.

\begin{lem}
\label{loopatleastorthogonaltolongsegment}
For $d,\eps>0$ there exists $l(\kappa,d,\eps)>0$ such that:

Let $\la$ be a (not necessarily shortest)
geodesic loop of length $\leq l$ based at $x$,
and let $\tau$ be a minimizing segment of length $\geq d$
with $x$ as initial point.
Then $\angle_x(\tau,\la)\geq\frac{\pi}{2}-\eps$.
\end{lem}

This means that the angle of $\tau$
with {\em both} initial directions of $\la$ is $\geq\frac{\pi}{2}-\eps$.

\begin{proof}
The assertion follows by applying angle comparison to the triangle with
sides $\tau,\la,\tau$.
This triangle has two minimizing sides, namely twice $\tau$,
and the non-minimizing side $\la$.
We apply comparison at the two angles adjacent to $\la$.
\end{proof}

\begin{lem}
\label{biganglebetweennonminseg}
For $0<L<\frac{\pi}{\sqrt{\kappa}}$ and $\eps>0$
exists $e=e(\kappa,L,\eps)>0$ such that:

If $q_{\pm},\,x\in X$ are points satisfying
$d(q_-,q_+)=L$ and $d(x,q_{\pm})\geq\frac{L}{4}$
and if
$\al_{\pm}=\ol{xq_{\pm}}$ are {\em not necessarily minimizing}
geodesic segments such that
\begin{equation}
\label{lengthexcess} length(\al_+)+length(\al_-)\leq d(q_-,q_+)+e,
\end{equation}
then $\angle_x(\al_+,\al_-)>\pi-\eps$.
\end{lem}

If we restrict to nonpositive curvature, we do not need an upper
bound for $L$.

\begin{proof}
{\em Step 1.}
Suppose that one of the segments $\al_{\pm}$ is minimizing.
We then can apply angle comparison to the geodesic triangle
$\De(q_-,x,q_+)$
which has $\al_{\pm}$ as two of its sides and as third side a minimizing
segment $\ol{q_-q_+}$.
Inequality
(\ref{lengthexcess}) implies that the comparison angle at $x$
is $\geq\phi(\kappa,L,e)$
with $\lim_{e\ra0}\phi(\kappa,L,e)=\pi$.
By choosing $e$ sufficiently small,
we obtain that $\angle_x(\al_+,\al_-)\geq\pi-\frac{\eps}{3}$.

{\em Step 2.}
The general case can be deduced
by considering minimizing segments $\bar\al_{\pm}$ from $x$ to $q_{\pm}$.
Then we have
$\angle_x(\bar\al_-,\bar\al_+)\geq\pi-\frac{\eps}{3}$
and
$\angle_x(\bar\al_{\pm},\al_{\mp})\geq\pi-\frac{\eps}{3}$.
It follows that
$\angle_x(\bar\al_{\pm},\al_{\pm})\leq \frac{2\eps}{3}$
and
$\angle_x(\al_+,\al_-)\geq\pi-\eps$,
as desired.
\end{proof}

We now focus our attention to {\em singular} minimizing segments
and investigate the cut locus with respect to their midpoints.

\begin{lem}
\label{cutclosetosingseg2}
For $L,\eps>0$
there exists $h=h(\kappa,L,\eps)>0$ such that:

Let $p$ be the midpoint of a minimizing
singular segment $\si=\ol{q_-q_+}$ of length $\geq L$.
Let $x$ be a point with $d(p,x)\leq h$.
Suppose that
there are at least 3 minimizing segments $\ol{px}$, or
that $x$ is singular and
there are at least 2 minimizing segments $\ol{px}$.

Then the cone angle at $\si$ is $\geq\pi-\eps$.
\end{lem}
\begin{proof}
We may assume that $L<\frac{\pi}{\sqrt{\kappa}}$,
say $L\leq\frac{\pi}{2\sqrt{\kappa}}$.
This is relevant only in the positive curvature case.

{\em Step 1.} We denote the cone angle at $\si$ by $\theta$. The
Dirichlet polyhedron ${\cal D}(p)$ associated to $p$ can be
regarded as a convex polyhedron in the model space
$\M^3_k(\theta)$ with singular axis of cone angle $\theta$, cf.\
Section \ref{sec:expcutinj}.

The minimizing segments $\si_i$ from $p$ to $x$ correspond to points
$\bar x_i$ in $\D {\cal D} (p)$. Each of the segments $\si_i$ determines a
so-called Voronoi cell $V_i$ in the link $\La_{x}X$.
By definition, $V_i$
consists of those directions at $x$ whose angle with $\si_i$ is strictly less
than the angle with all other minimizing segments from $x$ to $p$.
$V_i$ is an open convex spherical polygon.
The link $\La_{\bar x_i}{\cal D}(p)$ is canonically identified with
the closed spherical polygon $\hat V_i$ arising as the metric completion of
$V_i$.

The Dirichlet polyhedron ${\cal D} (p)$ contains at least the
subsegment $\ol{q_+q_-}$ of the singular axis, and maybe more.
Inside ${\cal D}(p)$ there are unique geodesic segments $\ol{\bar x_iq_{\pm}}$.
Note that the corresponding segments in $X$ need not be minimizing.
However, they are almost minimizing and
their initial directions $\eta_{\pm,i}\in\ol V_i$ at $x$
are almost gradient directions
for the distance functions $-d(q_{\pm},\cdot)$.
Namely,
Lemma~\ref{biganglebetweennonminseg} implies that
for any $\eps_1>0$ (a constant to be fixed later)
we have
\begin{equation}
\label{almostantip}
\angle(\eta_{+,i},\eta_{-,i})>\pi-\eps_1
\end{equation}
if $h=h(\kappa,L,\eps_1)>0$ is chosen sufficiently small.

For each pair of points $\bar x_i$ and $\bar x_j$, $i\neq j$,
there are several geodesic segments in ${\cal D} (p)$ connecting them.
They correspond to loops in $X$ with base point $x$.
At least two of the segments $\ol{\bar x_i\bar x_j}$ have length $<2h$.
Combining Lemmas
\ref{loopatleastorthogonaltolongsegment}, \ref{biganglebetweennonminseg}
and the proof of the latter one,
we obtain, after further decreasing $h$ if necessary,
that such segments $\ol{\bar x_i\bar x_j}$
have angle $\in(\frac{\pi}{2}-3\eps_1,\frac{\pi}{2}+3\eps_1)$
with both segments $\ol{\bar x_iq_{\pm}}$.
In this sense these segments $\ol{\bar x_i\bar x_j}$
are {\em almost horizontal},
where we regard the singular axis as vertical.

{\em Step 2.}
Let us consider the case that $x$ is smooth and
there exist (at least) three minimizing segments
$\si_1,\si_2,\si_3$ between $p$ and $x$.

We next construct an almost horizontal geodesic triangle $\De$
in ${\cal D}(p)$
with vertices $\bar x_i$ which winds
once around the singular axis.
Each point $\bar x_i$ lies
on a half plane, or hemisphere if $k>0$, $H_i$ in $\M^3_k(\theta)$
bounded by the singular axis.
Notice that the $H_i$ are pairwise different
because of the horizontality of the segments $\ol{\bar x_iq_{\pm}}$
and the convexity of the Dirichlet polyhedron.
Between $\bar x_i$ and $\bar x_j$
we choose the segment which does not intersect the third half plane
(hemisphere) $H_{6-i-j}$.
The resulting triangle $\De$ is contained in ${\cal D} (p)$ by convexity.
Notice that its sides are in general not minimizing.
Each side $\ol{\bar x_i\bar x_j}$ of $\De$
determines directions $\zeta_{ij}\in \ol V_i$,
respectively, $\hat\zeta_{ij}\in \hat V_i$.
We saw in step 1 that
\begin{equation}
\label{almostorth}
|\angle_x(\eta_{\pm,i},\hat\zeta_{ij})-\frac{\pi}{2}|<3\eps_1 .
\end{equation}
Let $\phi_i$ denote
the angle of the triangle $\De$ at the vertex $\bar x_i$ measured in
${\cal D} (p)$,
i.e.\ the angle between the two directions $\hat\zeta_{ij}\in\hat V_i$,
$j\neq i$.
In view of (\ref{almostantip}) and (\ref{almostorth})
we conclude that the convex spherical polygon $V_i$ almost contains a bigon
with angle $\phi_i$.
Hence, for $\eps_1$ (and accordingly $h$) sufficiently small,
$area(V_i)\geq 2\phi_i-\eps/3$.
Since $\sum_i area(V_i)\leq area (\La_xX)= 4\pi$, we obtain
\[ angle \;  sum (\De)   \leq 2\pi +\eps/2  .\]
On the other hand, $\De$ is almost horizontal in ${\cal D} (p)$.
We embed the 2-dimensional singular model $\M^2_k(\theta)$ as cross section of
$\M^3_k(\theta)$
so that it contains, say, $\bar x_1$.
Nearest point projection of $\De$ to $\M^2_k(\theta)$
almost preserves angles.
Due to Gau\3-Bonnet, horizontal triangles with small diameter have angle
sum $\simeq\pi+(2\pi-\theta)$, and therefore
\[  angle \; sum (\De)  \geq  3\pi -\theta -\eps/2 \]
if $h$ is sufficiently small.
It follows that
$\theta\geq\pi-\eps$ as claimed.

{\em Step 3.} The argument is analogous in the case when $x$ is
singular and there are 2 minimizing segments between $p$ and $x$:
$\De$ becomes an almost horizontal bigon winding once around the
singular axis; such bigons have angle sum $\simeq 2\pi-\theta$,
i.e.\ $\geq2\pi-\theta-\frac{\eps}{2}$ for $h$ sufficiently small; on the
other hand, since $area(\La_xX)\leq2\pi$, the angle sum must be
$\leq\pi+\frac{\eps}{2}$; therefore $\theta\geq\pi-\eps$, as claimed.
\end{proof}

\begin{cor}
\label{cutnearsingseg}
For $L>0$ and $0<\al<\pi$
exists $h=h(\kappa,L,\al)>0$ such that:

If $p$ is the midpoint of a minimizing singular segment of length
$\geq L$ and cone angle $\leq\al$, then the ball $B_h(p)$ contains
no point $x$ with at least 3 minimizing segments between $x$ and
$p$, and no singular point $x$ with at least 2 minimizing segments
between $x$ and $p$. The intersection $Cut(p)\cap B_h(p)$ is a
totally geodesic surface whose boundary (if non-empty) is geodesic
and consists of singular segments.
\end{cor}

Notice that there are no singular vertices close to $p$,
because the links at vertices have diameter $\leq\frac{\pi}{2}$
(Lemma \ref{fatlink})
whereas the links at points near $p$ have almost diameter $\pi$
(Lemma \ref{biganglebetweennonminseg}).

We come to the main result of this section.

\begin{prop}
[Umbilic cross sections]
\label{umbiliccrosssec}
For $L,d>0$ and $0<\al<\pi$
there is $i=i(\kappa,\al,L,d)>0$ such that:

Let $p$ be the midpoint of a minimizing singular segment $\si$ of length
$\geq L$ and cone angle $\leq\al$,
and assume that $r_{inj}(p)<i$.
Then there exists an umbilic turnover
$S$ through $p$
with $diam(S)\leq d$.
\end{prop}
\begin{rem}
Lemma \ref{stronglycurvedboundarythenball} implies
that the principal curvature of the cross section $S$
is bounded in terms of $\kappa$ and $L$.
Namely,
if $pc(S)$ were too large,
then $S$ would bound a singular ball of radius $<L/2$.
\end{rem}

\begin{proof}
Let $q_{\pm}$ denote the endpoints of $\si$ and $\theta$ its cone angle.
We study the Dirichlet polyhedron ${\cal D}(p)$
which we regard as a convex polyhedron in $\M^3_k(\theta)$.

{\em Step 1:
The edges of ${\cal D}(p)$ are almost vertical.}
By \ref{cutnearsingseg},
near $p$
the interior points on faces of ${\cal D}(p)$ correspond to smooth cut
points, and the points on boundary edges correspond to singular cut points.
(We must also allow degenerate edges with dihedral angle $\pi$.)
Let $\bar x$
be a point on a boundary edge $\bar\ga$ of ${\cal D}(p)$
with $d(p,\bar x)<h$,
$h$ as in \ref{cutnearsingseg},
and let $x$ be the corresponding point in $X$
which lies on a singular edge $\ga$.
The segments $\ol{xq_{\pm}}$ in $X$ corresponding to the segments
$\ol{\bar xq_{\pm}}$ have angle $\simeq\pi$ by \ref{biganglebetweennonminseg}.
This forces their directions at $x$ to be close to the singular poles of
$\La_xX$.
This means that $\bar\ga$ is almost vertical.

{\em Step 2:
The cross sections of ${\cal D}(p)$ are small.}
By assumption there exists a cut point $y\in Cut(p)$ with $d(p,y)<i$.
Let $\bar y$ be a corresponding point in $\D {\cal D}(p)$.
The cross section $C_{\bar y}$ of ${\cal D}(p)$ through $\bar y$ and
perpendicular to $\si$
is a convex polygon with cone point of angle $\theta$.
Observe that, since $\theta<\pi$,
the circumradius and the inradius of the polygon control each other.
Hence,
for $i>0$ sufficiently small,
the polygon $C_{\bar y}$ has small diameter $<<h$.
Notice that $C_{\bar y}$ has at least one vertex.

{\em Step 3:
The cross sections are bigons.}
Our discussion implies that, near $p$,
cross sections of ${\cal D}(p)$ are compact convex polygons and
$\D{\cal D}(p)$ is a union of finitely many geodesic strips which are almost
vertical.
According to \ref{cutnearsingseg},
the boundary identifications on $\D{\cal D}(p)$ inside $B_h(p)$
are given by an involutive isometry
$\iota:\D{\cal D}(p)\cap B_h(p)\ra\D{\cal D}(p)\cap B_h(p)$
which fixes the boundary edges.
It follows that there are exactly two strips which are exchanged by $\iota$.

{\em Step 4:
${\cal D}(p)$ is rigid.}
Let $v_1$ and $v_2$ be the vertices of the cross section $C_{\bar y}$,
$\tau$ and $\tau'$ its sides,
and let $\si_i$ be the edge of $\D{\cal D}(p)$ through $v_i$.
We denote by $H_i$ the half-plane (respectively, hemisphere if $k>0$)
in $\M^3_k(\theta)$ bounded by the singular axis and containing
$v_i$.
Since $\tau$ and $\tau'$ are exchanged by the isometry $\iota$,
they must have the same length it follows that $H_1$
and $H_2$ are opposite to each other in the sense that they meet at $\si$
with angle $\frac{\theta}{2}$.
Furthermore,
each edge $\si_i$ has equal angles with $\tau$ and $\tau'$
and thus $\si_i\subset H_i$.
It follows that $\iota$ extends to an
isometry of the whole model space $\M^3_k(\theta)$,
namely to the reflection at $H_1\cup H_2$.

{\em Step 5:
Conclusion.}
Consider the unique umbilic foliation ${\cal F}$ of $\M^3_k(\theta)$
orthogonal to $\si$ and $\si_1$.
It is also orthogonal to the boundary strips of ${\cal D}(p)$ and hence
to $\si_2$.
Let $L$ denote the leaf of ${\cal F}$ through $p$.
Then $L\cap{\cal D}(p)$ projects to an umbilic turnover $S$ in $X$.
If $i$ is chosen sufficiently small, we have $diam(S)<d$.
\end{proof}

\subsection{Proof of Theorem \ref{descthin}: Analysis of the thin part}
\label{sec:thinpart}

In this section we combine the previous results to analyze the thin part of
$X$.
The proof of Theorem \ref{descthin} is organized in five steps:

\medskip
{\em Step 1: Around singular vertices.}
Let $v\in X$ be a singular
vertex. The space of directions $\La_vX$ then has diameter
$\leq D(\al)<\frac{\pi}{2}$,
cf.\ Lemma \ref{fatturnover},
and by the Thick Vertex Lemma \ref{lem:thickvertex}, $v$ is the center of an
embedded (closed) standard ball with radius $r_1(\kappa,\al,D_0)>0$. To
make the balls around the various vertices disjoint, we define
$B_v$ as the closed ball of radius $\frac{r_1}{2}$ centered at $v$. The umbilic
boundary spheres $\D B_v$ are convex and have principal curvature
$\leq P_1(\kappa,r_1)=P_1(\kappa,\al,D_0)>\sqrt{\kappa}$.

{\em Step 2: Organizing small umbilic turnovers.} Let us now
consider the umbilic turnovers $S \subset X$ with $pc(S)<P_1$.
If $d_2=d_2(\kappa,\al,P_1)=d_2(\kappa,\al,D_0)>0$ is chosen small enough,
then according to Proposition \ref{thinpartnearsmallumbilic}
any such turnover $S$ with diameter $<d_2$ is a leaf in the natural foliation
of an embedded umbilic tube $T_{S}\subset X$.
Moreover,
$T_{S}$ has one or two boundary components
which are strictly convex with diameter $d_2$.

The argument used
to prove \ref{smallumbilicsdontintersect} shows that,
after decreasing $d_2$ sufficiently,
any two turnovers $S$ in consideration are either disjoint or coincide.
The same holds then for the tubes $T_{S}$.
It shows as well
that the $T_{S}$ are disjoint from the balls $B_v$
with $diam(\D B_v)>d_2$. On the other hand, the singular balls $B_v$ with
$diam(\D B_v)\leq d_2$ are contained in a tube $T_{S}$. In the
following, we forget about the balls $B_v$ contained in tubes
$T_{S}$. Denote by $V_1$ the union of the remaining balls $B_v$
and the tubes $T_{S}$; we saw that they are pairwise disjoint.

Notice that
for umbilic turnovers as considered here,
i.e.\ with cone angles $\leq\al$, controlled principal curvature $\leq P_1$
and small upper diameter bound $d_2$,
diameter and thickness control each other.
This is seen as in the proof of Lemma \ref{thicklink};
the lower bound on the cone angles follows from Gau\3-Bonnet and the fact
that the Gau\3 curvature is bounded.
As a consequence,
there is a lower bound for the thickness of the components of $\D V_1$.
We also get a lower bound for the thickness of the components of $V_1$ by
leaving out ``short'' umbilic necks, i.e.\ necks with central leaves of
diameter, say, $>\frac{d_2}{2}$.

{\em Step 3: Around short closed singular geodesics.}
We choose $l_1=l_1(\kappa,\al,D_0)>0$ small enough
so that \ref{injatclosedsinggeo} implies
that the normal injectivity radius of closed singular geodesics $\ga$
with period $\leq 2l_1$
is $>2n_1(\kappa,\al,D_0)>0$.
The closed tubular $n_1$-neighborhoods $\bar N_{n_1}(\ga)$ around these
geodesics are then pairwise disjoint,
and we denote their union by $V_2$.
Note that the injectivity radii of their boundaries
are everywhere $\leq i_1(\kappa,l_1,n_1)$
with $\lim_{l_1\to0}i_1(\kappa,l_1,n_1)=0$.

By choosing $l_1$ sufficiently small, we may achieve that
\[ V_1\cap V_2=\emptyset \]
This can be seen as follows:
The singular closed geodesics of period $\leq2l_1$ lie outside $V_1$.
If $V_1\cap V_2$ were non-empty,
then an embedded smooth 2-torus $T$ (equidistant to a
short singular closed geodesic) with controlled
principal curvatures and small area
would intersect (touch) an umbilic turnover $S'$ with controlled
principal curvature and lower diameter bound.
Moreover, the distance between $T$ and $\Si_X$ is bounded away from zero
so that, from the thickness of $S'$, we get a lower bound on the injectivity
radius at the touching point.
This contradicts the thinness of the torus $T$.

{\em Step 4: Bounding the injectivity radius on the rest of the
singular locus.} First we show that, in the spirit of
Lemma~\ref{raysatsingedge},
a singular edge either closes up with short period or
minimizes up to a certain length.

\begin{sublem}
\label{findsingseg}
There exists $l=l(\kappa,\al,D_0)>0$ such that
for every $l'\leq l$ holds:

If the singular edge $\si$
does not close up with period $\leq2l'$
then, for every point $x\in\si$,
there is a {\em minimizing} subsegment of
length $>l'$ with $x$ as initial point.
\end{sublem}
\begin{proof} Let $\si_1=\ol{xy_1}$ and $\si_2=\ol{xy_2}$ be the maximal
minimizing subsegments of $\si$ emanating from $x$ in the two
antipodal singular directions. Suppose that both have length $\leq
l'$. Due to our diameter assumption, there exists a
segment $\tau=\ol{xz}$ of length $\geq D_0/2$ starting in $x$.

We regard $\si_1,\si_2,\tau$ also as segments in the Dirichlet polyhedron
${\cal D}(x)$ and denote their respective endpoints on $\D {\cal D}(x)$ by
$\bar y_1,\bar y_2,\bar z$.

Since ${\cal D}(x)$ has a singular axis with cone angle $\leq\al$,
the convex hull $C(\bar z)$ of $\bar z$ in ${\cal D}(x)$
is a totally geodesic disc
which intersects $\si$ orthogonally in a cone point $c$
and which has geodesic boundary with a corner at $\bar z$.
The convexity of the Dirichlet polyhedron implies that
$\angle_{\bar y_i}(x,\bar z)\leq\frac{\pi}{2}$.
By considering the geodesic triangle $\De(x,\bar y_i,\bar z)$ in ${\cal D}(x)$
it follows that
$\angle_x(\tau,\si_i)\geq\phi(\kappa,D_0,l)$ with $\phi\ra\frac{\pi}{2}$ as
$l\ra0$.
We obtain furthermore that
$C(\bar z)$ contains a disc of radius $r(\kappa,\al,D_0,l)>0$
centered at $c$.

We investigate $\D{\cal D}(x)$ near the singular axis. For a point
$\bar w\in\D{\cal D}(x)$ near $\bar y_i$, say with $d(x,\bar
w)\leq 2l$, we consider the convex hull of $\bar w$ and $\bar z$.
It is contained in ${\cal D}(x)$, and we obtain that the space of
directions $\La_{\bar w}{\cal D}(x)$ contains a standard disc of
radius $\rho(\kappa,\al,D_0,l)>0$. Moreover,
$\lim_{l\ra0}\rho(\kappa,\al,D_0,l)=\frac{\pi}{2}$. Note hereby
that $r(\kappa,\al,D_0,l)$ increases as $l\ra0$ and in particular
remains uniformly bounded below by a positive constant. We choose
$l=l(\kappa,\al,D_0)$ small enough so that
$\rho(\kappa,\al,D_0,l)>\al/2$. Then the metric suspensions of
circles of length $\leq\al$, and hence the links of singular
points in $X$ cannot contain a {\em smooth} $\rho$-disc. It
follows that $\bar w$ cannot project to a singular point in
$Cut(x)$ unless it is singular itself, i.e.\ coincides with a
point $\bar y_i$. On the other hand, the points $y_i$ are
singular, and thus the boundary identifications on $\D {\cal
D}(x)$ can identify the points $\bar y_i$ at most with each other
and with no other points in $\D {\cal D}(x)$. Since all cut points
near the $y_i$ are smooth implies that the $y_i$ cannot be
singular vertices. Hence $\bar y_1$ and $\bar y_2$ have to be
identified, and $\si$ closes up with period $\leq2l'$. Not that
the argument also shows that $Cut(x)$ is totally geodesic near
$y_1=y_2$.
\end{proof}

We further decrease $l_1$ until $l_1 \leq l$
with the constant $l$ of Sublemma \ref{findsingseg}.
This amounts to removing from $V_2$ some of its components.

We then have that
every singular point $x$ outside $V_1\cup V_2$
is the initial point of a minimizing singular segment $\si$
of length $\geq l_1$.
Let $m$ be the midpoint of this segment.
If there were an umbilic turnover $S$ through $m$ and $\perp\si$,
then $pc(S)< P_1$ and $diam(S)\geq d_2$
because $m\not\in V_1$.
Proposition \ref{umbiliccrosssec} implies
a lower bound $i_2=i_2(\kappa,\al,l_1,d_2)=i_2(\kappa,\al,D_0)>0$
for $r_{inj}(m)$.

We now use the {\em lower} bound $\beta$ on cone angles
to control $r_{inj}(x)$ in terms of $r_{inj}(m)$.
There is a {\em smooth} standard ball of radius
$\geq r_2(\kappa,\al,D_0,\beta)>0$
embedded in the singular standard ball $B_{i_2}(m)$.
Proposition \ref{decayofinj} implies
a lower bound for $r_{cone-inj}(x)$.
Since $x\not\in V_1$,
this yields a lower bound
$i_3(\kappa,\al,D_0,\beta)>0$
for $r_{inj}(x)$.

After feeding the baby, if necessary, we choose
$r_3(\kappa,\al,D_0,\beta)$ with $0<r_3<\min(\frac{i_3}{3},n_1)$ and
define $V_3$ as the closure of the union of all balls of radius $r_3$
centered at singular points outside $V_1\cup V_2$. Since $r_3<n_1$, we
have $V_2\cap V_3=\emptyset$ because $N_{n_1}(\D V_2)$ contains no
singular points, and so
\[ (V_1\cup V_3)\cap V_2=\emptyset .\]
Our construction yields that
(i)
$V_1\cup V_2\cup V_3$ contains the tubular $\frac{r_3}{2}$-neighborhood
of $\Si_{X}$, and
(ii)
there is a lower bound for $r_{inj}$ on $\D(V_1\cup V_3)$.
On the other hand, $\D V_2$ can become arbitrarily thin.

{\em Step 5: Foliating the thin part away from the singular
locus.} We have $r_{geod}\geq\frac{r_3}{2}$ on the complement $Y$
of $V_1\cup V_2\cup V_3$. The local Margulis lemma, cf.\ the
discussion in Section \ref{seclocmarglem}, implies that there
exist constants $i_4=i_4(\kappa,r_3)=i_4(\kappa,\al,D_0,\beta)>0$,
$l_2=l_2(\kappa,r_3)=l_2(\kappa,\al,D_0,\beta)>0$, $l_2<<i_4$, and
an open subset $Y^{thin}\subseteq Y$ carrying a natural locally
homogeneous codimension-one foliation ${\cal F}$, possibly with
singular one-dimensional leaves, which enjoy the following
properties:
(i) $\{r_{inj}<i_4\}\cap Y\subset Y^{thin}$;
(ii) ${\cal F}$ is locally equivalent to a foliation of model space $\M^3_k$
by equidistant surfaces of a geodesic or a horosphere;
(iii) intrinsically,
the regular leaves of ${\cal F}$ are flat,
and they admit foliations by parallel geodesics of length $<l_2$
(without singular leaves since $X$ is orientable).
Note that these one-dimensional foliations need not be unique on compact leaves.
Since $r_{inj}$ is bounded below on $\D(V_1\cup V_3)$,
we can arrange by choosing the constants $i_4,l_2$ sufficiently small that
\[ Y^{thin} \cap (V_1\cup V_3) =\emptyset. \]
It can happen that $Y^{thin}$ intersects a component $C$ of $V_2$.
Then $\D C\subset\ol{Y^{thin}}$
and the natural foliation of $C$ by equidistant surfaces around the
singular core geodesic
extends ${\cal F}$.
We denote by $V_2'$ the union of all these components $C$ of $V_2$
and conclude that ${\cal F}$ extends to a foliation $\hat {\cal F}$ of
$Y^{thin}\cup \bar V_2'=:\hat Y^{thin}$.
Notice that the completeness of $X$ implies that the regular leaves of $\hat {\cal F}$ are
{\em complete},
since they stay away from $\D(V_1\cup V_3)$.

Our discussion implies that
the connected components of $\hat Y^{thin}$
are smooth cusps of rank one or two (i.e.\ quotients of horoballs)
or tubular neighborhoods of closed (smooth or singular) geodesics.
The singular leaves of ${\cal F}$ are short closed smooth geodesics.
The boundary of $\hat Y^{thin}$ is a union of complete leaves
and we have $r_{inj}\geq i_4$ on $\D \hat Y^{thin}$.
Each component of $\hat Y^{thin}$ has non-empty boundary,
because the leaves of $\hat{\cal F}$ are strictly convex
and $r_{inj}$ increases towards the convex side.

We define $X^{thin}$ as $V_1\cup V_2\cup\hat Y^{thin}$ with those components
omitted which are tubular neighborhoods of closed smooth geodesics with
length, say, $>i_4$.
We already removed short umbilic necks earlier.
All other components are uniformly thick.

This concludes the proof of Theorem \ref{descthin}.

\subsection{Totally geodesic boundary}
\label{sec:totgeod}

In this section we allow totally geodesic boundary for $X$ and
apply the discussion in Section \ref{secnecks} to investigate the
geometry near boundary components.

\begin{prop}
[$I$-bundle]
\label{Ibundlethm}
For $d_0>0$ there exists $\tilde h=\tilde h(\kappa,d_0)>0$ such that:

Let $X$ be a cone 3-manifold
of curvature $k\in[-\kappa,0]$
with totally geodesic boundary
and cone angles $\leq\pi$.
Suppose that $\D X$ contains a connected component $S$ with
$diam(S)\leq d_0$ and $d(S,Cut(S))\leq\tilde h$.
Then $\D X=S$,
$Cut(S)$ is a compact totally geodesic surface,
possibly with boundary,
$\D Cut(S)$ is a union of singular edges,
and $X$ carries a natural structure as a singular bundle
$X \ra Cut(S)$
with fiber a compact interval.
\end{prop}
\begin{proof}
Let us consider the cut locus $Cut(\D X)$ with respect to the entire boundary.
Note that $d(S,Cut(\D X))\leq d(S,Cut(S))$.
The discussion of Section \ref{secnecks} applies
since we can replace boundary components of $X$ by nearby equidistant umbilic
surfaces.
We choose $h>0$ arbitrarily
and then $\tilde h$ as the constant provided by Lemma~\ref{cutstaysclose2}.
Combining Lemma~\ref{cutstaysclose2}
and the description of the cut locus in Lemma \ref{cutnearsigma}
we obtain the assertion.
The bundle structure is given as follows:
The fiber over $x\in Cut(S)$ is the union of the (one or two) shortest
segments from $x$ to $S$.
\end{proof}

Note that the boundary identifications on ${\cal D}(S)$ correspond,
via the nearest point projection $\D{\cal D}(S)\ra S$,
to an isometric involution $\iota$ on $S$.

If $k<0$,
cut points at maximal distance from $S$ must be corners of $\D Cut(S)$,
and in particular $\D Cut(S)\neq\emptyset$ is a non-empty union of singular
edges.

In the case that $k\leq0$,
$X$ is orientable and
$S$ is a turnover,
$\iota$ reverses orientation and two possibilities can occur:
Either $\iota$ fixes all three
cone points and $Cut(S)$ is a triangle,
or $\iota$ fixes one cone point and
exchanges the other two.
In the latter case, $Cut(S)$ is a disc with one cone point and one corner.

\medskip
In subsection~\ref{sec:deg} we will need the
following version of Proposition \ref{Ibundlethm}.

\begin{cor}
\label{cor:Ibundle}
Suppose that $k<0$.
Then for $0<\beta<\pi$ there exists $\rho=\rho(k,\beta)>0$ such that:
Let $X$ be a cone manifold of curvature $k$ with totally geodesic boundary
and cone angles $\in[\beta,\pi]$.
Suppose that $\D X$ contains a turnover $S$.
If $X$ is $\rho$-thin,
i.e.\ contains no embedded smooth standard ball of radius $\rho$,
then the conclusion of Proposition \ref{Ibundlethm} holds.

Moreover, given $\eps>0$, there exists
$\rho_1=\rho_1(k,\beta,\eps)>0$ such that if $X$ is $\rho_1$-thin
then  the cone angles at the singular edges in the boundary of
$Cut(S)$ are $\geq\pi-\eps$.
\end{cor}
\begin{proof}
The lower bound $\beta$ on cone angles
yields an upper bound $d_0=d_0(k,\beta)$ for the diameter
of the turnover $S$.
We are done by Proposition~\ref{Ibundlethm}
if $d(S,Cut(S))\leq\tilde h(-k,d_0)$.
We suppose therefore that $d(S,Cut(S))>\tilde h(-k,d_0)$.

We certainly have $Cut(S)\neq\emptyset$ because $X$ is $\rho$-thin
and denote $n:=d(S,Cut(S))$.
We consider the family of embedded umbilic surfaces $S_r$, $0<r<n$,
equidistant to $S$ with distance $r$ from $S$.
They have uniformly bounded principal curvatures
$pc(S_r)<\sqrt{-k}$,
and the ratio $diam(S_r)/diam(S)$
is a function of $k$ and $r$, monotonically increasing in $r$.

We choose $h$ smaller than the constant $h(-k,\sqrt{-k})$ given in
Corollary~\ref{closecuttotgeod2},
and let
$d_1=d_1(-k,\sqrt{-k},h)$ and $\tilde h=\tilde h(-k,\sqrt{-k},h)$
be the constants from Lemma~\ref{cutstaysclose}.
By choosing $\rho$ sufficiently small
we assure that the surfaces $S_r$
are uniformly thin and, in view of the lower bound on cone angles,
have diameter $<d_1$.
Lemma~\ref{cutstaysclose}, applied to $S_r$ for $r\ra n$,
then implies the conclusion of Proposition \ref{Ibundlethm}.
The additional assertion regarding the cone angles at $\D Cut(S)$
follows from Lemma~\ref{cutnearsigma}
by choosing the constant $h$ sufficiently small so that
$\hat\psi(-k,\sqrt{-k},h)>\frac\pi2-\frac\eps2$,
cf.\ (\ref{openingangle}).
\end{proof}

\section{Proof of the main theorem}
\label{sec:main}

Let $\OO$ be a compact connected orientable small $3$-orbifold
with non-empty ramification locus $\Sigma$.

The singular locus $\Sigma$ is a trivalent graph properly embedded
in $\vert \OO\vert$. Let $\Sigma^{(0)}$ denote the set of vertices
of $\Sigma$ and $\Sigma^{(1)}=\Sigma-\Sigma^{(0)}$ the union of
open edges.
We regard circle components of the singular locus as edges which close up.

We consider the manifold $M=\vert \OO\vert-\Sigma^{(1)}-{\mathcal
N}(\Sigma^{(0)})$, i.e. we remove the edges of
$\Sigma$ and an open ball neighborhood of each vertex. The manifold $M$
is noncompact, with boundary $\partial M = (\partial
\NN(\Sigma^{(0)})\cup \partial \OO) -\Sigma^{(1)}$ a finite
collection of thrice punctured spheres.

\subsection{Reduction  to the case when the smooth part is hyperbolic}
\label{sec:red2hyp}

The following proposition allows to reduce the proof of the main
theorem
to the case where $M$ admits a complete hyperbolic structure with
finite volume and totally geodesic boundary.

\begin{prop}\label{prop:eitherM} Either the manifold $M$ has a complete
hyperbolic structure with finite volume and with totally geodesic
boundary, or $\OO$ admits a Seifert fibration, an $I$-bundle structure, or
a spherical structure (i.e. a quotient of $\mathbb S^3$ or $B^3$ by an
orthogonal action).
\end{prop}

\begin{proof} Let $\overline M= \OO-\NN(\Sigma)$ be a compact core of
$M$. The boundary $\partial \overline M$ is the union of compact pairs of
pants (which are a compact core of $\partial M$) together with a
collection $P \subset\partial\overline M$ of tori and annuli,
corresponding to the boundary of a neighborhood of edges  in $\Sigma$.

\begin{lem}\label{lem:pared}
    Either $\overline M$ is Seifert fibred or $(\overline M, P)$ is
    an atoroidal pared manifold.
\end{lem}
\medskip

We recall that an \emph{atoroidal pared manifold} is a pair
$(\overline M, P)$ such that:
\begin{itemize}
    \item[--] $\overline M$ is a compact orientable irreducible 3-manifold.
    \item[--] $P\subset\partial \overline M$ is a disjoint union of
incompressible tori and annuli such that  no
    two components of $P$ are isotopic in $\overline M$.
    \item[--] $\overline M$ is homotopically atoroidal and $P$ contains all
torus components of
    $\partial \overline M$.
    \item[--] There is no essential annulus $(A,\partial A)\subset
(\overline M, P)$.
\end{itemize}

We remark that with this definition, an atoroidal pared manifold
is never Seifert fibred.

\begin{proof} The manifold $\overline M$ is irreducible and topologically
atoroidal because so is $\OO$. Assuming that $\overline M$ is not
Seifert fibred, then $\overline M$ is homotopically atoroidal, and
we prove that $(\overline M, P)$ is an atoroidal pared manifold.
First we show that $P$ is incompressible in $\overline M$. A
compressible annulus in $P$ would give a teardrop in $\OO$,
contradicting irreducibility of $\OO$. If a torus component of $P$
was compressible, then the irreducibility of $\overline M$ would
imply that $\overline M$ is a solid torus, hence Seifert fibred.
It only remains to check that the pair $(\overline M, P)$ is
anannular. Let $(A,\partial A)\subset (\overline M, P)$ be an
essential annulus; we distinguish three cases according to whether
$\partial A$ is contained in a) torus components of $P$, b)
annulus components of $P$, or c) a torus and an annulus of $P$. In
the first  case, a classical argument using the atoroidality of
$\overline M$ implies that $\overline M$ is Seifert fibred
\cite[Lemma 7]{BSOne}. In case b), adding two meridian discal
orbifolds to $A$ along $\partial A$ would give a bad or an
essential spherical 2-suborbifold, contradicting the
irreducibility of $\OO$. Case c) reduces to case b), by
considering the  essential annulus obtained by gluing two parallel
copies of $A$ with the annulus $P_0-\NN(\partial A)$, where $P_0$
is the torus component of $P$ that meets $\partial A$. \end{proof}

\begin{proof}[End of the proof of Proposition~\ref{prop:eitherM}]  We consider
both possibilities of Lemma~\ref{lem:pared}. When $\overline M$ is Seifert
fibred, then the fibration of $\overline M$ extends to a fibration of the
orbifold $\OO$ by adding the components of $\Sigma$ as fibers, because
$\OO$ is irreducible.

When $(\overline M,P)$ is an atoroidal pared 3-manifold, since $\overline M$
is Haken, by Thurston's hyperbolization theorem for atoroidal Haken pared
$3$-manifolds (cf. \cite{ThuBull,Thu3,Thu4,Thu5}, \cite{McMOne},
\cite{Kap}, \cite{MB} \cite{OtaOne,OtaTwo}), the interior of $M$
admits a complete hyperbolic structure with parabolic locus $P$.
The convex core of this metric may have dimension two or three. If
it has dimension three, then this gives a hyperbolic metric on $M$
with totally geodesic boundary and cusp ends, because the boundary
is a union of three times punctured spheres, and therefore the
Teichm\"uller space of $\partial M$ is a point.

If the convex core has dimension two, then $M$ is an $I$-bundle. Since
$\partial M$ is a union of three times punctured spheres, and such a
punctured spheres do not have free involutions, it follows that $M$ is a
product of the interval with a three times punctured sphere. Hence there
are three possibilities. In the first case $\partial \OO=\emptyset$ and
$\OO$ is the suspension of a turnover. This turnover must be spherical and
it is clear that $\OO$ is spherical itself. If $\partial \OO$ has
precisely one component, then $\OO$ is a standard quotient of a ball
(hence spherical).  Finally, the last case happens when $\partial\OO$ has
two components. In this case $\OO$ is an $I$-bundle over a turnover, the
turnover is Euclidean or hyperbolic, and $\OO$ is also Euclidean or
hyperbolic.
\end{proof}

\subsection{Deformations of hyperbolic cone structures}\label{sec:deform}

From now on we assume that the manifold $M$ admits a complete
hyperbolic structure of finite volume with totally geodesic
boundary.

Starting with the hyperbolic metric on $M$,
we define in this section a deformation space of
hyperbolic cone structures on $\OO$ and prove an openness property.

\begin{dfn}
\label{dfn:conestruc}
A {\em hyperbolic cone structure} on $\OO$ is a hyperbolic cone
3-manifold $X$ with totally geodesic boundary
together with an embedding  $i \! : (X,\Sigma_X) \hookrightarrow
(\vert \OO\vert,\Sigma)$ such that $|\OO| - X$ is a (possibly empty)
collection of vertices, open ball neighborhoods of vertices and Euclidean
2-orbifolds in $\partial\OO$.

If $i$ is a homeomorphism, we call the cone structure
a {\em cone metric} on $\OO$.
\end{dfn}

The choice of a hyperbolic cone structure on $\OO$ assigns cone angles to the
edges of $\Si$.
Usually we will consider situations when the cone angles are less than or
equal to the orbifold angles.
In order to be able to work with small cone angles,
we allow in our definition of cone structure
small deviations between the topologies of $\OO$ and $X$,
i.e.\ we do not require $i$ to be a homeomorphism.
Namely, observe that in a
cone manifold the sum of the cone angles of the
singular edges adjacent to a vertex of $\Sigma$ is $>2 \pi$,
because the link of a vertex is a spherical turnover.
Thus,
if $v$ is a vertex component of $|\OO|-X$, we require the hyperbolic
structure on the punctured neighborhood of $v$ to be a cusp;
the cone angle sum for the singular edges adjacent to $v$ then equals $2\pi$.
If $|\OO|-X$ contains the open ball neighborhood $B_v$ of a vertex $v$,
we require the boundary turnover $\D B_v$ to be totally geodesic;
the cone angle sum is $<2\pi$ in this case.
If $S$ is a Euclidean 2-orbifold in $\partial\OO-X$
then we request that the hyperbolic structure near $S$ is also a cusp.

Notice that all boundary components of $X$ are turnovers because $\OO$ is
small.

There are analogous definitions for {\em Euclidean} and
{\em spherical} cone structures on $\OO$,
however in these cases it will be enough for us to consider cone metrics.

\medskip
We regard the complete hyperbolic structure of finite volume and
geodesic boundary on $M$ as a hyperbolic cone structure on $\OO$
with all cone angles equal to zero.

Let $m_1,\ldots,m_q$ be the ramification indices of the edges of $\Sigma$
(with respect to a fixed numbering). Throughout the proof of the orbifold
theorem we will consider the following set of hyperbolic cone structures
with fixed ratios for the cone angles. Define:

$$ \JJ(\OO)=\left\{ t\in [0,1]
\,\left\vert
\begin{array}{l}
 \text{there exists a hyperbolic cone structure on }\OO
\\ \text{with cone angles } \left(\frac{2\pi t}{m_1},\ldots,\frac{2\pi
t}{m_q}\right)
\end{array}
 \right. \right\} $$

A hyperbolic cone structure on $\OO$ induces a non-complete
hyperbolic structure on $M$. In particular it has a holonomy
representation $\pi_1(M)\to PSL_2(\mathbb C)$. The variety of
representations $\Hom (\pi_1(M),PSL_2(\mathbb C))$ is an affine
algebraic set, possibly reducible. The group $PSL_2(\mathbb C)$
acts on $\Hom (\pi_1(M),PSL_2(\mathbb C))$ by conjugation, and we
are interested in the quotient. The topological quotient is not
Hausdorff, and one therefore considers the algebraic quotient
 \[
\XX(M)=\Hom (\pi_1(M),PSL_2(\mathbb C))/\! / PSL_2(\mathbb C)
 \]
which is again an affine algebraic set. Note that the irreducible
representations form a Zariski open subset of $\Hom
(\pi_1(M),PSL_2(\mathbb C))$. Namely, $\Hom^{irr}
(\pi_1(M),PSL_2(\mathbb C))$ is the inverse image of a Zariski
open subset $\XX^{irr}(M)\subseteq \XX(M)$, and $\XX^{irr}(M)$ is
the topological quotient of $\Hom^{irr} (\pi_1(M),PSL_2(\mathbb
C))$. Notice that the holonomy representation $\rho_0$ of the
(metrically) complete hyperbolic structure on $M$ is irreducible.

The polynomial functions on $\XX(M)$ one-to-one correspond to the
polynomial functions on $\Hom (\pi_1(M),PSL_2(\mathbb C))$
invariant under the $PSL_2(\mathbb C)$-action. Given
$\gamma\in\pi_1(M)$, we define the trace-like function
$\tau_{\gamma}\!:\XX(M)\to\mathbb C$ as the function induced by
\begin{eqnarray*}
 \Hom (\pi_1(M),PSL_2(\mathbb C)) & \rightarrow & \mathbb C \\
   \rho & \mapsto & \operatorname{trace}(\rho(\gamma)^2).
\end{eqnarray*}

Let $\mu_1,\ldots,\mu_q$ be a family of meridian curves, one for
each component of $\Sigma^{(1)}$.

\begin{thm}[Local parametrization]\label{thm:bianalytic}
The map
\[
\tau_{\mu}=(\tau_{\mu_1},\ldots, \tau_{\mu_q})\!:\XX(M)\to\mathbb
C^q
\]
is locally  bianalytic at $[\rho_0]$.
\end{thm}

This result is the main step in the proof of Thurston's hyperbolic
Dehn filling theorem (see \cite[App. B]{BoP} or \cite{Kap} for the
proof). It implies in particular the following special case of
Thurston's Generalized Hyperbolic Dehn Filling Theorem.

\begin{cor}
\label{thm:GHDF} The set $\JJ(\OO)$ contains a neighborhood of $0$.
\end{cor}
\begin{proof} We have $\tau([\rho_0])=(2,\dots,2)$. Consider the path
\begin{equation}
\label{path}
\begin{array}{rcl}
\gamma\!:[0,\varepsilon) & \rightarrow & \mathbb C^q \\
  t & \mapsto & (2\cos \frac{2\pi t}{m_1},\ldots, 2\cos \frac{2\pi
  t}{m_k})
\end{array}
\end{equation}
where $\varepsilon>0$ is sufficiently small. The composition
$\tau_{\mu}^{-1}\circ\gamma$ gives a path of conjugacy classes of
representations. It can be lifted to a path $t\mapsto\rho_t$, because
there are slices to the action of $PSL_2(\mathbb C)$ on the representation
variety. (Existence of slices follows from Luna's theorem, as
$PSL_2(\mathbb C)\cong SO(3,\mathbb C)$ is an affine reductive group.) The
representations $\rho_t$ are the holonomies of incomplete hyperbolic
structures on $M$. By construction, the holonomies of the meridians are
rotations with angles $2\pi t/m_1,\ldots,2\pi t/m_q$. By a standard
result, the deformation of holonomies is, locally near $t=0$, induced by a
deformation of hyperbolic cone structures on $\OO$ with cone angles $2\pi
t/m_j$. \end{proof}

\begin{lem}\label{lem:curve} There exists a unique irreducible curve
${\mathcal D}\subset\mathbb C^q$ such that $\gamma([0,1])\subset
{\mathcal D}$.
\end{lem}

\begin{proof} For $n\in\mathbb N$, we consider the Chebyshev-like
polynomial
$$ p_n(x)=2\cos(n \arccos(x/2)). $$ It has the following property: $$
{\operatorname{trace}}(A^n)=p_n({\operatorname{trace}}(A)),
\qquad\forall A\in SL_2(\mathbb C), \ \forall n\in{\mathbb N}. $$
An easy computation shows that $p_n'(2)=n$, and therefore
$$ \{z\in\mathbb C^q\mid p_{m_1}(
z_1)=\cdots=p_{m_q}(  z_q)\} $$ is an algebraic curve with
$(2,\dots,2)$ as a smooth point. We take ${\mathcal D}$ to be the
unique irreducible component containing $(2,\dots,2)$. Then
$\gamma([0,\varepsilon))\subset {\mathcal D}$ for small
$\varepsilon>0$. Since $\ga$ is an analytic curve, it remains in
${\mathcal D}$. \end{proof}

We define the algebraic curve $\mathcal C\subset \XX(M)$ to be the
irreducible component of $\tau_{\mu}^{-1}(\mathcal D)$ that
contains $[\rho_0]$. By construction, $[\rho_t]\in\mathcal C$ for
small $t\geq0$.

For technical reasons, we define the following variant of
$\JJ(\OO)$. Here $v_0$ denotes the volume of the complete
hyperbolic structure on $M$.
\begin{equation}\label{eqn:J0}
\JJ_0(\OO)=\left\{ t\in [0,1] \,\left\vert
\begin{array}{l}
 \text{there exists a hyperbolic cone structure on }\OO
\\ \text{with cone angles } \left(\frac{2\pi t}{m_1},\ldots,\frac{2\pi
t}{m_q}\right) \text{, holonomy in } \mathcal C \\
\text{ and volume }\leq v_0
\end{array}
 \right. \right\}
\end{equation}
The condition that the holonomy is in the curve $\mathcal C$ will
be used in Theorem~\ref{thm:open}, because holomorphic maps on
curves are open (openness does not hold for maps on higher
dimensional varieties).

Note that $[0,\varepsilon)\subset \JJ_0(\OO)$ for small
$\varepsilon>0$ because, according to Schl\"afli's formula, the
volume of the continuous family of cone structures with holonomies
$\rho_t$ decreases.

\begin{thm}[Openness]
\label{thm:open}The set $\JJ_0(\OO)$ is open to the right. \end{thm}
\begin{proof} As remarked above, openness of $\JJ_0(\OO)$ at $t=0$ is a
consequence of Thurston's hyperbolic Dehn filling, and we only
prove openness at $t>0$.

Consider the path
\begin{eqnarray*}
\gamma:[t,t+\varepsilon)&\to&{\mathcal D}\subset\mathbb C^q\\
s&\mapsto&(  2\cos(s2\pi/m_1),\ldots,  2 \cos(s2\pi/m_q))
\end{eqnarray*}
defined for some $\varepsilon>0$. By construction, the image of
$\gamma$  is contained in the curve $\mathcal D\subset\mathbb C^q$
of Lemma~\ref{lem:curve}. Since $\tau_{\mu}:\mathcal C\to\mathcal
D$ is non constant, it is open, and therefore $\gamma$ can be
lifted to $\mathcal C$. We can lift it further to a path
\begin{eqnarray*}
\tilde \gamma:[t,t+\varepsilon)&\to&\Hom(\pi_1(M),PSL_2(\mathbb C))\\
s&\mapsto&\rho_s
\end{eqnarray*}
To justify this second lift, notice that the holonomy $\rho_t$ is
irreducible (because the corresponding cone structure has finite
volume) and therefore the $PSL_2(\mathbb C )$-action is locally
free.
 By construction,
$\rho_s(\mu_i)$ is a rotation of angle $\frac{2\pi s}{m_i}$.
Therefore the cone structure on $\OO$ with holonomy $\rho_t$ can
be deformed to a continuous family of cone structures on $\OO$
with holonomies $\rho_s$. By Schl\"afli's formula volume decreases
and thus $[t,t+\varepsilon)\subset\JJ_0(\OO)$  for $\varepsilon>0$
sufficiently small.\end{proof}

A straightforward consequence of Theorem~\ref{thm:open} is:

\begin{cor}
 If $\sup \JJ_0(\OO)    \in \JJ_0(\OO)$ then $1\in \JJ_0(\OO)$ and $\OO$ is
 hyperbolic.
\end{cor}

The next step in the proof is the analysis of degenerating
sequences of cone structures on $\OO$, namely sequences $(t_n)$ in
$\JJ_0(\OO)$ that converge to $t_{\infty} \not\in \JJ_0(\OO)$.
This analysis is carried out in the next section,
 using the results of Sections  \ref{sec:noncollapse},
 \ref{sec:sphericaluniformization} and \ref{sec:Seifert}.

\subsection{Degeneration of hyperbolic cone structures}
\label{sec:deg}

We  continue the discussion of deforming hyperbolic cone
structures on $\OO$ while keeping the ratios of the cone angles
fixed.

 Let $(t_n)$ be a
sequence in $\JJ_0(\OO)$. Let $X(t_n)$ be a hyperbolic cone
structure on $\OO$ corresponding to $t_n\in\JJ_0(\OO)$, with the
properties as in (\ref{eqn:J0}).

\begin{dfn}
We say that a sequence of
cone 3-manifolds $X_n$  \emph{collapses} if, for every $\rho>0$, only
finitely many $X_n$ are $\rho$-thick,
cf.\ Definition \ref{dfn:thick}.
\end{dfn}

\begin{prop}[Degeneration implies collapse]\label{cor:notcollapses} If
  $t_n\to t_{\infty}\not \in \JJ_0(\OO)$,
 then the sequence
 $(X(t_n))$  collapses.
\end{prop}

\begin{proof} Assume that $(X(t_n))$ does not collapse.
Then, up to passing
to a subsequence, the cone manifolds $X(t_n)$ are $\rho$-thick
for some $\rho>0$.
According to Corollary~\ref{thm:compactness},
$(X(t_n))$ subconverges to a hyperbolic cone 3-manifold $X_{\infty}$.
Since in the definition of $\JJ_0(\OO)$ we impose an upper volume bound on the
cone structures,
we have that $\operatorname{vol}(X_{\infty})<\infty$.

Moreover, if $t_{\infty}<1$, the cone angles of $X(t_n)$ are all bounded
away from $\pi$, and if  $t_{\infty}=1$, then they converge to the
orbifold angles of $\OO$. This  allows to apply   Theorem~\ref{thm:stab}
in Section~\ref{sec:noncollapse} below. According to this theorem,
$X_{\infty}$ yields a hyperbolic cone structure on $\OO$ and $t_{\infty}
\in \JJ_0(\OO)$, which contradicts the hypothesis.
\end{proof}


\medskip
We analyse now the situation when the sequence $(X(t_n))$ collapses.
We treat the cases with and without boundary separately.

The case
with boundary is handled by the following geometric fibration
result.

\begin{prop}[$I$-bundle] If
 $\partial X(t_n)\neq\emptyset$ for all $n$, then
$t_{\infty} = 1$ and $\OO$ is a twisted $I$-bundle over the
quotient of a turnover by an orientation reversing involution.
\end{prop}

\begin{proof}
Since the sequence $(X(t_n))$ collapses and $\partial X(t_n)$ is a
collection of turnovers, for $n$ sufficiently large
Corollary~\ref{cor:Ibundle} applies to $X(t_n)$.
It shows that the boundary $\D X(t_n)$ consists of a single turnover
and that $X(t_n)$ is a singular interval bundle
over the cut locus $Cut(\D X(t_n))$ with respect to $\D X(t_n)$.
The cut locus $Cut(\D X(t_n))$ is naturally homeomorphic to the quotient of
the turnover $\D X(t_n)$ by an orientation reversing isometric involution.

Moreover,
since the $X(t_n)$ are negatively curved,
$\D Cut(\D X(t_n))$ is a non-empty collection of singular edges,
cf.\ the discussion after Proposition \ref{Ibundlethm}.
By Corollary~\ref{cor:Ibundle},
their cone angles converge to $\pi$ as $n\to\infty$.
Thus $t_{n}\to 1$.
It follows that the boundary turnovers of the $X(t_n)$ correspond to a
hyperbolic boundary turnover of $\OO$.
Hence the $X(t_n)$ provide not only cone structures but cone {\em metrics} on
$\OO$,
and $\OO$ is an $I$-bundle over a quotient of a turnover,
as claimed.
\end{proof}

\begin{rem}
A turnover has always an isometric involution which fixes all three
cone points and reverses orientation.
The quotient is a triangular 2-orbifold.

The turnover is the double along the boundary of a geodesic triangle, and
if the triangle has a reflection symmetry, then the turnover has a
corresponding orientation reversing involution which fixes one cone point
and exchanges the other two. In this case the quotient is a disc with a
corner and one cone point. Figure~\ref{fig:7.1} illustrates the orientable
$I$-bundles over such 2-orbifolds.
\end{rem}

\begin{figure}
\begin{center}
 \includegraphics[scale=.75]{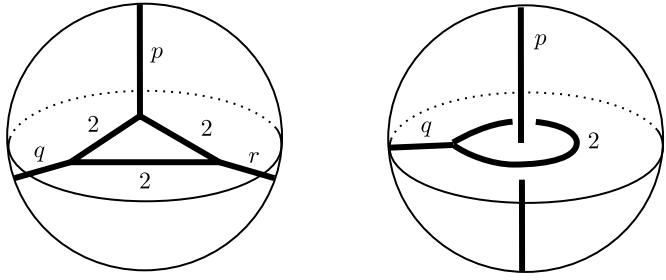}
 \end{center}
   \caption{The two $I$-bundles over the
    quotient of a turnover}\label{fig:7.1}
\end{figure}

Now we discuss the case without boundary,
i.e.\ we assume that $\D X(t_n)=\emptyset$ for all $n$.
Note that all cusps of the $X(t_n)$ are singular
because $\OO$ is small.
Since, without loss of generality,
the sequence $(t_n)$ is strictly increasing,
only finitely many $X(t_n)$ can have cusps
and the cusps correspond to singular vertices of $\OO$.
After passing to a subsequence,
we may assume that the $X(t_n)$ have no cusps at all,
i.e.\ they are closed cone manifolds,
$\OO$ is a closed orbifold
and the $X(t_n)$ provide cone {\em metrics} on $\OO$.
In particular,
the pairs
$(\vert \OO\vert,\Sigma)$ and $(X(t_n),\Sigma_{X(t_n)})$
are homeomorphic for all $n$.

\begin{dfn}
\label{dfn:thicklinks}
We say that the cone 3-manifold $X$ has
{\em $\om$-thick links}, $\om>0$,
if in every point $x\in X$ the link $\La_xX$ is $\omega$-thick,
cf.\ Definition \ref{dfn:thick}.
\end{dfn}

\begin{thm}[Fibration]\label{thm:fibration}
For $\om,D_0>0$ there exists $\de=\de(\om,D_0)>0$ such that:

Suppose that $\OO$ is closed and admits a cone metric $(X,i)$
of constant curvature $k\in[-1,0)$
with  cone angles less than or equal to the orbifold angles.
If the cone 3-manifold
$X$ has $\om$-thick links, $\Diam(X) \geq D_0$ and if $X$ is $\de$-thin,
then $\OO$ is Seifert fibred.
\end{thm}

We postpone the proof of this theorem to Section~\ref{sec:Seifert}.
Assuming it, we continue our argument.

\begin{prop}\label{cor:collapses}
If  $\partial X(t_n) = \emptyset$ for all $n$, then one of the
following holds:
\begin{itemize}
    \item[--] either there exists a Euclidean cone metric on $\OO$ with cone angles
      strictly less than the orbifold angles,
     \item[--] or  $\OO$ has a Euclidean or a Seifert
      fibred structure.
\end{itemize}
\end{prop}

\begin{proof}
{\em Case 1: Collapse with lower diameter bound.}
Assume that $diam(X(t_n))\geq D_0>0$ for all $n$.
Since $t_n<t_{\infty}$ and the cone angles for the cone metrics $X_n$ are
proportional to $t_n$,
we have a uniform lower bound $2\pi+\eps$, $\eps>0$,
for the cone angle sum of the singular edges adjacent to any vertex
with respect to all cone metrics $X(t_n)$.
By Gau\3-Bonnet, we obtain a lower bound for the area of links at singular
vertices of the $X(t_n)$,
and this converts into a lower bound for the thickness of these links.
Since the cone angles of the $X(t_n)$ are bounded away from zero,
the cone manifolds $X(t_n)$ have uniformly thick links.
Hence Theorem~\ref{thm:fibration} applies and we obtain that $\OO$ is Seifert
fibred.

{\em Case 2: Collapse to a point.} Assume that $\Diam(X(t_n))\to 0$. Then
we consider the sequence
\[
\overline X(t_n)= \frac1{\Diam(X(t_n))}X(t_n)
\]
of rescaled cone 3-manifolds with constant curvature $k_n =
-\Diam(X(t_n))^2 \in [-1,0)$ and diameter equal to $1$.
If the sequence $(\ol X(t_n))$ collapses, then for $n$ sufficiently large
Theorem~\ref{thm:fibration} applies as above to  show that $\OO$
is Seifert fibred.

If the rescaled sequence $\overline X(t_n)$ does not collapse,
then by the compactness result \ref{thm:compactness},
a subsequence converges
geometrically to a closed Euclidean cone 3-manifold $X_{\infty}$
with diameter 1,
and $X_{\infty}$ yields a Euclidean cone metric on $\OO$.
The cone angles of $X_{\infty}$ are either strictly less than the orbifold
angles of $\OO$,
or  $X_{\infty}$ corresponds to a Euclidean {\em orbifold} structure on
$\OO$ with the same branching indices as $\OO$, and so $\OO$ is
Euclidean.
\end{proof}

The last step of the proof of the main theorem is given by:

\begin{thm}[Spherical uniformization]\label{thm:su}
Let $\OO$ be a closed orientable
small 3-orbifold.
If there exists a Euclidean cone metric on $\OO$ with cone angles
strictly less than the orbifold angles, then $\OO$ is
spherical.
\end{thm}

\begin{figure}
\begin{center}
 \includegraphics[scale=1]{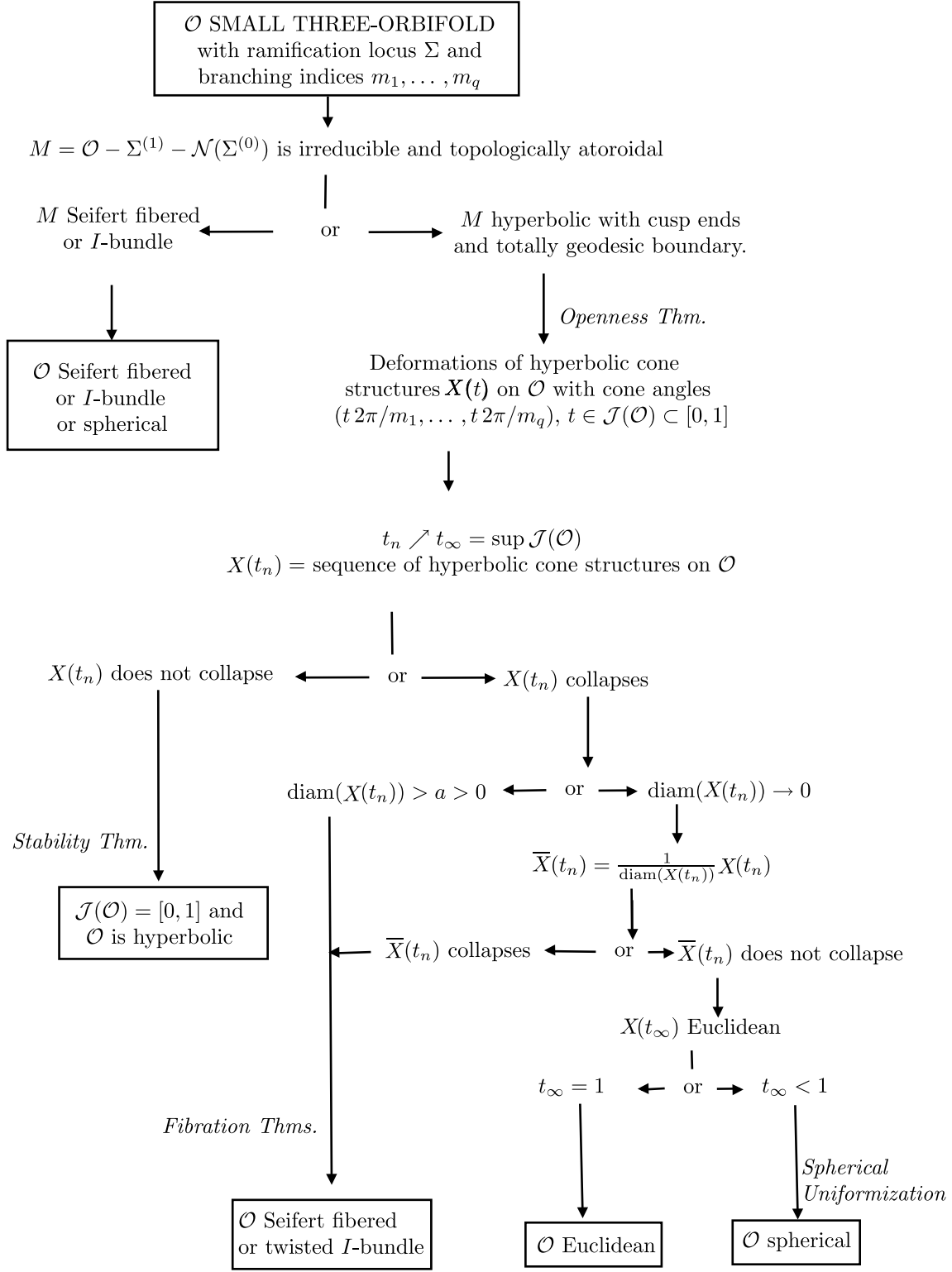}
 \end{center}
    \caption{Flowchart of the proof}\label{fig:fc}
\end{figure}


This theorem is proved in Sections \ref{sec:sphericaluniformization} and
\ref{sec:sphericalcone}. The flowchart in Figure~\ref{fig:fc} represents
the logic of the proof of the main theorem given in this section.

\section{Topological stability of geometric limits}
\label{sec:noncollapse}

In this section we discuss the change of the topological type of
cone manifolds under geometric limits. More specifically,
we consider a sequence of compact 
hyperbolic cone $3$-manifolds $X_n$ with cone angles $\leq\pi$. We
suppose furthermore that the sequence $(X_n)$ is {\em
non-collapsing}, i.e.\ that for a uniform radius $\rho>0$ the
$X_n$ contain embedded smooth standard balls $B_{\rho}(x_n)$.

Due to the compactness theorem (Corollary~\ref{thm:compactness}),
after passing to a subsequence, the pointed cone manifolds
$(X_n,x_n)$ geometrically converge
\begin{equation}
\label{geomconv} (X_n,x_n) \lra (X_{\infty},x_{\infty})  .
\end{equation}
Their limit $X_{\infty}$ is again a
hyperbolic cone $3$-manifold with cone angles $\leq\pi$. We know furthermore, that the singular sets converge,
$\Si_n\ra\Si_{\infty}$, and the cone angles converge as well.

If the limit $X_{\infty}$ is compact, then topological stability
holds, that is, $(X_n,\Si_n)$ is homeomorphic to
$(X_{\infty},\Si_{\infty})$ for sufficiently large $n$. In the
following, we will study the situation when $X_{\infty}$ is {\em
noncompact}. In order to obtain topological stability, we need to
impose further assumptions. The main result of this section is:

\begin{thm}[Stability in the noncompact limit case]
\label{thm:stab} Let $(X_n)_{n\in\mathbb N}$ be a sequence of
compact oriented hyperbolic cone 3-man\-ifolds with
totally geodesic boundary and with cone angles $\in[\om,\pi]$,
$\om>0$, which geometrically converges, as in (\ref{geomconv}), to
a {\em noncompact} cone 3-manifold . Assume that:
\begin{enumerate}
\item
The cone manifolds $X_n$ yield hyperbolic cone structures on the
same compact {\em small} orbifold $\OO$.
\item
There is a uniform upper volume bound $vol(X_n) \leq v$.
\item
Either the cone angles of the $X_n$ are $\leq\alpha<\pi$, or they
converge to the orbifold angles of $\OO$.
\end{enumerate}
Then, $X_{\infty}$ yields a hyperbolic cone structure on $\OO$ as
well.
\end{thm}
\begin{rem}
One can show that the second assumption is implied by the first
one using a straightening argument for triangulations under the
developing map, cf.\ \cite{Koj} for the cyclic case. However, this
is irrelevant in our later applications, because there we will
have a uniform volume bound by construction.
\end{rem}
The proof occupies the entire section and will be divided into
several steps.

\subsection{Case of cone angles $\leq\al<\pi$}
\label{subsect:coneangleslesspi}

We consider now the case that the $X_n$ have cone angles
$\leq\alpha<\pi$. Here the geometric results of Subsection
\ref{sec:thin} will come to bear.

\subsubsection*{Geometric finiteness of the limit}

The uniform upper volume bound for $X_n$ implies that the limit
cone manifold has finite volume:
\[ vol(X_{\infty}) \leq v \]
Since $\OO$ is small, the boundary components of $X_n$ are
turnovers. The lower bound on cone angles yields an upper bound on
their diameters, and it follows that the boundary components of
$X_{\infty}$ are also turnovers. Since the cone angles of
$X_{\infty}$ are bounded away from $\pi$, finite volume implies
geometric finiteness, i.e.\ $X_{\infty}$ has finitely many ends
and all ends are cusps. Here, we apply the finiteness corollary
(Cor.~\ref{finiteness}).

In each cusp $\CC_i\subset X_{\infty}$ we fix a horospherical
cross section $H_i$ far out in the thin part, along which we
truncate $\CC_i$. Denote by $N_{\infty}$ the resulting compact
core of $X_{\infty}$. The $H_i$ are Euclidean tori or turnovers.
By geometric convergence, for sufficiently large $n$, there are
$(1+\eps_n)$-bilipschitz embeddings
\[ f_n:N_{\infty}\hookrightarrow X_n ,\]
$\eps_n\ra0$, and we can arrange that $H_{i,n}:=f_n(H_i)$ is a
leaf of the canonical locally-homogeneous foliation of the thin
part of $X_n$ (namely, either a horosphere in a cusp or a torus
equidistant to a short geodesic, possibly singular, cf. Theorem
\ref{descthin}).

\subsubsection*{Hyperbolizing the smooth part}

We denote by $Y$ the manifold obtained from
$X_{\infty}-\NN(\Si^{(0)}_{\infty})-\Si^{(1)}_{\infty}$ by
truncating the singular cusps along horospherical turnovers. Here,
we remove {\em open} singular cusps. The boundary of $Y$ is a
union of thrice-punctured spheres because, due to our assumption
on cone angles, the cross sections of singular cusps are
turnovers.


\begin{prop}\label{Yhyperbolic}
The manifold $Y$ admits a metrically complete hyperbolic structure with
totally geodesic boundary and finite volume.
\end{prop}
\begin{proof} We first deal with the (easy) case that $Y$ has empty
boundary. This happens if and only if $X_{\infty}$ has no singular
vertices, no singular cusps and empty boundary. There is nothing to show
if $X_{\infty}$ has no singular locus. If there is a singular locus, one
can perturb the hyperbolic metric on the smooth part to a complete
Riemannian metric with upper negative curvature bound and finite volume.
Since $Y$ is Haken, it follows from Thurston's uniformization theorem that
$Y$ admits a complete hyperbolic metric of finite volume.

We assume from now on that $\D Y\neq\emptyset$. Let $\overline Y$
be a compact core of $Y$.  The boundary $\D \overline Y $
 is the union
of compact pairs of pants   (which are a compact core of $\partial
Y$) together with a collection $P \subset\partial\overline Y$ of
tori and annuli, corresponding to the boundary of a neighborhood
of edges and circles in $\Sigma_{\infty}$ and to cross-sections of
smooth cusps. We shall prove that $(\overline Y, P)$ is an
atoroidal pared manifold (see Lemma \ref{lem:pared} for the
definition). Then Proposition~\ref{Yhyperbolic} shall follow from
Thurston's hyperbolization theorem for atoroidal Haken pared
3-manifolds.

The argument will involve several steps. Part of the information
on the topology of $Y$ we obtain from putting a weak geometric
structure, part of it we deduce from the fact that $X_{\infty}$
arises as a limit of cone manifold structures on small orbifolds.

{\em Step 1: $\overline  Y$ is homotopically atoroidal.}
 We shall show in Proposition~\ref{cat0oncompactcore} (in Subsection~\ref{subsec:cat})
that $\overline  Y$ admits a non-positively curved metric,
possibly singular, which is negatively curved away from the
boundary tori. Consider a $\pi_1$-injective map
$h:T\hookrightarrow \overline  Y$ of a 2-torus. The group
$\pi_1(T)\cong\Z^2$ acts discretely on $\widetilde{\overline Y}$
and hence it preserves a 2-flat  (cf. the preliminaries of
\cite{LeebScott}). This 2-flat covers a boundary torus and it
follows that $h$ can be homotoped into a boundary torus.

{\em Step 2: $\overline  Y$ is irreducible.} Let $S\subset \overline  Y$
be an embedded 2-sphere. We may assume that we have the inclusion of
compact cores $\overline Y\subset N_{\infty}$.
The sphere $f_n(S)$ bounds a smooth ball $B_n$ in $X_n$ (for $n$
large). We are done, if $B_n\subset f_n(N_{\infty})$ for some $n$.
Otherwise, each $B_n$ contains a smooth cross section $H_{i,n}$,
which is a leaf of the canonical foliation of the thin part (see
above). After passing to a subsequence, we may assume that $B_n$
contains $H_{i_0,n}$ for a fixed $i_0$. This is absurd because
$H_{i_0}$ has non-trivial holonomy, and so has $H_{i_0,n}$ for
large $n$.

{\em Step 3:  $P$ is incompressible.} Each component  $P_i\subset
P$ corresponds to either cross sections of smooth cusps or tubular
neighborhoods of singular edges and circles of the cone manifold
$X_{\infty}$. Hence the holonomy of primitive elements of
$\pi_1(P_i)$ is non trivial.

{\em Step 4:   The pair $(\overline Y, P)$ is anannular.}
 Let $(A,\partial A)\subset (\overline Y, P)$ be an
essential annulus; we distinguish three cases according to whether
$\partial A$ is contained in a) torus components of $P$, b)
annulus components of $P$, or c) a torus and an annulus of $P$. In
the first  case a classical argument using the atoroidality of
$\overline Y$ implies that $\overline Y$ would be Seifert fibred
\cite[Lemma 7]{BSOne}, contradicting Step 1. In case b), the
annulus $A$ gives rise to an embedded 2-sphere with two cone
points in $X_{\infty}$ and hence (via geometric convergence
(\ref{geomconv})) in $\OO$. Such a 2-sphere bounds an embedded
ball with a singular axis because $\OO$ is small. It follows that
$A$ is parallel to a component of $P$. Case c) reduces to case b),
as in the proof of Lemma \ref{lem:pared}.
 \end{proof}

\subsubsection*{Controlling the geometry of the approaching cone manifolds
globally.} The Gromov-Hausdorff convergence (\ref{geomconv}) gives
us uniform control on the geometry of all $X_n$ a priori only on
the portions $f_n(N_{\infty})$. Taking into account the structure
of the thin part of cone manifolds and using the smallness of
$\OO$, we will be able to also describe the geometry of the
complements $X_n-f_n(N_{\infty})$ and see that it is very
restricted.

\begin{lem}
\label{topapproxends} Each component of $X_n-f_n(N_{\infty})$ is
contained in the thin part of $X_n$, and it is either
\begin{itemize}
\item
a singular ball,
\item
a singular neck containing a turnover $\subseteq\partial X_n$,
\item
or a (singular or smooth) solid torus.
\end{itemize}
\end{lem}

\begin{rem}
A singular cusp cannot occur in the conclusion of
Lemma~\ref{topapproxends} because we assume the $X_n$ to be
compact.
\end{rem}

\begin{proof} If $H_i$ is a turnover, then $H_{i,n}$ is an umbilic
turnover. We go through the three possible cases: If $H_{i,n}$  is
spherical, then it bounds a singular ball in $X_n$. The turnover
$H_{i,n}$ cannot be horospherical because $X_n$ is compact. If
$H_{i,n}$ is hyperspherical, then there is an umbilic tube bounded
by $H_{i,n}$ and a totally geodesic  turnover in $\partial X_n$.
Here we use that $\OO$ is small.

If $H_i$ is a torus, then $H_{i,n} \subset X_n$ is an almost
horospherical torus. It cannot be precisely horospherical, again
because $X_n$ is compact. Hence $H_{i,n}$ bounds  a (smooth or
singular)  solid torus $V_{i,n} \subset X_n$. \end{proof}

\medskip
If $H_i$ is a torus, denote by $\lambda_{i,n}\subset H_{i,n}$ a
geodesic which is a meridian of (i.e.\ compresses in) the solid
torus $V_{i,n}$ bounded by $H_{i,n}$. Furthermore, denote by
$\tilde\lambda_{i,n}\subset H_i$ a geodesic such that
$f_n(\tilde\lambda_{i,n})$ is homotopic to $\lambda_{i,n}$. The
lengths of $\tilde\lambda_{i,n}$ and $\lambda_{i,n}$ are
comparable in terms of the bilipschitz constant of $f_n$.

\begin{lem}
\label{meridiangrows} For all $i$, we have
$\lim\limits_{n\ra\infty}length(\lambda_{i,n})=\infty$.
\end{lem}
\begin{proof} The radii of the solid tori $V_{i,n}$ tend to $\infty$ as
$n\ra\infty$. Using the lower bound on cone angles if $V_{i,n}$ is
singular, the assertion follows. \end{proof}

\subsubsection*{Comparing the topology with the limit}
Using these geometric observations, we now describe the change of
topology during the transition
$(X_{\infty},\Si_{\infty})\leadsto(X_n,\Si_n)=(|\OO|,\Si)$. We
remove from the cone manifolds $X_n$ and the truncated cone
manifold $N_{\infty}$ disjoint small open standard balls around
the singular vertices, and denote by $\Si'_n$ and $\Si'_{\infty}$
the singular loci which remain in the resulting manifolds
$X_n-\NN_{\eps}(\Si^{(0)}_n)$ and
$N_{\infty}-\NN_{\eps}(\Si^{(0)}_{\infty})$. Lemma
\ref{topapproxends} implies that the transition
\begin{equation}
\label{transtoptype}
\bigl(N_{\infty}-\NN_{\eps}(\Si^{(0)}_{\infty}),\Si'_{\infty}\bigr)\leadsto
        \bigl(X_n-\NN_{\eps}(\Si^{(0)}_n),\Si'_n\bigr)
\end{equation}
is done by gluing (smooth or singular) solid tori to the boundary
tori of $N_{\infty}$ (i.e.\ to the smooth $H_i$).

Correspondingly, the transition
\begin{equation}
\label{transsmoothpart}
\underbrace{N_{\infty}-\NN_{\eps}(\Si^{(0)}_{\infty})-\Si_{\infty}}_{
      =:Y_{\infty}}
\leadsto X_n-\NN_{\eps}(\Si^{(0)}_n)-\Si_n
\end{equation}
between the smooth parts is done by gluing (smooth) solid tori to
some of the boundary tori and removing the remaining boundary
tori.

Every smooth cusp of $X_{\infty}$ corresponds to a smooth boundary
torus of $N_{\infty}$. According to Lemma~\ref{meridiangrows}, we
may pass to a subsequence such that the Dehn fillings at every
such torus (by smooth or singular solid tori) are pairwise
different.

We distinguish three cases:

{\em Case 1: Infinitely many Dehn fillings by smooth solid tori.}
We look at the smooth parts. Then by Proposition~\ref{Yhyperbolic}
we are in the situation that infinitely many Dehn filings at the
same finite volume hyperbolic manifold produce homeomorphic
manifolds. This contradicts Thurston's hyperbolic Dehn filling
theorem, which implies that those manifolds have different
hyperbolic volumes.

{\em Case 2: $X_{\infty}$ has smooth cusps and all Dehn fillings
use singular solid tori.} Let $Y_{\infty}:=
N_{\infty}-\NN_{\eps}(\Si^{(0)}_{\infty})-\Si_{\infty}$ and let
$Y_n:=f_n(Y_{\infty})$ be its image under the bilipschitz
embedding $f_n$. By our description of (\ref{transsmoothpart}),
$f_n|_{Y_{\infty}}$ may be isotoped to an embedding
\[ f'_n:Y_{\infty}\hookrightarrow X_n-\NN_{\eps}(\Si^{(0)}_n)-\Si_n \]
whose image is the complement of a union of open tubular
neighborhoods of singular circles. We use our assumption that the
cone manifolds $X_n$ have the same topological type as the
orbifold $\OO$, i.e. we have an embedding
$(X_n,\Sigma_{X_n})\hookrightarrow (\vert\mathcal
O\vert,\Sigma_{\mathcal O})$, that we compose with $f'_n$ to
obtain a new embedding
\begin{equation}
\label{embedintoO} \iota_n: Y_{\infty}\hookrightarrow M :=
|\OO|-\NN(\Sigma^{(0)})-\Sigma^{(1)}
\end{equation}
onto the complement of a disjoint union of singular solid tori.

The $\iota_n$ are homotopy equivalences. Note  that $Y_{\infty}$
is homotopy equivalent to a complete hyperbolic manifold with
finite volume and totally geodesic boundary; this follows, without
using Proposition~\ref{Yhyperbolic},  directly from the original
assumption that $M$ is hyperbolic.
 We consider the homotopy equivalence
 \[
 \iota_n^{-1}\circ\iota_1: Y_{\infty}\lra Y_{\infty}
 \]
where $\iota_n^{-1}$ denotes a homotopy inverse for $\iota_n$.
After passing to a subsequence, and possibly replacing the first
embedding $\iota_1$, we may assume that the
$\iota_n^{-1}\circ\iota_1$ preserve the toral boundary components.
Lemma \ref{meridiangrows} implies that the induced self homotopy
equivalences of each boundary torus are pairwise distinct. On the
other hand, by Mostow Rigidity, there are up to homotopy only
finitely many self homotopy equivalences of $Y_{\infty}$, and we
obtain a contradiction.

{\em Case 3: $X_{\infty}$ has no smooth cusp}. In this situation,
(\ref{transtoptype}) says that $X_{\infty}$ and $X_n$ have the
same topological type. This finishes the proof of
Theorem~\ref{thm:stab} when cone angles are bounded above away
from $\pi$. \end{proof}

\subsection{Case when cone angles approach the orbifold angles}

In this case   $X_{\infty}$ is a hyperbolic orbifold that has
 a thin-thick decomposition, by the
hypothesis about the cone angles and the bound on the volume.
 Let $N_{\infty}$ be a compact
core of $X_{\infty}$ obtained by truncating its cusps along
horospherical cross sections $H_i$. The $H_i$ are now smooth tori,
pillows or turnovers. By geometric convergence, for $n$ large
enough, we have
 $(1 + \varepsilon_n)$-bilipschitz embeddings
\[
f_n\! : (N_{\infty}, \Sigma_{\infty}\cap N_{\infty}) \rightarrow
(X_n, \Sigma_n),
\]
 with $\varepsilon_n \rightarrow 0$. Since $X_n$ is a cone structure on $\OO$,
 we view the image $f_n(N_{\infty})\subset X_n$ as a
suborbifold of $\OO$, which we denote by $N_n \subset \OO$. As a
3-orbifold, $N_n$ is homeomorphic to $N_{\infty}$.

\begin{lem} For $n$ sufficiently large, each component of
$\OO-\Int(N_n)$ is irreducible.
\end{lem}

\begin{proof} We assume that $\OO-\Int(N_n)$ contains a spherical
$2$-suborbifold $F^2$ which is essential. Since $\OO$ is
irreducible, $F^2$ bounds a discal $3$-orbifold $\Delta^3$, which
contains $N_n$, for $n$ sufficiently large. Let $\rho_n$ denote
the holonomy representation of $X_n$. We have that $\rho_n\circ
f_{n*}$ fixes a point of $\mathbb H^3$, because $f_n(N_{\infty})$
is contained in a discal 3-orbifold. This is impossible, because
$\rho_n\circ f_{n*}$ converges to the holonomy representation of
$N_{\infty}$, which cannot fix
 a point of $\overline{\mathbb H^3}$. \end{proof}

From smallness of $\OO$ and the previous lemma we obtain:

\begin{cor}\label{cor:solidtori} Every component of $\OO-\Int(N_n)$ is
either  a
  finite quotient of a solid torus (i.e. a solid torus or a
solid pillow, possibly singular) or a singular neck containing a
Euclidean turnover in $\D\OO$.
\end{cor}

According to this corollary, if  all the horospherical sections
 $H_i$ are turnovers then $N_{\infty}\cong \OO$
and $X_{\infty}$ is a cone structure on $\OO$.

Now we assume that some of the horospherical sections
$H_1,\ldots,H_p$ are tori or pillows and we look for a
contradiction.
 For $i=1,\ldots,p$, let $\lambda_{i,n}$ be an
essential curve on $H_i$ so that $f_n( \lambda_{i,n})$  bounds a
finite quotient of a solid torus component of $\OO-\Int(N_n)$.
First we prove the analogue of Lemma~\ref{meridiangrows}.

\begin{lem}\label{lem:curvesinfty} For each $i$,
$\lim\limits_{n\to\infty}\lambda_{i,n}=\infty$.
\end{lem}

\begin{proof} Suppose that  the lemma is false. Then, by passing to a
subsequence and changing the indices of the $H_i$, we can assume
that the curves $ \lambda_{1,n}$ represent a fixed class $
\lambda_1\in\pi_1(H_1)$ independent of $n$. Let $\rho_n$ and
$\rho_{\infty}$ denote the holonomy representation of $X_n$ and
$X_{\infty}$, respectively.
 Since the curves
$f_n(\lambda_{1})$ are compressible in $\OO$, their holonomies
$\rho_n(f_n(\lambda_{1}))$ are either trivial or elliptic with an
angle that does not converge to zero.
 The holonomy $\rho_{\infty}(\lambda_1)$ is
non-trivial and parabolic. Thus we obtain a contradiction because
$\rho_{\infty}(\lambda_1)$ is the limit of $\rho_n(f_n(
\lambda_{1}))$. \end{proof}

According to this lemma, Dehn fillings on the hyperbolic orbifold
$N_{\infty}$ along infinitely many different meridian curves
produce always the same orbifold $\OO$. This contradicts the
hyperbolic Dehn filling theorem for orbifolds
 \cite{DuM} (cf. \cite[App. B]{BoP}) because the results of surgery along
those curves can be distinguished either by an estimate of the
volume (obtained from Schl\"afli's formula) or the length of the
shortest geodesics. \qed

\subsection{Putting a CAT(-1)-structure on the smooth
part of a cone manifold}
\label{subsec:cat}

Let $X$ be a hyperbolic cone $3$-manifold with cone angles
$\leq\pi$ and totally geodesic boundary. To simplify our discussion, we
will assume that $X$ is {\em geometrically finite} in the sense that it
has finitely many ends and all ends are cusps.
The next proposition is used in
Subsection~\ref{subsect:coneangleslesspi} to show that the smooth part of
$X$ is homotopically atoroidal.

\begin{prop}
\label{cat0oncompactcore} The compact core of $X^{smooth}$ admits a
metric of nonpositive curvature with piecewise totally geodesic boundary.
Moreover, the metric is strictly negatively curved away from the boundary
tori corresponding to the smooth cusps and singular closed geodesics of
$X$.
\end{prop}

\begin{proof} {\em Step 1: Removing neighborhoods of the singular
vertices and truncating singular cusps.} Consider first a vertex
$v\in\Si_{X}^{(0)}$. We choose a small positive number $\rho_v<\half
r_{inj}(v)$ and denote by $w_1,w_2,w_3\in\Si^{(1)}$ the three singular
points at distance $\rho_v$ from $v$. We take the convex hull of
$\{w_1,w_2,w_3\}$ inside the closed ball $\ol{B_{\rho_v}(v)}$, and denote
by $U_v$ the interior of the convex hull. Its closure $\bar U_v$ can be
obtained by doubling a hyperbolic simplex along three of its faces; $\D
U_v$ is the union of two geodesic triangles glued along their boundaries.

Consider now a singular cusp $\CC\subset X$ with horospherical cross
section $H$. Since $X$ is orientable, $H$ is a Euclidean cone sphere with
three or four cone points. As before we form the convex hull of
$H\cap\Si_{X}$ inside $\CC$ and denote its interior by $U_{\CC}$. Then $\D
U_{\CC}$ is piecewise geodesic with vertices in $H\cap\Si_{X}$.

By taking out the neighborhoods $U_v$ around all singular vertices $v$ and
truncating all singular cusps $\CC$, we obtain a hyperbolic cone manifold
\[ X':=X-\bigcup_{v}U_v-\bigcup_{\CC}U_{\CC} \]
with piecewise totally geodesic concave boundary.

{\em Step 2: Removing neighborhoods around the singular edges.} The cone
manifold $X'$ has no singular vertices any more. Its singular locus
$\Si_{X'}:=\Si_{X}\cap X'$ consists of closed singular geodesics and of
singular segments with endpoints in the boundary. We now treat the latter
ones.

Consider a singular edge $\si=\ol{ww'}$ in $X'\cap\Si_{X}$ with
endpoints in $w,w'\in\D X'$. We will work inside a tubular
$\rho_{\si}$-neighborhood $T$ of $\si$ in $X'$ with small radius
$\rho_{\si}$. Choose an interior point $m$ of $\si$ and a little
totally geodesic disc $\De$ orthogonal to $\si$ and centered at
$m$. Consider one of the endpoints of $\si$, say $w$. The boundary
of $X'$ is concave at $w$. (This includes the possibility of being
totally geodesic, which we regard as weak concavity). The link
$\L_w$ of $w \in X'$ is a spherical polygon (in most cases a
bigon) with one cone point and concave boundary. We denote its
vertices by $\xi_1, \xi_2,\ldots, \xi_k$. (In the totally geodesic
case, $\D\L_w$ is a circle and the $\xi_i$ are not well-defined;
we then choose two opposite points $\xi_1$ and $\xi_2$ on the
circle.) For $0<\de_i<<\rho_{\si}$, let $y_i\in T\cap\D X'$ be the
points with $d(w,y_i)=\de_i$ and $\stackrel{\lra}{wy_i}=\xi_i$. If
the $\de_i$ are sufficiently small, then there exist boundary
points $z_i\in\D X'\cap T$ near $w'$ and geodesic segments
$c_i=\ol{y_iz_i}$ inside $T$ intersecting $\De$ orthogonally.
Exchanging the roles of $w$ and $w'$, we construct analogously
geodesic segments $c_1',\ldots,c_l'$, with $l\geq 2$. We
coordinate the $\de_i$ and $\de'_i$ so that the $k+l$ segments
$c_1,\ldots,c_k,c_1',\ldots,c_l'$ intersect $\De$ at the same
distance from $m$. Now we form the convex hull of
$c_1,\ldots,c_k,c_1',\ldots,c_l'$ inside $T$ and denote its
interior by $U_{\si}$. Since the $k+l$ segments intersect $\De$
orthogonally, the boundary of $U_{\si}$ is the union of $k+l$
totally geodesic quadrilaterals.

We perform this construction for all singular edges $\si$ so that the
closed neighborhoods $\bar U_{\si}$ are disjoint. Removing the
neighborhoods $U_{\si}$ for all singular edges with endpoints in $\D X'$
yields a compact cone manifold
\[ X^{carved}:=X'-\bigcup_{\si} U_{\si} \]
with piecewise totally geodesic boundary.

\begin{lem}
\label{linkscat1} The cone manifold $X^{carved}$ is locally CAT(-1) near
the boundary.
\end{lem}
\begin{proof}
Since $X^{carved}$ is everywhere locally conical, we have to check
that the links at all boundary points are CAT(1). This is
nontrivial only at the boundary vertices. In our notation above,
these are the points $y_i$ and $z_i$. The link of both $y_i$ and
$z_i$ in $X'$ is a concave bigon when they lie on an edge of
$\partial X'$, and a half sphere when they lie on a totally
geodesic piece of $\partial X'$. In the first case, the link in $
X^{carved}$ is the intersection of two concave bigons, such that
each bigon contains precisely one vertex of the other one. This
link is a quadrilateral, with two opposite vertices of angle
$>\pi$ and two other vertices of angle $<\pi$. Since each bigon
contains precisely one vertex of the other one, the secant between
the concave vertices divides the quadrilateral into two convex
triangles. Hence the link contains no nonconstant closed geodesics
of length $<2\pi$ -- actually no nonconstant closed geodesics at
all -- and therefore it is CAT(1) due to a criterion by Gromov,
cf. \cite[Sec.\ 4.2 A and B]{hypgroups}. In the other case, when
the vertex lies in a totally geodesic piece of $\partial X'$, its
link in $X^{carved}$ is the intersection of a half sphere with a
concave bigon (with precisely one vertex in the half sphere), and
the argument is similar.
\end{proof}

{\em Step 3: Modification near closed singular geodesics and at smooth
cusps.} We truncate the smooth cusps along horospherical torus cross
sections, respectively, remove open tubular neighborhoods of small radii
around the singular closed geodesics. Then, by a standard double warped
product construction, we perturb the metric locally near the new boundary
components to a nonpositively curved Riemannian metric with totally
geodesic flat boundary. We obtain a compact nonpositively curved
Riemannian manifold $X^{core}$ with piecewise totally geodesic boundary.
The boundary components are either totally geodesic (flat or hyperbolic),
or piecewise totally geodesic hyperbolic. Lemma \ref{linkscat1} implies
that $X^{core}$ is CAT(0). Topologically, $X^{core}$ is a compact core for
$X^{smooth}$. This concludes the proof of Proposition
\ref{cat0oncompactcore}. \end{proof}

\section{Spherical uniformization}\label{sec:sphericaluniformization}

This section and the following one are devoted to the proof of the
spherical uniformization theorem:

\medskip

\begin{sphtheorem} \textbf{\em {(Thm.~\ref{thm:su}).}}
Let $\OO$ be a closed orientable connected small 3-orbifold. If
there exists a Euclidean cone structure on $\OO$ with cone angles
strictly less than the orbifold angles of $\OO$, then $\OO$ is
spherical.
\end{sphtheorem}

\medskip

Note first that a Euclidean cone structure $X$ on $\OO$ could have
boundary. If $\D X$ is non-empty, then it consists of totally
geodesic turnovers. Due to our assumptions, $X$ has cone angles
$<\pi$, and the classification of noncompact Euclidean cone
manifolds (cf.\ Section \ref{sec:euc} above) implies that $X$ is
the product of a Euclidean turnover with an interval. In this case
$\OO$ is a suspension of a spherical turnover and therefore
obviously spherical.

From now on we assume that the Euclidean cone structure $X$ on
$\OO$ is closed. Then $X$ and $\OO$ have the same homeomorphism
type (as pairs of topological space and singular locus), and we
can consider $X$ as a Euclidean cone {\em metric} on $\OO$.

For the proof of Theorem~\ref{thm:su} we distinguish three cases which are
listed in Definition \ref{dfn:orbitype}.

\begin{dfn}
A singular vertex of a (locally) orientable 3-orbifold is called
\begin{itemize}
\item {\em dihedral} if its local isotropy group is a dihedral group,
and
\item {\em platonic} otherwise. (Its local isotropy group then is
the group of orientation preserving isometries of a platonic
solid.)
\end{itemize}
\end{dfn}

\begin{dfn}
\label{dfn:orbitype} A (locally) orientable 3-orbifold is called
\begin{itemize}
\item of \emph{cyclic type}
if its singular locus is not empty and has no vertex,
\item of \emph{dihedral type} if it has singular vertices
and all vertices are dihedral,
\item of \emph{platonic type} if it has a platonic singular vertex.
\end{itemize}
\end{dfn}

The cyclic  case relies on Hamilton's theorem, as in \cite{BoP}.
We reduce the dihedral  case to the cyclic one, by using a finite
covering argument. The platonic  case relies on a deformation
argument of spherical cone structures on $\OO$ given in Section
\ref{sec:sphericalcone}.

Previous to all these cases, in next subsection we prove that
$\pi_1\OO$ is finite.

\subsection{Nonnegative curvature and the fundamental group}

In this section, we relax the condition on cone angles and consider
Euclidean cone 3-manifolds with cone angles $\leq\al<2\pi$. In particular,
these are Alexandrov spaces of nonnegative curvature, and we will prove
the following result in the spirit of the Bonnet-Myers bounded diameter
theorem for Riemannian manifolds with lower positive curvature bound.

\begin{prop}\label{prop:finitegroup}
Let $X$ be a connected Euclidean cone 3-manifold with cone angles
$\leq\al<2\pi$ and non-empty singular set $\Si_{X}$. Suppose that
$\Ga$ is an infinite discrete group acting properly
discontinuously on $X$ by isometries. Then $\Ga$ is virtually
cyclic.

If $\Ga$ acts moreover cocompactly, then $X$ splits isometrically
as the product $Y\times\R$ of the real line with a closed
Euclidean cone surface $Y$.
\end{prop}
\begin{proof} {\em Step 1: There are no singular vertices.} Suppose that
$X$ contains singular vertices and consider a diverging sequence
$x_1,x_2,\dots$ of distinct vertices in the same $\Ga$-orbit. We
fix a base point $p\in X$ and denote by $v_i\in\L_p$ the direction
of the segment $\ol{px_i}$. We may assume without loss of
generality that the $v_i$ converge,
\begin{equation}
\label{dirconv} v_i\lra v_{\infty},
\end{equation}
and that $\frac{d_{i+1}}{d_i}\ra\infty$. Applying triangle
comparison (with Euclidean plane) to the triangles
$\De(p,x_i,x_{i+1})$ for large $i$, we obtain that
\[ \liminf_{i\ra\infty}\bigl(\angle_p(x_i,x_{i+1}) +
\angle_{x_i}(p,x_{i+1})\bigr)
 \geq\pi .\]
In view of (\ref{dirconv}), this means that
$\angle_{x_i}(p,x_{i+1})\ra\pi$. On the other hand, by the
Diameter Rigidity Theorem for CAT(1)-spaces, we have that
$diam(\L_{x_i})=d<\pi$, a contradiction. This shows that
$\Si_{X}^{(0)}$ is empty.

{\em Step 2: $\Si_{X}$ has finitely many connected components.}
The argument is similar. Assume that there are infinitely many
connected components $\si_i$ of $\Si_{X}$ and consider shortest
segments $\ol{px_i}$ from $p$ to $\si_i$. Let $w_i$ be the
direction of $\ol{px_i}$ at $x_i$. After passing to a subsequence,
we conclude as before that $\angle_{x_i}(p,x_{i+1})\ra\pi$. This
is absurd because, due to our upper bound on cone angles, we have
$rad(\L_{x_i}, w_i)\leq\max(\frac{\al}{2},\frac{\pi}{2})$ where
the radius $rad(\L_{x_i}, w_i)$ of $\L_{x_i}$ at $w_i$ is defined
as the Hausdorff distance of $\L_{x_i}$ and the one point subset
$\{w_i\}$. Thus $\Si_{X}$ has finitely many connected components.

{\em Step 3.} It follows that a finite index subgroup of $\Ga$
preserves one (each) singular component. The discontinuity of the
action implies that the components of $\Si_{X}$ are complete
geodesics and that $\Ga$ is virtually cyclic.

If the action of $\Ga$ is in addition cocompact then $X$ is
quasi-isometric to $\R$ and has therefore two ends. We apply the
Splitting Theorem, cf.\ \cite[Theorem 10.5.1]{burago}, to conclude
that $X$ splits as a metric product of $\R$ and an Alexandrov
space $Y$ of nonnegative curvature. Then $Y$ must be a closed
Euclidean cone surface. \end{proof}

\medskip
The first step in the proof of Theorem~\ref{thm:su} is   the
following lemma:

\begin{lem}\label{lem:pi1finite} Let $\OO$ be a closed orientable {\rm
irreducible} 3-orbifold. If there exists a closed Euclidean cone
structure on $\OO$ with cone angles strictly less than the
orbifold angles of $\OO$, then $\pi_1(\OO)$ is finite.
\end{lem}

\begin{proof} The Euclidean cone structure $X$ on $\OO$ (with cone angles
strictly less than the orbifold angles of $\OO$) lifts to a
Euclidean cone structure $\tilde X$ on the universal cover
$\widetilde \OO$. The Euclidean cone manifold $\tilde X$ has cone
angles $\leq \omega < 2\pi$, for some constant $0 < \omega< 2\pi$.
In addition, the fundamental group $\pi_{1}(\OO)$ acts
isometrically on $\tilde X$.

If $\pi_{1}(\OO)$ is infinite, then Proposition
\ref{prop:finitegroup} shows that $\pi_{1}(\OO)$ is virtually
cyclic and that $\tilde X$ splits as a metric product $\mathbb R
\times \tilde Y^2$, where $\tilde Y^2$ is a closed Euclidean cone
2-manifold. Since the action of $\pi_{1}(\OO)$ preserves this
product, $\tilde Y^2$ covers a totally geodesic surface $Y$ in
$X$. The cone surface $Y$ is a Euclidean cone structure on a
spherical turnover. This turnover is essential in $\OO$,
contradicting the irreducibility. \end{proof}

\subsection{The cyclic case}\label{subsec:cyclic}

Suppose that $\OO$ is of cyclic type.

The following lemma due to M.\ Feighn allows to apply
geometrization results for manifolds.

\begin{lem}[\cite{Fei}]\label{lem:verygood}
 If a closed orientable 3-orbifold of cyclic or dihedral type has finite
fundamental group, then it is very good.
\end{lem}

\begin{rem} A closed, orientable, irreducible very good
3-orbifold with  finite fundamental group is small. This is a
consequence of the equivariant Dehn Lemma (cf.
\cite{JR},\cite{MYOne}).
\end{rem}

We lift the Euclidean cone metric to  the universal covering of
$\OO$, which is a compact manifold denoted by $\widetilde \OO$.
Thus we have a $\pi_1(\OO)$-invariant Euclidean cone metric on
$\widetilde \OO$ with cone angles $< 2\pi$. This metric can easily
be desingularized to a $\pi_1(\OO)$-invariant smooth Riemannian
metric of non-negative sectional curvature,  because the singular
components are circles. More precisely, the singular locus is
locally a product of a singular disc with $\mathbb R$. This disc
is isometric to the neighborhood of the tip of a cone in Euclidean
space described by the equation $z= k\sqrt{x^2+y^2}, z\geq 0$. It
suffices to round the tip of the cone with a metric of
non-negative curvature invariant by rotations, that only depends
on the distance to the singular locus.

By Hamilton's theorem \cite{HamOne,HamTwo}, it follows that
$\widetilde \OO$ admits a $\pi_1(\OO)$-invariant smooth metric
locally modelled on $S^3$, $S^2\times\mathbb R$ or $\mathbb R^3$.
Only the spherical case is possible because $\pi_1(\OO)$ is
finite. Thus $\OO$ is spherical in the cyclic case.

\subsection{The dihedral case}

Suppose that $\OO$ is of dihedral type.

There exists a singular edge $e_0\subset\Sigma$ with the following
properties:
\begin{itemize}
\item[--] $e_{0}$ has two different endpoints, and
\item[--] the
branching index of every other edge of $\Sigma $ adjacent to $e_0$
is $2$.
\end{itemize}
To prove the existence of $e_0$, notice that any edge with label
$>2$ satisfies the properties. In addition, if all labels are
$=2$, there is always an edge with different endpoints.

\medskip
The covering provided in the following result will be useful for
the reduction to the cyclic case.

\begin{prop}\label{prop:virtuallycyclic}
There exists a finite regular covering
\[ p\co\hat \OO\to \OO \]
of orbifolds such that $\hat\OO$ is small of cyclic type and its
branching locus is $\hat \Sigma = p^{-1}(e_0)$. In addition $p$
preserves the ramification index of $e_0$.
\end{prop}

\begin{proof} Let $\OO'$ be the orbifold obtained from $\OO$ by removing
the open edge $e_0$ from the branching locus. (This change of the
orbifold structure amounts to putting the label $1$ on $e_0$.)
Since all edges of $\OO$ adjacent to an endpoint of $e_0$ have
label $2$, $\OO'$ is still an orbifold.

\begin{lem} The fundamental group of $\OO'$ is finite. \end{lem}

\begin{proof} We use the same Euclidean cone metric on $\OO$ as above and
consider it as singular metric of nonnegative curvature on $\OO'$.
We argue by contradiction as in the proof of
Lemma~\ref{lem:pi1finite}. If $\pi_1(\OO')$ were infinite,
$\widetilde{\OO'}$ would split metrically as a product $\mathbb
R\times Y'$. This is absurd because the metric on
$\widetilde{\OO'}$ has singular vertices. \end{proof}

Since $\OO'$ is of dihedral type and $\pi_1(\OO')$ is finite,
$\OO'$ is  very good by Lemma~\ref{lem:verygood}. The universal
covering of $\OO'$ induces a finite regular covering $p\!:\! \hat
\OO\to \OO$, where $\hat \OO$ is a closed orientable 3-orbifold
with underlying space the universal covering of $\OO'$,
 and
branching locus $\hat \Sigma = p^{-1}(e_0)$. Notice that $\hat
\Sigma$ is a finite collection of disjoint embedded circles, by
the choice of $e_0$. Therefore $\hat \OO$ is  of cyclic type.

\begin{lem}
$\hat\OO$ is small.
\end{lem}

\begin{proof} Since $\OO$ is irreducible and very good, by the
equivariant sphere theorem (cf. \cite{DD}, \cite{JR},
\cite{MYThree}) the universal cover $\widetilde\OO$ of $\OO$
(which is also the universal cover of $\hat \OO$) is an
irreducible 3-manifold.

Suppose that $F\subset\hat\OO$ is a spherical 2-suborbifold. Since
$\widetilde\OO$ is irreducible, $F$ bounds a ball quotient $Q$
which is a cyclic 3-suborbifold. The Smith conjecture \cite{MB}
implies that $Q$ is a discal suborbifold. (It also follows from
the orbifold theorem in the cyclic case whose proof we already
completed, respectively from Corollary~\ref{cor:finact} in the
Introduction.) Hence $\hat\OO$ is irreducible.
\begin{rem}
Irreducibility of $\hat \OO$ also follows from the irreducibility
of $\OO$ by an equivariant spherical 2-orbifold theorem
(\cite{Mai}, \cite{TY}), whose proof relies on PL least area
techniques for 2-orbifolds that generalize the notion of PL least
area surfaces introduced in \cite{JR}. However, for completeness,
we use here the fact that $\OO$ is very good with finite
fundamental group to give a direct argument.
\end{rem}

To see that $\hat\OO$ is small, suppose that $F \subset\hat\OO$ is
an essential 2-suborbifold. By irreducibility, $F$ cannot be
spherical or bad and therefore has infinite fundamental group.
Since $\pi_1(\hat \OO)$ is finite, $F$ lifts to a compressible
surface in the universal covering. The equivariant Dehn Lemma
implies that $F$ has a compressing disc. \end{proof}

This finishes the proof of Proposition~\ref{prop:virtuallycyclic}.
\end{proof}

We consider the compact 3-orbifold $\OO_0 = \OO-\NN(\Sigma_0)$
obtained by removing a regular neighborhood of the singular edge
$e_0$.

\begin{lem} The orbifold $\OO_0 = \OO-\NN(e_0)$ is
Haken and topologically atoroidal.
\end{lem}

\begin{proof}
 We first prove the irreducibility of
$\OO_0$. Let $S\subset \OO_0$ be a spherical 2-suborbifold. It
bounds a discal 3-suborbifold in $\OO$, since $\OO$ is
irreducible. If it does not bound a discal 3-orbifold in $\OO_0$,
then a neighborhood $\NN(e_0)$ is contained in the interior of a
discal 3-suborbifold of $\OO$. This is impossible, since $e_0$ is
a singular edge with two distinct vertices in $\Sigma$.

It is clear that $\OO_0$ does not contain any Euclidean or
hyperbolic turnover, because $\OO$ cannot contain  such a turnover
by smallness.

To see that $\OO_0$ is topologically atoroidal, suppose that
$T\subset\OO_0$ is a non-singular torus or a pillow which is
incompressible in $\OO_0$. Since $\OO$ is small, $T$ must be
compressible in $\OO$,  and the compression discal 2-suborbifold
must meet $e_0$. Hence, by irreducibility of $\OO$, $T$ bounds a
solid pillow containing $e_0$ and thus is parallel to
$\partial\OO_0$. The case that $T$ is a smooth torus cannot occur.
\end{proof}

Now we can apply Thurston's hyperbolization for Haken orbifolds,
cf.\ \cite[Ch. 8]{BoP}. It follows that $\OO_0$ is Seifert or
hyperbolic.

The Seifert case is quickly treated: $\OO$ is obtained from
$\OO_0$ by gluing a solid pillow $P$ to its boundary. If the
meridian of the pillow $\D\OO_0$ is homotopic to the fiber of
$\OO_0$ then the irreducibility of $\OO$ implies that $\OO_0$
contains no essential annulus. It follows that $\OO_0$ is a solid
pillow itself. Solid pillows admit many Seifert fibrations, and we
can modify the Seifert fibration so that the fiber is not a
meridian of the solid pillow $P$. Hence the Seifert fibration
extends to $\OO$. Since a Seifert fibred 3-orbifold with finite
fundamental group is spherical, $\OO$ is spherical, i.e.\ the
orbifold theorem holds in this case. Hence from now on we make the
following:

\begin{ass}
The orbifold $\OO_0$ admits on its interior a complete hyperbolic
structure of finite volume.
\end{ass}

We proceed now with the proof as in Section~\ref{sec:deform}
 by starting to
increase the cone angle along the singular edge $e_0$ and by
analyzing the degenerations.

We fix some notation. Let $m_0$ be the ramification index of the
edge $e_0$. For $t\in[0,1]$, let $X(t)$ denote a hyperbolic cone
structure on $\OO$, having the same prescribed cone angles as the
orbifold $\OO$ along the edges and circles of $\Sigma-e_0$ and
cone angle $\frac{2\pi}{m_0}t$ along the edge $e_0$. $X(0)$
denotes the complete hyperbolic structure of finite volume on the
interior of $\OO_0$.

In order to study the deformations of the hyperbolic cone
structure $X(t)$ while increasing $t$, we consider the variety of
representations $\Hom (\pi_{1}(\OO_0),PSL_2(\mathbb C))$ and the
variety of characters $\XX(\OO_0)$. As in Theorem
\ref{thm:bianalytic}, the (square of the) trace of the meridian
around $e_0$ is a local parameter for $\XX(\OO_0)$ near the complete
structure, cf.\ \cite[Theorem B.2.7]{BoP}. Hence the irreducible
component $\mathcal C$ of $\XX(\OO_0)$, that contains the holonomy
of $\OO$, is a curve.

As in Section \ref{sec:deform}, we define
    \[
    \II(\OO):=\left\{ t\in [0,1] \,\left\vert \begin{array}{l}
\text{there exists a hyperbolic cone structure on  }\OO \\
\text{with cone angle }\frac{2\pi t}{m_0}\text{ along } e_0\
\text{and cone angles equal} \\
\text{to the orbifold angles at all other edges,} \\
\text{with holonomy in }\mathcal C \text{ and volume } \leq v_0
\end{array}
\right. \right\}
    \]
where $v_0$ denotes the volume of the complete structure.

By hypothesis, $0\in\II(\OO)$, hence $\II(\OO)\neq\emptyset$.
Exactly as in Theorem \ref{thm:open} one proves that $\II(\OO)$ is
open to the right. Since $\pi_1(\OO)$ is finite, $\OO$ is not
hyperbolic and $1\not\in \II(\OO)$. Let $t_{\infty}:=\sup
\II(\OO)$. We have that $t_{\infty}\not\in  \II(\OO)$ by (right)
openness.

\begin{lem}\label{lem:cncollapses}
For any sequence $(t_n)$ in $\II(\OO)$ with $t_n\to t_{\infty}$,
the sequence of cone manifolds $(X(t_n))_{n \in \N}$ collapses.
\end{lem}

\begin{proof} Assume that $X(t_n)$ does not collapse. After choosing base
points $x_n$ in the thick parts of $X(t_n)$ and passing to a
subsequence, $(X(t_n),x_n)$ converges geometrically to a pointed
hyperbolic cone 3-manifold $(X_{\infty},x_{\infty})$ with finite
volume. $X_{\infty}$ cannot be compact because $t_{\infty}\not\in
\II(\OO)$.

We use the finite cover $p\!:\hat \OO\to \OO$ of
Proposition~\ref{prop:virtuallycyclic} and denote by $\hat
X(t)=p^{-1}(X(t))$ the lifted hyperbolic cone structure on $\hat
\OO$. The $\hat X(t_n)$ converge to a finite regular covering
$\hat X_{\infty}$ of $X_{\infty}$. Since all cone angles of each
$\hat X(t_n)$ are equal, the stability theorem
(Thm.~\ref{thm:stab}) applies and $\hat X_{\infty}$ is a
hyperbolic cone structure on $\hat\OO$. Now $\hat X_{\infty}$ is
not compact, and its ends are singular cusps which correspond to
singular vertices of $\hat\OO$. But $\hat\OO$ is of cyclic type,
contradiction. \end{proof}

Let $(t_n)$ be a sequence as above. We distinguish two cases
according to whether the sequence $\Diam(X(t_n))$ is bounded below
away from zero or not.

If $\Diam(X(t_n))\geq D>0$ for some uniform $D$, then the
fibration theorem   implies that $\OO$ is a Seifert fibred
3-orbifold. Since $\pi_1(\OO)$ is finite, it follows that $\OO$ is
spherical.

Otherwise, up to taking a subsequence, we can assume that
$\Diam(X(t_n)) \to 0$. 
Then we consider the rescaled sequence $\frac{1}{\Diam( X(t_n))}
X(t_n)$. If this rescaled sequence collapses, then the fibration
theorem still implies that $\OO$ is   Seifert fibred, and hence
spherical.

If the rescaled sequence does not collapse, then a subsequence converges
to  a closed Euclidean cone manifold $X(t_{\infty})$ homeomorphic to
$\OO$. We have $t_{\infty}<1$, because $\pi_1(\OO)$ is finite. Thus
$X(t_{\infty})$ lifts to a $\pi_1(\OO)$-invariant Euclidean cone metric on
the universal covering $\tilde{\mathcal \OO}$ with singular locus a link
and cone angle $t_{\infty}2\pi < 2\pi$. We conclude as in the cyclic case
that $\OO$ is spherical. \qed

\subsection{The platonic case}
\label{sec:plato}

Suppose that $\OO$ is of platonic type. The proof of
Theorem~\ref{thm:su} in this case is by induction on the number of
platonic vertices in the branching locus.

At each platonic vertex we have one singular edge with label 2 and
two edges with label $3,4$ or $5$. We fix a singular edge $e$ of
$\OO$ with label $n_e > 2$ such that at least one of its adjacent
vertices is platonic.

Let $\OO'$ be the orbifold obtained from $\OO$ by replacing the
branching index $n_e$ of  $e$
 by  2. We want to apply the induction hypothesis
to $\OO'$, because it has less platonic vertices than $\OO$. To do
it, we need the following lemma:

\begin{lem}\label{lem:o'small} The orbifold $\OO'$ is small.
\end{lem}

\begin{proof} By assumption,  $\OO$ and hence also $\OO'$ are closed. The
lemma follows from the fact that, for closed orbifolds, smallness
is a property independent of the labels of the branching locus.
\end{proof}

Notice that $\OO'$ has orbifold angles greater than or equal to
those of $\OO$ because one label has decreased. Thus the Euclidean
cone structure on $\OO$ given in the hypothesis of Theorem
\ref{thm:su}, when viewed on $\OO'$ still has cone angles strictly
less than the orbifold angles. It follows from the induction
hypothesis that $\OO'$ is spherical.

The induction step, and hence the conclusion of the proof of
Theorem \ref{thm:su}, is due to the following result which we will
prove in Section \ref{sec:defospace}.

\begin{prop}
\label{spericalcanbedeformed} The spherical structure on $\OO'$
can be deformed, through a continuous family of spherical cone
metrics, to a spherical structure on $\OO$.
\end{prop}

\section{Deformations of spherical cone structures}
\label{sec:sphericalcone}

This section is devoted to proving
Proposition~\ref{spericalcanbedeformed}.
 Some
preliminaries on varieties of representations are required.

\subsection{The variety of representations into $SU(2)$}

Let $\Gamma$ be a finitely generated group. The variety of
representations $\Hom(\Gamma,SU(2))$ is compact. The group $SU(2)$
acts on   $\Hom (\Gamma,SU(2))$ by conjugation, and the quotient
$$ \XX(\Gamma,SU(2))= \Hom (\Gamma,SU(2))/SU(2).$$
is also compact, but in general not algebraic.

\paragraph{The variety of characters.} The action by conjugation of
$SL(2,\mathbb C)$ on the variety of representations
$\Hom(\Gamma,SL(2,\mathbb C))$ is algebraic. The quotient
$$ \XX(\Gamma,SL(2,\mathbb C))=\Hom(\Gamma,SL(2,\mathbb C)) /\!/
SL(2,\mathbb C) $$ of this action provided by geometric invariant theory
carries a natural structure as an affine algebraic subset of $\mathbb C^N$
defined over $\mathbb Q$ \cite{MoS}, points can be interpreted as
characters of representations $\Gamma\ra SL(2,\mathbb C)$. The conjugacy
class of a representation $\Gamma\ra SU(2)$ is determined by its
character, and we have a natural inclusion
\begin{equation}
\label{natinc} \XX(\Gamma,SU(2)) \subseteq \XX(\Gamma,SL(2,\mathbb C))
.
\end{equation}

\paragraph{Trace functions.} For $\gamma\in\Gamma$, the trace function
\begin{eqnarray*}
 \Hom(\Gamma,SL(2,\mathbb C))
&\to&\mathbb C\\ \rho \quad &\mapsto&\operatorname{trace}
(\rho(\gamma))
\end{eqnarray*}
induces an algebraic function $$I_{\gamma}\!:\!
\XX(\Gamma,SL_2(\mathbb C))\to\mathbb C.$$

\paragraph{The ambient space $\mathbb C^N$.}
The embedding of $ \XX(\Gamma,SL(2,\mathbb C))$ into $\mathbb C^N$
is realized by taking as coordinates the functions $I_{\gamma}$,
where $\gamma$ runs through the words of length at most three in
the generators of $\Gamma$ \cite{GM}.

Since traces of matrices in $SU(2)$ are real, we have the
embedding
$$
\XX(\Gamma,SU(2))\subset \XX(\Gamma,SL_2(\mathbb C))\cap \mathbb R^N.
$$ of $\XX(\Gamma,SU(2))$ into a real algebraic variety. One can
show that $\XX(\Gamma,SL_2(\mathbb C))\cap \mathbb
R^N=\XX(\Gamma,SU(2))\cup \XX(\Gamma,SL_2(\mathbb R))$ \cite{MoS}. The
following lemma implies that $\XX(\Gamma,SU(2))$ is locally $\mathbb
R$-algebraic at characters of non-abelian representations.

\begin{lem}[{\cite[lemme 5.25]{PoOne}}]
\label{lem:SU(2)} Let $[\rho]\in \XX(\Gamma,SU(2))$ be the conjugacy
class of a nonabelian representation. There exists an open
neighborhood $U\subset \mathbb R^N $ of $[\rho]$ such that: $$
\XX(\Gamma,SU(2)) \cap U= \XX(\Gamma,SL_2(\mathbb C))\cap \mathbb
R^N\cap U. $$
\end{lem}

\subsection{Lifts of holonomy representations into $SU(2)\times SU(2)$
and spin structures}

We recall that an element $(A,B)\in SU(2)\times SU(2)\cong
Spin(4)$ acts on $SU(2)\cong S^3$ by
\[ x\mapsto A x B^{-1}.\]
In particular, the kernel of the action is the order two subgroup
$\{\pm(\operatorname{Id}, \operatorname{Id})\}$. The following
lemma is classical:

\begin{lem}\label{lem:rotations} The element $(A,B)\in SU(2)\times SU(2)$
acts on $SU(2)\cong S^3$ as a rotation (i.e.\ has fixed points)
iff $\trace (A)=\trace (B)$. Trace and rotation angle are related
by
\begin{equation}
\label{traceangle} \trace(A)=\pm2\cos(\frac{\alpha}{2}) .
\end{equation}
\end{lem}

One can view $SO(4)$ as the frame bundle on $S^3=SO(4)/SO(3)$. The
unique spin structure on $S^3$ is given by the canonical
projection $Spin(4)\ra SO(4)$. Given a spherical manifold $N$, not
necessarily complete, with holonomy representation
$\rho:\pi_1(N)\ra SO(4)$, the spin structures on $N$ correspond to
lifts of $\rho$ to a representation into $Spin(4)$. (The lift can
be obtained from a spin structure by developing it, using the
natural connection, onto $Spin(4)$ and taking its holonomy.)

Let $p_i \co SU(2)\times SU(2)\to SU(2)$ denote the projection to
the $i$-th factor, for $i=1,2$.

\begin{lem}\label{lem:angles} Let
$ \phi:\pi_1(M)\to SU(2)\times SU(2)$ be the lift of the holonomy of a
spherical cone manifold. Then both $p_1\circ\phi$ and $p_2\circ\phi$ are
non-abelian unless $\Sigma $ is a link and $M=X-\Sigma$ is Seifert fibred.
\end{lem}

\begin{proof} Assume for instance that $p_1\circ\phi$ is abelian. This
means that the image of $p_1\circ\phi$ is contained in a
diagonalizable subgroup $\cong S^1$. Therefore $\phi$ preserves
the corresponding Hopf fibration on $S^3$. It follows easily from
this that $\Sigma $ is a link and $M$ is Seifert fibred (cf.
\cite[Lemma~9.1]{PoTwo}). \end{proof}

Let $\mu_1,\ldots\mu_q$ be meridians for the singular edges and
circles in $\Sigma$.

\begin{thm}[Local parametrization]\label{thm:LR} Let $\OO$ be a spherical
orbifold such that $M=\OO-\Sigma$ is not Seifert fibred. If
$\phi:\pi_1(M)\to SU(2)\times SU(2)$ is a lift of the holonomy,
then both $[p_1\circ\phi]$ and $[p_2\circ\phi]$ are smooth points
of $\XX(M,SU(2))$. Moreover,
$$ (I_{\mu_1},\ldots,I_{\mu_q}) : \XX(M,SU(2))\to
\mathbb R^q $$ is a local diffeomorphism at both points
$[p_1\circ\phi]$ and $[p_2\circ\phi]$.
\end{thm}

This theorem is an infinitesimal rigidity result for spherical
orbifolds, and its proof is postponed   to the last subsection
\ref{sec:prooflocrig}. It will be obtained from a cohomology
computation, using the fact that spherical orbifolds are finitely
covered by $S^3$.

\subsection{Deformation space of spherical structures}
\label{sec:defospace}

The aim of this section is to prove Proposition
\ref{spericalcanbedeformed}.

We adopt the notation of Section \ref{sec:plato}. Let $\mu_e$
denote the meridian of the edge $e$. We recall that $\XX(M,SU(2))$
is contained in the real algebraic set $V=\XX(M,SL_2(\mathbb
C))\cap \mathbb R^N$, cf.\ (\ref{natinc}). Let $\phi_0$ be the
lift to $SU(2)\times SU(2)$ of the holonomy representation of
$\OO'$ corresponding to the choice of spin structure on $M$. Its
conjugacy class $[\phi_0]$ is contained in
\[ \XX(M,SU(2)\times SU(2))=\XX(M,SU(2))\times \XX(M,SU(2))\subseteq V\times V .\]
According to Theorem \ref{thm:LR}, $\XX(M,SU(2)\times SU(2))$ is
locally bianalytic to $\mathbb R^{2q}$ at $[\phi_0]$ and, due to
Lemma \ref{lem:SU(2)}, is a neighborhood of $[\phi_0]$ in $V\times
V$.

Let  $p_1:SU(2)\times SU(2)\to SU(2)$ denote the projection to the
first factor. Consider the algebraic subset
\begin{equation*}
 W =\left\{([\rho_1],[\rho_2])\in V\times V\ \left\vert
\begin{array}{l}
I_{\mu_i}([\rho_1])= I_{\mu_i}([\rho_2]), \text{ for each meridian
} \mu_i
\\
I_{\mu_i}([\rho_1]) = I_{\mu_i}([ p_1(\phi_0)]), \text{ for each
meridian } \mu_i\neq\mu_e
\end{array}
\right.\right\}
\end{equation*}
By Lemma~\ref{lem:rotations}, these equations are the algebraic
conditions for a representation of $\pi_1(M)$ in $SU(2)\times
SU(2)$ to be the lift of a representation in $SO(4)$ with the
properties:
\begin{itemize}
\item[--] the images of the meridians are rotations;
\item[--] the rotation angles of all meridians are fixed
except for the meridian $\mu_e$ of the edge $e$.
\end{itemize}

Let $W_0$ be the irreducible component of $W$ containing
$[\phi_0]$. By Theorem \ref{thm:LR}, $W_0$ is a real algebraic
curve and $I_{\mu_e}\circ p_1$ restricted to $W_0$ is a smooth
local parameter near $[\phi_0]$. In particular it is non-constant
on $W_0$.

A neighborhood of $[\phi_0]$ in $W_0$ can be lifted to a curve of
representations $\pi_1(M)\ra SU(2)\times SU(2)$ which are lifts of
holonomies of spherical cone structures on $\OO'$. The fact that
the trace $I_{\mu_e}\circ p_1$ is a smooth local parameter on
$W_0$ near $[\phi_0]$ implies that the cone angle at $e$ is a
(continuous) local parameter for the family of spherical cone
structures near the orbifold structure. It takes values in a
neighborhood of $\pi$.

We take $\mathcal S^{\pm}$ to be the connected component of the
semi-algebraic set
$$ \{([\rho_1],[\rho_2])\in W_0\mid 0 \leq \pm I_{\mu_e}([\rho_1])\leq
2\cos(\pi/n_e)\} $$ that contains $[\phi_0]$. One of these two
sets, say $\mathcal S^+$, contains representations arising from
cone metrics with cone angle $<\pi$ at $e$.

\begin{lem}
\label{lem:arehol} All representations in $\mathcal S^+$ are lifts
of holonomy representations for spherical cone structures on
$\OO'$ with cone angle at $e$ in $[\frac{2\pi}{n_e},\pi]$.
\end{lem}

\begin{proof} Let $A$ be the subset of representations in $\mathcal S^+$
that are such lifts. By our previous discussion, $A$ contains a
neighborhood of the endpoint $[\phi_0]$ of $\mathcal S^+$.

{\em Openness of $A$.} Lemmas~\ref{lem:SU(2)} and~\ref{lem:angles}
imply that perturbations of representations in $A$ still take
values in $SU(2)\times SU(2)$. Moreover, perturbations of holonomy
representations are induced by perturbations of cone structures,
cf. \cite{CEG,PoTwo}.

{\em Closedness of $A$.} Consider a sequence of points
$[\phi_n]\in A$ converging to $[\phi_{\infty}]\in \mathcal S^+$.
We have to show that $[\phi_{\infty}]\in A $.

For $n\in\mathbb N$, let $X_n$ be the spherical cone manifold with
holonomy lift $\phi_n$. The $X_n$ are Alexandrov spaces with
curvature $\geq1$ and therefore have diameter $\leq\pi$. Hence if
the sequence $(X_n)_{n\in\N}$ does not collapse then, up to a
subsequence, it geometrically converges to a spherical cone
manifold $X_{\infty}$ with the same topological type and holonomy
lift $\phi_{\infty}$. Thus $[\phi_{\infty}]\in A$ in this case.

Assume now that the sequence $(X_n)$ collapses. In the spherical
case, this is equivalent to $\operatorname{vol}(X_n)\ra0$, because
$$
\operatorname{vol}(X_n)\leq 2 \pi\, \operatorname{inj} (X_n).
$$
This formula can be proved by using the Dirichlet domain centered
at some point $x\in X_n$ of smallest injectivity radius, since
this domain embeds isometrically in a lens in $S^3$ of with
$2\operatorname{inj}(x)$.

 Denote by $\al_n$ the cone angle of $e$ for
the cone structure $X_n$. We may assume that $(\al_n)$ converges
and distinguish two cases.

{\em Case 1: $\al_n\ra\pi$.} We will apply Schl\"afli's formula,
relying on the algebraic structure of our deformation space.
$\mathcal S^+$ is contained in the algebraic curve $W_0$. Hence,
after passing to a subsequence, we may assume that the $[\phi_n]$
lie on an analytic path with endpoint $[\phi_{\infty}]$. Since the
trace $I_{\mu_e}\circ p_1$ is a non-constant analytic function on
$W_0$, its critical points do not accumulate and we can assume
that it is monotonic on this path. So the cone angle at $e$ is
monotonically increasing towards $\pi$. Schl\"afli's formula then
implies that the volume of the cone structure {\em increases}, as
we approach $[\phi_{\infty}]$. This contradicts collapse.

{\em Case 2: $\al_n\ra\al_{\infty}<\pi$.} Since at least one of the
singular vertices in $\OO$ adjacent to $e$ is platonic each $X_n$ has at
least one vertex with two cone angles uniformly bounded away from $\pi$.
Since $(X_n)$ collapses by assumption, the thick vertex lemma
(Lemma~\ref{lem:thickvertex}) implies that $\operatorname{diam}(X_n)\to
0$.

Now we rescale $X_n$ by the (inverse of the) diameter. Applying the thick
vertex lemma again, the rescaled sequence converges to a compact Euclidean
cone manifold, which is a Euclidean structure on $\OO$ with cone angles
greater than or equal to the orbifold angles of $\OO$.

\begin{sublem}
\label{sublem:cat0}
This Euclidean cone structure on $\OO$ has
curvature $\leq0$ in the orbifold sense (i.e. it is locally the
quotient of a CAT(0) space by a finite group of isometries).
\end{sublem}

\begin{proof}
By Gromov's criterion \cite{Bal}, all we have to check is that the
link of each point is CAT(1) in the orbifold sense (i.e. a
quotient of a CAT(1) space by the isotropy group). The links are
quotients of the unit sphere except at the points of $e$. The
CAT(1) property for the links of interior points of $e$ is clear
because the cone angle is greater than or equal to the orbifold
angle. It is clear for the same reason at a dihedral endpoint of
$e$. (There may be at most one.) At platonic endpoints of $e$ it
follows from Proposition~\ref{prop:polarofplatonic} if $e$ has
different endpoints, and from its addendum if $e$ is a loop, to be
proved in next subsection.
\end{proof}

It follows from the sublemma and from Haefliger's version of the
Cartan-Hadamard theorem \cite{Hae} that $\pi_1(\OO)$ is infinite.
We obtain a contradiction, because $\pi_1(\OO)$ is finite by
Lemma~\ref{lem:pi1finite}. Thus $(X_n)$ does not collapse. This
finishes the proof of Lemma \ref{lem:arehol}. \end{proof}

\begin{lem}
 \label{lem:coneanglesurj}
 The map
\[ \al:\mathcal S^+\lra\left[\frac{2\pi}{n_e},\pi\right] \]
given by the cone angle at $e$ is surjective.
\end{lem}

\begin{proof} Recall that $\mathcal S^+$ is a closed connected subset of
the real algebraic curve $W_0$. Moreover, $\mathcal S^+$ is
compact because it is contained in the subset of conjugacy classes
of $SU(2)\times SU(2)$-valued representations. As an algebraic
curve, $W_0$ is homeomorphic to a graph with finitely many
vertices. It follows that $\mathcal S^+$ is a compact graph.

Consider the subset
\[ \mathcal S^+_0 :=\al^{-1}\left((\frac{2\pi}{n_e},\pi)\right) \subset
\mathcal S^+ .\] The complement $\mathcal S^+-\mathcal S^+_0$ is
finite, because the nonconstant analytic trace function has
discrete level sets.

Since $\mathcal S^+_0$ is locally algebraic, Sullivan's local
Euler characteristic theorem \cite{Sul} implies that vertices of
$\mathcal S^+_0$ have even valency. It follows that $\mathcal
S^+_0$ is a (noncompact) graph with an even number of ends. Recall
that a neighborhood of $[\phi_0]$ in $\mathcal S^+$ is a curve
initiating in $[\phi_0]$, and therefore precisely one end of
$\mathcal S^+_0$ converges to $[\phi_0]$. As a consequence, there
exists another end of $\mathcal S^+_0$ converging to a point
$[\phi_1]\in \mathcal S^+-\mathcal S^+_0$ different from
$[\phi_0]$. We have that
\[ \al([\phi_1])\in\{\frac{2\pi}{n_e},\pi\} .\]
We are done if $\al([\phi_1])=\frac{2\pi}{n_e}$. If
$\al([\phi_1])=\pi$, we obtain a contradiction from the following
result.

\begin{thm}[de Rham, cf.\ \cite{rothenberg}]
\label{prop:deRham} A spherical structure on a closed orientable
connected smooth 3-orbifold is unique up to isometry.
\end{thm}

Namely, since a spin structure has been fixed on $M$, the
non-conjugate representations $\phi_0$ and $\phi_1$ correspond to
non-isometric spherical structures on the orbifold $\OO'$.

This concludes the proof of   Lemma~\ref{lem:coneanglesurj}.
\end{proof}

Proposition \ref{spericalcanbedeformed} is a direct consequence of
the results in this section.

\subsection{Certain spherical cone surfaces with the CAT(1)
property}

The results of this subsection have been used in
Sublemma~\ref{sublem:cat0}  to show that a certain Euclidean cone
structure on the orbifold $\OO$ is a metric that satisfies locally
the CAT(0) property, by showing that the links are CAT(1).

Let $\L$ be a spherical cone surface with cone angles $>2\pi$.
Such a surface has curvature $\leq1$ in the {\em local} sense.
We discuss some examples when
$\L$ has the CAT(1) property, i.e.\
satisfies {\em global} triangle comparison with {\em upper}
curvature bound $1$.

Due to a general criterion for piecewise spherical complexes, $\L$
is CAT(1) if and only if it contains no nonconstant closed
geodesic with length $<2\pi$,
cf. \cite[sec.\ 4.2 A and B]{hypgroups}.
An elegant way to check the CAT(1) property was
discovered by Rivin in his thesis:

\begin{thm}
[Rivin, cf.\ \cite{rivinhodgson}] The polar dual of a compact
convex polyhedron in $\H^3$ is CAT(1).
\end{thm}

The result extends to ideal polyhedra, cf.\ \cite[Thm.\
4.1.1]{CharneyDavis}.

The {\em polar dual} or {\em Gau\3 image} $G(P)$ of a convex
polyhedron $P$ in $\H^3$ is constructed as follows, generalizing
the Gau\3 map for convex polyhedra in Euclidean space. For every
vertex $v$ of $P$ we take the set $G(v)$ of all outer unit normal
vectors at $v$. Equipped with its natural metric as a subset of
the unit tangent sphere, $G(v)$ is a spherical polygon. The sides
of $G(v)$ correspond to the edges of $P$ adjacent to $v$. If the
vertices $v_1$ and $v_2$ are joined by an edge $e$, we glue the
polygons $G(v_1)$ and $G(v_2)$ along their sides corresponding to
$e$. The resulting spherical complex is $G(P)$, and it is easily
seen to be a spherical cone surface with cone points of angles
$>2\pi$.

\medskip
As an example, relevant in Sublemma~\ref{sublem:cat0}, we
determine the polar dual of the platonic solids. Let $P$ be a
regular polyhedron in $\H^3$. ($P$ can be a tetrahedron with face
angles $<\frac{\pi}{3}$, a cube with face angles $<\frac{\pi}{2}$,
an octahedron with face angles $<\frac{\pi}{3}$, a dodecahedron
with face angles $<\frac{3\pi}{5}$, or an icosahedron with face
angles $<\frac{\pi}{3}$.) The isometry group $Isom(P)$ acts on the
dual $G(P)$ as a reflection group. Let $\De(P)$ be the quotient
2-orbifold. It is a triangle with reflector boundary equipped with
a spherical metric. Each one of the three vertices of $\De(P)$
corresponds respectively to the vertices, edge midpoints and face
centers of $P$. The angles of the triangle equal the orbifold
angles at the first two vertices, but the third angle is bigger
than the corresponding orbifold angle. In other words, in the
orbifold sense, the metric is smooth everywhere except at the
third vertex, and there it has concentrated negative curvature.

For instance, if $P$ is a (possibly ideal) icosahedron with face
angles $<\frac{\pi}{3}$, then $G(P)$ is a piecewise spherical
dodecahedron composed of regular pentagons with angles
$\in(\frac{2\pi}{3},\pi]$. The orbifold $\De(P)$ has angles
$\frac{\pi}{5},\frac{\pi}{2}$ and third angle
$\in(\frac{\pi}{3},\frac{\pi}{2}]$. Going through the various
cases yields the following examples of piecewise spherical metrics
on $S^2$ with the CAT(1) property.

\begin{prop}
\label{prop:polarofplatonic} Let $\De$ be a spherical 2-orbifold which is
topologically a triangle with reflector boundary. Suppose that
$\De$ is equipped with a spherical metric so that the boundary is
geodesic, all angles are $\leq\frac{\pi}{2}$, and the metric is
smooth everywhere except at one vertex where it has concentrated
negative curvature (in the orbifold sense).

Then the pull-back of the metric to the universal covering
orbifold $\widetilde\De\cong S^2$ satisfies the CAT(1) property.
\end{prop}

Only the case when $\De$ is cyclic is not given by the platonic
solids. But in this case, $\widetilde\De$ is the suspension of a
circle with length $>2\pi$, hence also CAT(1).

The next example is related.

\begin{add}
Let $\De$ be as in the proposition 
and with local isotropy groups $D_2,D_3$ and $D_3$ at the
vertices. Suppose that the spherical metric has angle
$\frac{\pi}{2}$ at the $D_2$-vertex and equal angles
$\in(\frac{\pi}{3},\frac{\pi}{2}]$ at the two vertices with
$D_3$-isotropy. Then the same conclusion holds.
\end{add}
\begin{proof} By folding these orbifolds along their symmetry axis, one
obtains index-two ramified coverings over the orbifolds with
vertex isotropies $(D_2,D_3,D_4)$. The angles of the quotient
orbifolds satisfy the assumptions in the proposition. The
quotients have the same universal cover, and the assertion follows
from the previous proposition. \end{proof}

\subsection{Proof of the local parametrization theorem}
\label{sec:prooflocrig}

\begin{proof}[Proof of Theorem~\ref{thm:LR}] Let $\rho=p_i\circ \phi\co
\pi_1(M)\to SU(2)$, for $i=1$ or $2$. The proof has three main
steps. In Step~1 we show that the dimension of the Zariski tangent
space of $\XX(M,SU(2))$ at $[\rho]$ equals the number of meridians.
In Step~2 we prove that $[\rho]$ is a smooth point of
$\XX(M,SU(2))$. In Step~3 we check that the differential forms $\{
dI_{\mu_1},\ldots,dI_{\mu_q} \}$ are a basis for the cotangent
space.

\medskip
\emph{Step 1: Computation of the Zariski tangent space of
$\XX(M,SU(2))$ at $[\rho]$.} The Zariski tangent space is given by
the cohomology group $H^1(M;Ad\circ\rho)$ where we work with
coefficients in the lie algebra $su(2)$ twisted by the adjoint
representation $Ad\circ\rho$.

First compute the cohomology $H^*(\OO;Ad\circ\rho)$. Notice that
$Ad\circ\rho:\pi_1(M)\ra Aut(su(2))$ factors through
$\pi_1(M)\ra\pi_1(\OO)$ because we may compose the holonomy
representation
of the spherical structure on $\OO$ (which in general does not
lift to $Spin(4)$) with the (first or second) projection $SO(4)\ra
PSU(2)$ and $Ad$:
\[ \pi_1(\OO)\lra SO(4) \lra PSU(2) \lra Aut(su(2)) \]
We use simplicial homology. Let $K$ be a triangulation of the
underlying space of $\OO$ compatible with the branching locus and
let $\tilde K$ be its lift to the universal cover $\tilde \OO\cong
S^3$. We consider the following chain and cochain complexes:
\begin{eqnarray*}
 C_*(K;Ad\circ\rho)&=& su(2)\otimes_{\pi_1(\OO)} C_*(\tilde K;\mathbb
Z)  \\ C^*(K;Ad\circ\rho)&=& \Hom_{\pi_1(\OO)}(C_*(\tilde
K;\mathbb Z), su(2))
\end{eqnarray*}
 The homology of $C_*(K;Ad\circ\rho)$ is denoted by $H_*(\OO;Ad\circ\rho)$ and
the cohomology of $C^*(K;Ad\circ\rho)$ by $H^*(\OO;Ad\circ\rho)$.
From the differential point of view, $H^*(\OO;Ad\circ\rho)$ is the
cohomology of the $su(2)$-valued differential forms on $S^3$ which
are $\pi_1(\OO)$-equi\-variant.

\begin{lem}\label{lem:HO=0} $H_*(\OO;Ad\circ\rho)\cong 0$. \end{lem}

\begin{proof} There is a canonical projection $C_*(S^3,su(2))\ra
C_*(\OO,Ad\circ\rho)$ and, since $\pi_1(\OO)$ is finite, an averaging map
$C_*(\OO,Ad\circ\rho)\ra C_*(S^3,su(2))$ which is a section for the
projection. According to \cite[Prop.~10.4 in Ch.\ 3]{Brown}, the induced
map
\[ H^*(\OO,Ad\circ\rho) \lra \left(H^*(S^3,su(2))\right)^{\pi_1(\OO)}
 \cong H^*(S^3,\mathbb R) \otimes su(2)^{\pi_1(\OO)}
\]
is an isomorphism. The Lie algebra $su(2)$ does not have
nontrivial elements invariant by $Ad\circ\rho$, because $\rho$ is
non-abelian (Lemma~\ref{lem:angles}). \end{proof}

Let $\NN(\Sigma)$ be a tubular neighborhood of $\Sigma$ and let
$N=\OO-\NN(\Sigma)$. We remark that $N$ is compact and that the
inclusion $N\subset M$ is a homotopy equivalence, so
$H_*(M;Ad\circ\rho)\cong H_*(N;Ad\circ\rho)$.

From the previous lemma and by applying Mayer-Vietoris to the pair
$(N,\NN(\Sigma))$, it follows that there is a natural isomorphism
\begin{equation}
\label{mayervietoris} H^*(N; Ad\circ\rho)\oplus
H^*(\NN(\Sigma);Ad\circ\rho) \cong H^*(\partial N; Ad\circ\rho)
\end{equation}
induced by inclusion.

Consider the following piece of the exact sequence of the pair
$(N,\partial N)$:
\begin{equation}\label{eqn:pair}
0\lra H^1(N;Ad\circ\rho) \overset{i^*}\lra H^1(\partial
N;Ad\circ\rho) \overset{\delta}\lra H^2(N,\partial N;Ad\circ\rho)
\end{equation}
The injectivity of $i^*$ comes from (\ref{mayervietoris}).
Poincar\'e duality implies that $H^1(N;Ad\circ\rho)$ and
$H^2(N,\partial N;Ad\circ\rho)$ are dual, $H^1(\partial
N;Ad\circ\rho)$ is dual to itself, and moreover $\delta$ and $i^*$
are dual maps. Hence $\delta$ is surjective and
\begin{equation}
\label{dimH1M=1/2}
 \dim H_2(N,\partial N;Ad\circ\rho)=\dim H_1(N;Ad\circ\rho) =
\frac12 \dim H_1(\partial N;Ad\circ\rho).
\end{equation}

\begin{prop}\label{cor:dimH1} If $\Sigma$ has $q$ edges and circles, then $\dim
H^1(M;Ad\circ\rho)=q$.
\end{prop}
\begin{proof} We have to show that $\dim H^1(\partial N;Ad\circ\rho)=2q$. Since
$\partial N=\partial \NN(\Sigma)$ and the homology is the direct
sum of the homology of the connected components, it suffices to
compute $\dim H_1(\partial \NN(\Sigma);Ad\circ\rho)$ assuming that
$\Sigma$ is connected.

If $\Sigma$ is a circle, then $\partial \NN(\Sigma)$ is a torus
and $\rho\vert_{\pi_1(\partial \NN(\Sigma))}$ is abelian and
non-trivial. Thus $H^0(\partial \NN(\Sigma);Ad\circ\rho)\cong
su(2)^{Ad\circ\rho}\cong\mathbb R$. By duality, $\dim H^2(\partial
\NN(\Sigma);Ad\circ\rho)=1$. Since the torus has zero Euler
characteristic, it follows that $\dim H^1(\partial
\NN(\Sigma);Ad\circ\rho)=2$.

When $\Sigma$ is a trivalent graph with $v$ vertices, then it has
$q=3v/2$ edges. Note that the restriction
$\rho\vert_{\pi_1(\partial \NN(\Sigma))}$ is irreducible. (Namely,
the holonomy lift of the neighborhood of a vertex takes values in
the stabilizer of a point $\cong SU(2)$ and has irreducible image.
The stabilizer is a diagonal subgroup in $Spin(4)\cong SU(2)\times
SU(2)$ and projects isomorphically onto both factors.) Thus $
H^0(\partial \NN(\Sigma);Ad\circ\rho)\cong
su(2)^{Ad\circ\rho}\cong 0$. By duality, $H^2(\partial
\NN(\Sigma);Ad\circ\rho)\cong 0$. Since $\chi(\partial
\NN(\Sigma))= 2\chi(\Sigma)= 2(v-q)= -v=-\frac{2q}3$, we get $
\dim H^1(\partial \NN(\Sigma);Ad\circ\rho)=- \chi(\partial
\NN(\Sigma))\cdot \dim( su(2))=2 q$.
\end{proof}

\medskip
\emph{Step 2: Smoothness of $\XX(M,SU(2))$ at $[\rho]$.} The
cohomology group $H^1(M;Ad\circ\rho)$ is naturally identified with
the Zariski tangent space of $\XX(M,SU(2))$ at $[\rho]$. For an
element in $H^1(M;Ad\circ\rho)$, there is an infinite sequence of
obstructions to be integrable (cf. \cite{Dou}, \cite{Ab},
\cite{HPS},  \cite{Gol}). These obstructions are natural and live
in the second cohomology group.

Note that $H^2(\NN(\Sigma);Ad\circ\rho)=0$ because the orbifold
$\NN(\Sigma)$ is finitely covered by a manifold homotopy
equivalent to a 1-dimensional complex. Hence (\ref{mayervietoris})
provides an isomorphism
$$ H^2(M;Ad\circ\rho)\cong H^2(\partial N;Ad\circ\rho) .$$
The obstructions to integrability vanish for $\D N$ because $\D N$
is a surface (\cite{Gol} and \cite{HPS}). Naturalness of the
obstructions then implies that they also vanish for $M$. Thus
every element in $H^1(M;Ad\circ\rho)$ is formally integrable as a
power series, and by a theorem of Artin \cite{Ar} it is actually
integrable. Since $H^1(M;Ad\circ\rho)$ has dimension $q$, it
follows that $[\rho]$ is a smooth point of $\XX(M,SU(2))$ with local
dimension $q$.

\medskip
\emph{Step 3: A basis for $T^*_{[\rho]}\XX(M,SU(2))$.} By Step 2,
there is a natural isomorphism
\[ T^*_{[\rho]}\XX(M,SU(2)) \cong H_1(M;Ad\circ\rho).\]

\begin{lem}\label{prop:basis} The set $ \{ d I_{\mu_1},\ldots,d I_{\mu_q} \}$ is a
basis for $H_1(M;Ad\circ\rho)$.
\end{lem}
\begin{proof} As before we identify $H_*(M;Ad\circ\rho)\cong
H_*(N;Ad\circ\rho)$. Dual to (\ref{mayervietoris}), we have the
isomorphism
$$ H_1(\partial N; Ad\circ\rho)\cong H_1(N;
Ad\circ\rho)\oplus H_1(\NN(\Sigma);Ad\circ\rho). $$ The map
$\iota_*:H_1(\partial N; Ad\circ\rho)\ra
H_1(\NN(\Sigma);Ad\circ\rho)$ of cotangent spaces is induced by
the restriction map $\XX(\NN(\Sigma),SU(2))\ra \XX(\D N,SU(2))$. Since
the meridians have finite order in $\pi_1(\NN(\Sigma))$, the trace
functions $I_{\mu_j}$ are constant on $\XX(\NN(\Sigma),SU(2))$.
Hence $\iota_*(d I_{\mu_j})=0$. Since $\dim H_1(N;Ad\circ\rho)=q$
by Prop.~\ref{cor:dimH1}, the proof reduces to the following
lemma:

\begin{lem}\label{lem:LI}
The differential forms $\{ d I_{\mu_1},\ldots,d I_{\mu_q}\}$ are
linearly independent in $H_1(\partial N;Ad\circ\rho)$.
\end{lem}

\proof For each meridian $\mu_j$, we construct a deformation of
the restriction
  $\rho\vert_{\D N}$
that is parameterized by the trace functions $I_{\mu_j}$ and
leaves all other meridians constant. We proceed on each component,
assuming that the deformation on the other components is trivial.

First we consider the case  that $\mu$ is a meridian around a singular
circle, so the corresponding boundary component $T^2\subset \partial N$ is
a torus. We wish to deform $\rho|_{T^2}$. We choose $\lambda$ so that
$\lambda$ and $\mu$ generate $\pi_1(T^2)$. The elements $\rho(\mu)$ and
$\rho(\lambda)$ commute and can therefore be simultaneously diagonalized.
Thus $\rho(\mu)$ and $\rho(\lambda)$ can be varied independently inside a
circle subgroup of $SU(2)$. The matrix $\rho(\mu)$ has eigenvalues $\pm
e^{\pm i\frac{\alpha}{2}}$ and trace $\pm2\cos\frac{\alpha}{2}$ where
$\alpha$ is the cone angle, cf.\ (\ref{traceangle}). Since
$0<\alpha\leq\pi$ there are variations of $\rho(\mu)$ with non-zero
derivative of the trace function $I_{\mu}$.

Now  we deal with
  a surface of genus $g\geq 2$ in $\partial N$.
We take a decomposition of this surface in pairs of pants, so that
the curves of the decomposition are the meridians (each pair of
pants $P$ corresponds to a singular vertex and $\D P$ consists of
meridians). The deformations we require are easily constructed if
we prove that $\XX(P, SU(2))$ is locally parameterized by the trace
functions of the components of $\partial P$. Choose generators $a$
and $b$ of $\pi_1(P)$, so that $a$, $b$ and $ab$ represent the
three components of $\D P$. Since $a$ and $b$ generate a free
group, it is well known that
\[
(I_a,I_b,I_{ab}):\XX(P,SL_2(\mathbb C))\lra\mathbb C^3
\]
is an isomorphism of algebraic varieties defined over $\mathbb Q$,
see for instance \cite{GM}. This isomorphism implies that, for
each $\gamma\in\pi_1(P)$, $I_{\gamma}$ is a polynomial on
$I_a,I_b,I_{ab}$ with  coefficients in $\mathbb Q$. Thus
 \[
(I_a,I_b,I_{ab}):\XX(P,SL_2(\mathbb C))\cap\mathbb R^N\lra\mathbb
R^3
\]
is also an isomorphism. The irreducibility of
$\rho\vert_{\pi_1(P)}$ and Lemma~\ref{lem:SU(2)} imply that $
(I_a,I_b,I_{ab}):\XX(P,SU(2))\ra\mathbb R^3 $ is a local
diffeomorphism at the conjugacy class of the restriction
$\rho\vert_{\pi_1(P)}$.\end{proof}

This finishes the proof of Theorem \ref{thm:LR}. \end{proof}

\section{The fibration theorem}
\label{sec:Seifert}

Throughout this section, $\OO$ denotes a closed orientable small
3-orbifold $\OO$ and $X$ a cone \emph{metric} on $\OO$ of constant
curvature in $[-1,0)$ and with  cone angles  less than or equal to
the orbifold angles of $\OO$. Since we assume that $X$ is a cone
metric on $\OO$ and not only a cone structure, cf.\ Definition
\ref{dfn:conestruc}, the pairs
 $(\vert \OO\vert,\Sigma)$ and $(X,\Sigma_{X})$
are homeomorphic,
 and in particular $X$ is closed.
The main result of this section is:

\medskip

\begin{ftheorem}\textbf{\em{ (Thm.~\ref{thm:fibration})}.}
For $\om,D_0>0$ there exists $\de=\de(\om,D_0)>0$ such that: If
$X$ has $\om$-thick links (cf.\ Definition \ref{dfn:thicklinks}),
$\Diam(X) \geq D_0$ and if $X$ is $\delta$-thin, then $\OO$ is
Seifert fibred.
\end{ftheorem}

\medskip

\subsection{Local Euclidean structures}

The local geometry of thin cone manifolds
is modelled on noncompact Euclidean cone manifolds,
cf.\ \cite[part 2, Proposition 3.4]{CheegerGromov} in the case of manifolds.

Recall that by Corollary \ref{cor:soul} every noncompact Euclidean cone 3-manifold
$E^3$ with cone angles $\leq\pi$ has a \emph{soul} $S$.

\begin{lem}
\label{lem:LocalSoul} For every $\varepsilon>0$ and $R>1$, there exists
$\delta_0=\delta_0(\varepsilon,R,\omega)>0$ such that:  If $\de<\delta_0$,
$X$ is a cone manifold of curvature in $[-1,0)$ with $\omega$-thick links,
$\Diam(X)\geq D_0$ and $X$ is $\delta$-thin, then each $x\in X$ has a
neighborhood $U_x\subset X$, and a $(1+\varepsilon)$-bilipschitz
homeomorphism
\[
f\!: U_x\to\NN_{\nu_x}(S)
\]
where $\NN_{\nu_x}(S)$ is the normal cone fiber bundle, of radius
$\nu_x\in(0,1)$ depending on $x$, of the soul $S$ of a noncompact
Euclidean cone 3-manifold. In addition $\dim S=1$ or $2$, and
\[
 \max\big(d(f(x),S),\Diam(S)\big)\leq \nu_x/R.
\]
\end{lem}

\begin{proof}
Assume that the assertion were false. Then there exist
$\varepsilon>0$, $R>1$ and a sequence of cone manifolds $X_n$ with
diameter $\geq D_0$, curvature in $[-1,0)$ and $\om$-thick links
such that $X_n$ is $\frac1n$-thin, and there exist points $x_n\in
X_n$ for which the conclusion of the lemma does not hold.

The fact that $X_n$ is $\frac1n$-thin and has $\om$-thick links
implies that also the radii of embedded {\em singular} standard
balls in $X_n$ are $\leq r_n\ra0$. Let $\lambda_n>0$ be the
supremum of all radii $r$ such that $B_r(x_n)$ is contained in a
(smooth or singular) standard ball. We have that $\lambda_n\ra0$.
The sequence of rescaled cone manifolds
$(\frac1{\lambda_n}X_n,x_n)$ with base points subconverges to a
limit space $(E,x_{\infty})$. Observe that the balls
$B_1(x_n)\subset\frac1{\lambda_n}X_n$ are uniformly thick, as
shown by the following sublemma.


\begin{sublem}
Assume that the cone $3$-manifold $X$ has curvature $k\in[-\kappa,\kappa]$,
cone angles $\leq\pi$ and $\om$-thick links.
Suppose that the distance ball $B_r(p)\subset X$ is contained in a standard
ball.
Then $B_r(p)$ is $\de$-thick
with $\de=\de(\kappa,\om,r)>0$.
\end{sublem}
\begin{proof}
We may assume that $X$ is a complete singular cone. The smooth points in
$X$ with injectivity radius $<\de$ are contained in a tubular neighborhood
of radius $\rho(\om,\de)$ around the singularity where
$\lim_{\de\ra0}\rho(\om,\de)=0$, because $X$ has $\omega$-thick links. We
choose $\de<\frac{r}2$ sufficiently small such that
$\rho(\om,\de)<\frac{r}2$.
\end{proof}
It follows with the compactness theorem
(Corollary~\ref{thm:compactness}) that $E$ is a $3$-dimensional
Euclidean cone manifold.

Since $E$ is not compact, we can apply Theorem \ref{thm:euclid}
which gives the classification of noncompact Euclidean cone
$3$-manifolds with cone angles $\leq\pi$. The space $E$ is not a
complete cone because balls around $x_{\infty}$ with radii $>1$
are not contained in a standard ball. Hence the soul $S$ of $E$
has dimension 1 or 2. Let $N\subset E$ be a tubular neighborhood
around $S$ with radius $\rho>R\cdot diam(S\cup\{x_{\infty}\})$. We
use now that the convergence is bilipschitz. For sufficiency large
$n$, there exists a $(1+\eps)$-bilipschitz embedding
$(N,x_{\infty})\hookrightarrow (\frac1{\lambda_n}X_n,x_n)$, and
hence $(\lambda_n N,x_{\infty})\hookrightarrow (X_n,x_n)$. Hence
$x_n$ satisfies the conclusion of the lemma with
$\nu_{x_n}=\lambda_n\rho\ra0$, a contradiction.
\end{proof}

We apply this lemma to each point of $X$ with some constants
$R>1$, $\varepsilon>0$ to be specified later. Consider the
thickening
\[ W_x:=f^{-1}(\overline{\NN_{\lambda\nu_{x}}(S)}) \]
of the soul of $U_x$ where $0<\lambda<\frac1R$. We will also view
$W_x$ as a suborbifold of $\OO$.

The local models $E$ have 1- or 2-dimensional soul and therefore belong to
the following list by Theorem~\ref{thm:euclid}:
\begin{itemize}
    \item[--] When $S$ is 2-dimensional and orientable, then $E$ is
isometric to
    the product $S\times \mathbb R$. The possible Euclidean cone surfaces $S$
    are a torus $T^2$, a pillow $S^2(\pi,\pi,\pi,\pi)$, i.e.\ a Euclidean
    surface homeomorphic to $S^2$ with four cone points of angle $\pi$, and
    a Euclidean turnover $S^2(\alpha,\beta,\gamma)$ with cone angles
    $\al,\beta$ and $\ga$ satisfying $\alpha+\beta+\gamma=2\pi$.
    \item[--] When $S$ is 2-dimensional but non-orientable (possibly with
    mirror boundary), then $E=\tilde S\times\mathbb R/\iota $, where
    $\tilde S$ is the orientable double covering of $S$ and $\iota $ is an
    involution that preserves the product structure and reverses the
    orientation of each factor. Hence $E$ is a twisted line bundle over $
    S$.
    \item[--] When $\dim(S)=1$, then either $S=S^1$ or $S$ is an interval with
    mirror boundary (a quotient of $S^1$). In the former case, $E$ is
    either a solid torus or a singular solid torus. In the latter, $E$ is
    either a solid pillow or a singular solid pillow.
\end{itemize}

Not all of these possibilities can occur:
$W_x$ contains no turnover, because $\OO$ is small
and singular vertices of $X$ have thick links (i.e. the addition of
cone angles of edges adjacent to a singular vertex is
$>2\pi+\omega'$).

\begin{lem}
If $\eps=\eps(\om)>0$ is small enough,
then the soul $S$ of the local model for $W_x$ is neither a turnover
nor the quotient of a turnover.
\end{lem}
\begin{proof}
 Assume the contrary. By hypothesis, $\OO$ is closed and small.
Hence every turnover in $\OO$ bounds a discal suborbifold. Using
the bilipschitz homeomorphism from $W_x$ to the local model, it
follows that the sum of cone angles of the turnover is close to
$2\pi$. This contradicts the $\om$-thickness of links.
\end{proof}

By looking at the remaining possibilities,
we deduce:

\begin{cor} \label{lem:WSF}
Each $W_x$ admits a Seifert fibration (in the orbifold sense).
In particular, $\D W_x$ is a  union of smooth tori and pillows.
\end{cor}

\begin{lem}
\label{lem:toriO0}
(i) If $W_x$ contains a singular vertex,
then $\OO-int(W_x)$ is Haken.

(ii) Assume that $\OO$ is cyclic. If $\eps=\eps(\om)>0$ is small
enough, $R=R(\om)$ is large enough and if $W_x\cap
\Sigma\neq\emptyset$, then $\OO-int(W_x)$ is Haken.
\end{lem}

\begin{proof}
The main point is to prove that $\OO-int(W_x)$ is irreducible.
For this we have to show that $W_x$ is not contained in a discal suborbifold.
Since the pairs $(\OO,W_x)$ and $(\OO,\bar U_x)$ are homeomorphic,
this amounts to showing that $U_x$ is not contained in a discal suborbifold,
and we can use the metric properties of $U_x$.

If $U_x$ contains a singular vertex, this vertex and at least one singular
edge lie in the soul. Hence $U_x$ contains an entire singular edge or loop
and therefore cannot be included in a discal suborbifold. This proves
irreducibility in case (i).

Now we proceed with case (ii). Suppose now that $\OO$ is of cyclic type,
$U_x$ meets the singular locus and is contained in the discal suborbifold
$\De$. Topologically, $\De$ is a singular ball with one axis $a$. Hence
$U_x$ cannot contain an entire singular edge. By looking at the possible
local models it follows that $U_x$ contains at least two singular segments
of length $>\nu_x$ whose midpoints $m_1$ and $m_2$ have distance
$<\frac{\nu_x}R(1+\eps)$. By developing the smooth part of $X$ into model
space, and composing the developing map with the projection onto the axis
fixed by the holonomy representation, we find a 1-Lipschitz function on
$X$ whose restriction to $a$ is linear with slope 1. It follows that $a$
is distance minimizing inside $\De$ and hence $d_{\De}(m_1,m_2)>\nu_x$, a
contradiction. This finishes the proof of irreducibility in case (ii).

Since $\OO-int(W_x)$ has boundary, all that remains to check for
the Haken property is that there are no Euclidean or hyperbolic
turnovers. This follows from the smallness of $\OO$ in both cases
(i) and (ii).
\end{proof}

\subsection{Covering by virtually abelian subsets}

In this section,
we study general properties of coverings by virtually abelian subsets.
The arguments in this section closely follow \cite{Gro}
and \cite{BoP}.

We assign a special role to one of the subsets $W_x$ along which
we will cut $\OO$ later on. Namely, we choose $x_0\in X$ as
follows: If $\OO$ has singular vertices, we require that $W_{x_0}$
contains a singular vertex and that its radius $\nu_{x_0}$ is
almost maximal:
\[
\nu_{x_0}\geq \frac1{1+\varepsilon}\sup\{\nu_x\vert
W_x\cap\Sigma^{(0)}\neq\emptyset\}.
\]
If $\OO$ is cyclic, we make an analogous choice for  $W_{x_0}$
among all  $W_{x}$ that intersect the singular locus.
We denote $W_0=W_{x_0}$, $\OO_0=\OO-int(W_0)$, $\nu_0=\nu_{x_0}$.
In view of Lemma~\ref{lem:toriO0}, $\OO_0$ is Haken.

\begin{dfn}We say that a subset $S\subset \OO$ is \emph{virtually abelian in
$\OO_0$} if,
for each connected component $Z$ of $S\cap\OO_0$,
the image of
$\pi_1(Z)\to\pi_1(\OO_0)$
in the fundamental group of the corresponding component of $\OO_0$
is virtually abelian. Moreover,
for $x\in X$ we define:
\[
va(x) =\sup\{ r>0\mid B_r(x) \text{ is virtually abelian in
}\OO_0\}
\]
and
\[ r(x)=\inf(\frac{va(x) }8,1) .\]
\end{dfn}

\begin{lem}\label{lem:cov} Let $x,y\in X$. If $B_{r(x)}(x)\cap
B_{r(y)}(y)\neq\emptyset$, then
\begin{itemize}
\item[(a)]  $3/4\leq{r(x)}/{r(y)}\leq 4/3$;
\item[(b)]  $B_{r(x)}(x)\subset B_{4r(y)}(y)$.
\end{itemize}
\end{lem}
\begin{proof}
We may assume that $r(x)\leq r(y)$,
and moreover $r(x)=\frac18va(x)<1$.
By the triangle inequality we have
\[ va(x) \geq va(y)-r(x)-r(y) ,\]
so
\[ 8r(x)=va(x)\geq 8r(y)-r(x)-r(y) \geq 6r(y).\]
This shows part (a).
Part (b) follows because $2r(x)+r(y)<4r(y)$.
\end{proof}

\begin{lem}
\label{lem:covc} For $R$ sufficiently large it holds $ W_0\subset
B_{\frac{r(x_0)}9}(x_0)$.
\end{lem}
\begin{proof}
This follows because
$va(x_0)\geq\frac1{1+\eps}\nu_{x_0}(1-\frac1R)$ and $W_0$ is
contained in the ball of radius $3(1+\eps)\frac{\nu_x}R$ around
$x_0$.
\end{proof}

We proceed to construct coverings of $X$.
We have already distinguished a point $x_0\in W_0$.
Consider sequences $\{x_0,x_1,\ldots\}$ starting with $x_0$,
such that:
\begin{equation}\label{eqn:property}
 \text{ the balls }
B_{\frac14 {r(x_0)}}(x_0),B_{\frac14 {r(x_1)}}(x_1),\ldots \text{ are
pairwise disjoint.}
\end{equation}
A sequence satisfying (\ref{eqn:property}) is finite, by
Lemma~\ref{lem:cov} and compactness of $X$.

\begin{lem}\label{lem:seq} If the sequence $\{x_0,x_1,\ldots,x_p\}$ is
maximal for property (\ref{eqn:property}), then the balls
\(B_{\frac23 r(x_0)}(x_0),\ldots, B_{\frac23 r(x_p)}(x_p)\) cover
$X$.
\end{lem}
\begin{proof}
Let $x\in X$ be an arbitrary point. By maximality, there exists a
point $x_j$ such that $B_{\frac14 {r(x)}}(x)\cap B_{\frac14
{r(x_j)}}(x_j)\neq\emptyset$. By Lemma~\ref{lem:cov} we have
$r(x)\leq\frac43r(x_j)$ and $d(x,x_j)\leq\frac14(r(x)+r(x_j))\leq
\frac7{12}r(x_j)$. Thus $x\in B_{\frac23 r(x_j)}(x_j)$.
\end{proof}

We fix now a sequence $x_0,x_1,\ldots,x_p$ maximal for property
(\ref{eqn:property}) and consider
the covering of $X$ by the open sets
\begin{itemize}
\item $V_0=B_{r(x_0)}(x_0)$ and
\item  $V_i=B_{r(x_i)}(x_i)-W_0$ for  $i=1,\ldots,p$.
\end{itemize}

Lemmas \ref{lem:covc} and \ref{lem:seq} imply that the open
sets $V_0,\ldots,V_p$ cover $X$.
We denote $r_i:=r(x_i)$ and $B_i:=B_{r(x_i)}(x_i)$.

\begin{lem}
There is an a universal bound $N$ on the number of balls $B_i$
that can intersect a fixed ball $B_k$.
\end{lem}
\begin{proof}
 For every ball $B_i$ intersecting $B_k$ it holds
 $B_i\subset B_{2 r_i+r_k}(x_k)\subseteq B_{4 r_k}(x_k)$.
  On the other hand, the points $x_i$ are separated
 from each other, since
 $d(x_{i_1},x_{i_2})\geq\frac14(r_{i_1}+r_{i_2})\geq\frac38r_k$. Thus
 the number of such $x_i$ is bounded above by:
 \[
     \frac{\operatorname{vol}(B_{4
r_k}(x_k))}{\operatorname{vol}(B_{\frac3{16} r_k}(x_i))}
     \leq
     \frac{\operatorname{vol}(B_{8
r_k}(x_i))}{\operatorname{vol}(B_{\frac3{16} r_k}(x_i))}
     \leq
     \frac{v_\kappa(8 r_k)}{v_\kappa(\frac3{16} r_k)}.
 \]
Here $v_{\kappa}(r)$ denotes the volume of the ball of radius $r$
in the  3-space of constant curvature $\kappa\in[-1,0)$, and the
last inequality follows from Bishop-Gromov. Since $r_k\leq 1$, the
ratio $v_{\kappa}(8 r_k)/v_{\kappa}(\frac34 r_k)$ is bounded.
\end{proof}

In particular,
the dimension of our covering is universally bounded by $N$.
We want to decrease the dimension of the covering $\{V_0,\ldots,V_p\}$
while keeping the properties that the covering sets are virtually abelian
and only one of them meets $W_0$.

Using a partition of unity $(\phi_i)$ subordinate to $(V_i)$
one can construct a map
\begin{equation}
\label{tonerve} f=\frac1{\sum_i\phi_i}(\phi_0,\ldots,\phi_p)\! :X\ra
\Delta\subset \mathbb R^{p+1},
\end{equation}
where $\Delta$ is the unit simplex. The image of $f$ is contained
in the nerve of the covering $K\subset \Delta$, which has
dimension $\leq N$. We start by controlling the Lipschitz constant
of the map $f\!:X\to K$.

\begin{lem}
\label{lem:boundlip} There exists a constant $L_N>0$ such that the
partition of unity can be chosen so that the restriction $f\vert_{V_k}$ is
$\frac{L_N}{r_k}$-Lipschitz.
\end{lem}

We first need the following geometric property of our covering:

\begin{sublem}
Every $x\in X$ belongs to an open set $V_k$ such that
    $d(x,\partial V_k) \geq{r_k}/3$.
\end{sublem}
\begin{proof}
If $x\in B_{\frac23r_0}(x_0)$ we choose $k=0$.
Suppose that $x\not\in B_{\frac23r_0}(x_0)$.
Then there exists $k$ with $x\in B_{\frac23r_k}(x_k)$.
The assertion is trivial if $B_k$ and $B_0$ are disjoint.
We assume therefore also that $B_k\cap B_0\neq\emptyset$.
Then by Lemma~\ref{lem:covc}:
\[ d(x,W_0)\geq d(x,x_0)-\frac19r_0
\geq \frac23r_0-\frac19r_0 \geq \frac34\cdot\frac59r_k>\frac13r_k.
\]
Hence $d(x,\D V_k)\geq \frac13r_k$.
\end{proof}

\begin{proof}[Proof of Lemma~\ref{lem:boundlip}.] Let $\tau:[0,1]\ra[0,1]$
be an auxiliary  4-Lipschitz function which vanishes in a
neighborhood of $0$ and satisfies $\tau|_{[\frac13,1]}\equiv1$. We
put $\phi_k:=\tau(\frac1{r_k}d(\D V_k,\cdot))$ on $V_k$ and extend
it trivially to $X$. Then $\phi_k$ is $\frac4{r_k}$-Lipschitz.

Let $x\in V_k$. Then at most $N+1$ functions $\phi_i$ are non-zero
in $x$, and all of them have Lipschitz constant $\leq
\frac43\cdot\frac4{r_k}$. The claim follows since the functions
\[ (x_0,\dots,x_N)\mapsto \frac{x_k}{\sum_{i=0}^N x_i} \]
are Lipschitz on $\{x\in {\mathbb R}^{N+1} \mid x_i\geq0 \;\forall i\land
\sum_{i=0}^N x_i\geq1\}$.
\end{proof}

We now homotope $f$ into the $3$-skeleton $K^{(3)}$ by an inductive procedure
while controlling the local Lipschitz constant.

\begin{lem}
\label{lem:homotopelowerskel}
For $d\geq4$ and $L>0$ exists $L'=L'(d,L)>0$ such that the following is true:

Suppose that
$g:X\ra K^{(d)}$ is a continuous map which is
$\frac{L}{r_k}$-Lipschitz on $V_k$
and has the property that the inverse image of the open star of the vertex
$v_{V_k}\in K^{(0)}$ is contained in $V_k$.
Then $g$ can be homotoped to a map $\tilde g:X\ra K^{(d-1)}$
with the same properties,
$L$ being replaced by $L'$.
\end{lem}
\begin{proof}
It suffices to find a constant $\theta>0$
such that every $d$-dimensional simplex $\si\subset K$ contains a point
$z$ at distance $\geq\theta$ from both $\D\si$ and the image of $g$.
To push $g$ into the $(d-1)$-skeleton
we compose it on $\si$ with the central projection from $z$.
This will increase the Lipschitz constant by a factor bounded in terms of $d$,
and it reduces the inverse images of open stars of vertices.

If $\theta$ does not satisfy the desired property for some $d$-simplex $\si$,
then $image(g)\cap int(\si)$ must contain a subset
of at least $C(d)\cdot\frac1{\theta^d}$ points
with pairwise distances $\geq\theta$.
Let $A\subset X$ be a set of inverse images, one for each point.
Let $v_{V_k}$ be a vertex of $\si$.
Then $A\subset V_k\subseteq B_k$.
Since $f$ is $\frac{L}{r_k}$-Lipschitz continuous on $V_k$,
the points of $A$ are separated by distance $\frac1{L}r_k\theta$.
Since $r_k\leq1$,
volume comparison implies that $A$ contains at most
$C\cdot(\frac{L}{\theta})^3$ points.
The inequality $C(d)\cdot\frac1{\theta^d}\leq C\cdot(\frac{L}{\theta})^3$
yields a positive lower bound $\theta_0(d,L)$ for $\theta$.
Hence any constant $\theta<\theta_0$ has the desired property.
\end{proof}

\begin{lem}
\label{lem:volsubcubic}
For sufficiently small $\eps>0$ there exists a constant $C=C(\eps)>0$
such that
\[ \Vol(V_i)\leq C\frac1R r_i^3.\]
for all $i$.
\end{lem}

\begin{proof}
We first show that $W_0$ does not enter to far into the other
sets $U_{x_i}$.

\begin{sublem}
There exists a constant $c=c(\eps)>0$ such that, if $R>1$ is
sufficiently large, then $d(x_i,W_0)\geq c\, \nu_{x_i}$ for all
$i\neq0$.
\end{sublem}
\begin{proof}
By Lemma~\ref{lem:covc}, $W_0\subset B_{\frac{r_0}9}(x_0)$ and we
obtain:
\[ d(x_i,W_0)\geq d(x_i,x_0)-\frac19 r_0
\geq \frac14r_0-\frac19r_0 >\frac18r_0 \geq \frac1{64(1+\eps)}\nu_{x_0}.
\]
For the last estimate,
we use that $va(x_0)\geq\frac1{1+\eps}\nu_{x_0}$ and,
since $\nu_{x_0}\leq1$,
$r_0=\inf(\frac{va(x_0)}8,1)\geq\frac1{8(1+\eps)}\nu_{x_0}$.

We now assume that $W_0$ intersects $U_{x_i}$
because otherwise there is nothing to show.
If $W_0\subset U_{x_i}$
then, according to our choice of $W_0$, we can compare the radii
$\nu_{x_0}$ and $\nu_{x_i}$ by
$\nu_{x_0}\geq\frac1{1+\eps}\nu_{x_i}$,
and the assertion holds with $c<(8(1+\eps))^{-2}$.

We are left with the case
that $W_0\not\subset  U_{x_i}$
but intersects the ball of radius, say, $\frac{\nu_{x_i}}4$ around $x_i$.
Then we can bound the ratio $\frac{diam(W_0)}{\nu_{x_i}}$ from below by:
\[ (1+\eps)\frac{\nu_{x_i}}R+\frac{\nu_{x_i}}4+diam(W_0)
\geq\frac{\nu_{x_i}}{1+\eps}. \] By definition of $W_0$ we have
$diam(W_0)\leq(1+\eps)\frac2R\nu_{x_0}$. Combining these estimates, we
obtain a lower bound for $\frac{d(x_i,W_0)}{\nu_{x_i}}$, as claimed.
\end{proof}

The sublemma implies that $r_i\geq c\,\nu_{x_i}$. By Bishop-Gromov
inequality,
\[
 \Vol(V_i)\leq
 \Vol\big(B_{r_i}(x_i)\big)\leq
\Vol\big(B_{c\,\nu_{x_i}}(x_i)\big)
\frac{\V_{\kappa}(r_i)}{\V_{\kappa}(c\,\nu_{x_i})} \leq
\Vol\big(B_{c\,\nu_{x_i}}(x_i)\big)\,c_1\,
\frac{r_i^3}{\nu_{x_i}^3}
\]
for some uniform $c_1>0$, where $\V_{\kappa}(r)$ denotes the
volume of the ball of radius $r$ in the space of constant
curvature $\kappa\in [-1,0)$. Using the geometry of the local
models, it follows that
\[
\Vol\!\big(B_{c\,\nu_{x_i}}(x_i)\big)\leq  \Vol (U_{x_i})\leq
 c_2\frac{\nu_{x_i}^3}R
\]
for some $c_2>0$. Thus
 \(
 \Vol(V_i)\leq c_1\,c_2\,\frac{r_i^3}R.
 \)
\end{proof}

Now we can further homotope $f$ into the $2$-skeleton.

\begin{prop}
For suitable constants $\eps>0$ and $R>1$, the map $f$ in
(\ref{tonerve}) is homotopic to a map
\[ \tilde f:X\ra K^{(2)} \]
with the property that
the inverse image of the open star of the vertex
$v_{V_k}\in K^{(0)}$ is contained in $V_k$.
\end{prop}
\begin{proof}
The inverse image under $f$ of the open star of $v_{V_k}$
is contained in $V_k$.
Using Lemma~\ref{lem:homotopelowerskel} repeatedly,
we can homotope $f$ to a map
$\hat f:X\ra K^{(3)}$
which is locally Lipschitz
and satisfies $\hat f^{-1}(star(v_{V_k}))\subset V_k$.
More precisely,
there is a universal constant $\hat L$
such that $\hat f$ is $\frac{\hat L}{r_k}$-Lipschitz on $V_k$.

It suffices to show that no $3$-simplex $\si\subset K$
is contained in the image of $\hat f$.
The inverse image
$\hat f^{-1}(int(\si))$ lies in the intersection of sets $V_j$
where $v_{V_j}$ runs through the vertices of $\si$.
Let $V_k$ be one of these.
With Lemma~\ref{lem:volsubcubic} it follows that
\[
\Vol(image(\hat f)\cap\si)\leq \Vol(\hat f(V_k))
 \leq ( \frac{\hat L}{r_k})^3 \Vol(V_k)\leq C \hat L^3\frac1R
\]
with uniform constants $C$ and $\hat L$. So, if $R$ is large
enough, $\Vol(image(\hat f)\cap\si)<\Vol(\si)$.
\end{proof}

Note that $f$ maps $W_0$ to the vertex $v_{V_0}$
because $W_0$ intersects none of the sets $V_j$ with $j\neq0$.
The proposition therefore implies the following properties
which will be crucial below.

\begin{itemize}
\item[(i)] $\tilde f(W_0)=\{v_{V_0}\}$.
\item[(ii)]
For every vertex $v$ of $K$,
$\tilde f^{-1}(\Star v)$ is virtually abelian in $\OO_0$.
\end{itemize}

\subsection{Vanishing of simplicial volume}\label{subsection:simplicial}

The orbifold $\OO_0$ is Haken and therefore has a JSJ-splitting into
Seifert and hyperbolic suborbifolds.

\begin{prop}
All components in the JSJ splitting of $\OO_0$ are Seifert.
\end{prop}
\begin{proof}
Since the orbifold $\OO_0$ is Haken, it is very good
and there is a finite covering
\[ p\! : N\to\OO_0 \]
by a manifold \cite{McCMi}. The boundary $\partial N$ is a union
of tori. The JSJ splitting of $\OO_0$ pulls back to the JSJ
splitting of $N$. We have to show that no hyperbolic components
occur in the JSJ splitting of $N$.

We may assume that the boundary of $N$ is incompressible because otherwise
$N$ is a solid torus and the assertion holds. We construct a closed
manifold $\bar N$ by Dehn filling on $N$ as follows. Let $Y\subset N$ be a
component of the JSJ splitting which meets the boundary, $Y\cap\D
N\neq\emptyset$. If $Y$ is hyperbolic we choose, using the hyperbolic Dehn
filling theorem, the Dehn fillings at the tori of $Y\cap\D N$ such that
the resulting manifold $\bar Y$ remains hyperbolic. If $Y$ is Seifert, we
fill in such a way that $\bar Y$ is Seifert and the components of $\D Y-\D
N$ remain incompressible. This can be done because the base of the Seifert
fibration of $Y$ is neither an annulus nor a disc with zero or one cone
point. The manifold $\bar N$ has a JSJ splitting along the same tori as
$N$ and with the same number of hyperbolic (and also Seifert) components.

It suffices to show that $\bar N$ has zero simplicial volume,
because then \cite[Sec. 3.5]{Gro} and \cite{thilo} imply
that $\bar N$ contains no hyperbolic component in its JSJ splitting.
To this purpose we will apply Gromov's vanishing theorem,
see \cite[Sec. 3.1]{Gro}, \cite{Iva}.

We compose $\tilde f$ with the projection $p$
and extend the resulting map $N\ra K^{(2)}$ to a map
\[ h\!:\bar N\to K^{(2)} \]
by sending the filling solid tori to the vertex $v_{V_0}$. Note
that $h$ is continuous because $\tilde f(\D\OO_0)=\{v_{V_0}\}$.
The inverse images under $h$ of open stars of vertices are
virtually abelian as subsets of $\bar N$, because they are already
virtually abelian in $N$ and the filling tori intersect only one
of the subsets. These subsets yield an open covering of $\bar N$
with covering dimension $\leq2$. By Gromov's theorem, the
simplicial volume of $\bar N$ vanishes.
\end{proof}

{\em Conclusion of the proof of Theorem~\ref{thm:fibration}.}
Since $\OO_0$ is graphed and $\OO$ results from $\OO_0$ by gluing
in a Seifert orbifold it follows that $\OO$ is graphed. Since
$\OO$ is moreover atoroidal, it must be Seifert. \qed

Laboratoire \'Emile Picard, CNRS UMR 5580, Universit\'e Paul
Sabatier, 118 Route de Narbonne, F-31062 Toulouse Cedex 4, France,
boileau@picard.ups-tlse.fr

Mathematisches Institut, Universit\"at M\"unchen, Theresienstra\ss
e 39, D-80333 M\"un\-chen, Germany,
leeb@mathematik.uni-muenchen.de

Departament de Matem\`atiques, Universitat Aut\`onoma de
Barcelona, E-08193 Be\-lla\-ter\-ra, Spain, porti@mat.uab.es

\end{document}